\documentclass[11pt,oneside,reqno]{amsart}

\usepackage{epsfig,psfrag,subfigure}
\usepackage{epstopdf}
\usepackage{amsmath,amssymb,amsfonts,latexsym, mathrsfs,euscript}
\usepackage{graphicx,graphics,cite}
\usepackage{wrapfig}
\usepackage[usenames]{color}
\usepackage{geometry}                % See geometry.pdf to learn the layout options. There are lots.
\usepackage{colortbl,booktabs,ctable,multirow}
\usepackage{epstopdf}
\usepackage{subfig}

\DeclareGraphicsExtensions{.pdf,.png,.jpg}

\newcommand{\beq}{\begin{equation}}
\newcommand{\eeq}{\end{equation}}
\newcommand{\beqa}{\begin{eqnarray}}

\newcommand{\eeqa}{\end{eqnarray}}

\newcounter{nfig}

\newcommand{\ignore}[1]{}

\newcommand\ds{\displaystyle}

\newcommand\bs{\boldsymbol}

\newcommand{\sech}{{\rm sech}}
\newcommand{\sgn}{{\rm sgn}}

\newcommand{\ti}{{\tilde{t}}}

\newcommand{\uti}{{\tilde{\bs{u}}}}
\newcommand{\mti}{{\tilde{\bs{m}}}}
\newcommand{\nati}{{\tilde{\nabla}}}

\newtheorem{note}{Note}
\title[Wave dynamics for Euler-Poincar\'e equations]{Solitary waves and $N$-particle algorithms for a class of Euler-Poincar\'e equations}

\author{Roberto Camassa}
\address{Department of Mathematics, University of North Carolina, Chapel Hill, 27599, USA}
\email{camassa@amath.unc.edu\\ TEL:919-962-8476}

\author{Dongyang Kuang}
\address{Department of Mathematics, University of Wyoming, Laramie, WY 82071-3036, USA}
\email{dkuang@uwyo.edu\\ TEL:307-766-4221}

\author{Long Lee}
\address{Department of Mathematics, University of Wyoming, Laramie, WY 82071-3036, USA}
\email{llee@uwyo.edu\\ TEL:307-766-4368}
\date{\today}

\begin{document}

\maketitle

\begin{abstract}
We study a class of partial differential equations (PDEs) in the family of the 
so-called Euler-Poincar\'e differential systems, with the aim of developing a foundation for numerical algorithms of their solutions. This requires  
particular attention to the mathematical properties of this system when the associated  class of elliptic operators possesses non-smooth kernels. By casting the system in its Lagrangian (or characteristics) form, we first formulate 
a particles system algorithm in free space with homogeneous Dirichlet boundary conditions for the evolving fields. We next examine the deformation of the system when non-homogeneous ``constant stream" boundary conditions are assumed. We show how this simple change at the boundary deeply affects the nature of the evolution,  from hyperbolic-like to dispersive with a non-trivial dispersion relation, and examine the potentially regularizing properties of singular kernels offered by this deformation. 
From the particle algorithm viewpoint, kernel singularities affect the existence and uniqueness of solutions to the corresponding  ordinary differential equations systems.   We illustrate this with the case when the operator kernel assumes a conical shape over the spatial variables, and examine in detail two-particle dynamics under the resulting lack of Lipschitz-continuity. Curiously, we find that for the conically-shaped kernels the motion of the related two-dimensional waves can become completely integrable under appropriate initial data. This reduction projects the two-dimensional system to the one-dimensional completely integrable Shallow-Water equation [Camassa, R. and Holm, D. D., Phys. Rev. Lett., 71, 1961-1964, 1993], while retaining the full dependence on two spatial dimensions for the single channel solutions. 
Finally, by comparing with an operator-splitting  
pseudospectral method we illustrate the performance of the particle algorithms with respect to their Eulerian counterpart for this class of non-smooth kernels.

\end{abstract}

\begin{description}
\item[{\footnotesize\bf keywords}]
{\footnotesize Euler-Poincar\'e differential equations, diffeomorphisms, Lagrangian formulation, dispersive, particle algorithms, completely integrable, Shallow-Water equation}
\end{description}

\section{Introduction}

The Euler-Poincar\'e differential equations, also called the Euler equations for planar diffeomorphisms, originate in models of template matching and are of general interest as evolution equations on Riemannian manifolds endowed with Sobolev metrics \cite{bib:bmty05, NP, bib:hrty04, bib:my01, bib:jm00}. In one spatial dimension, the system of the partial differential  equations (PDEs) we consider in this paper may reduce to a completely integrable equation arising as a model of long wave evolution in shallow water, derived in~\cite{bib:ch93, chh94} (hereafter referred to as the SW -- for Shallow-Water -- equation). 
In this physical context, these PDEs can be used as  a model of the competition between nonlinear and dispersive effects, whose intertwined properties contribute to the rich dynamics exhibited by this class of nonlinear evolution equations.

One notable feature of the model PDEs under study (also known as the ``EPDiff'' differential equations in some literature) is that they admit traveling-wave {\it weak}-solutions, for which the momentum-like variable may be viewed as concentrated at a single point as if it were a ``particle.'' In fact, these particles are reminiscent of point vortices in Euler equations, which are widely studied  in the literature both for their inherent interest as dynamical systems and as a foundation for numerical algorithms for the evolution of general Euler solutions. Similarly to this latter case, once projected onto the particle solution class  the evolution of the PDEs can be written in the form of a finite-dimensional particle system of ordinary different equations (ODEs).
%similar to the point vortices of the Euler equations governing the dynamics of an incompressible fluid. 
We will refer to this system of ODEs  as the { $N$-particle finite-dimensional dynamical system}, or {\it $N$-particle system}.

For nonlinear dispersive equations, the interplay between nonlinearity and dispersion is often understood as the mechanism underlying the existence of traveling wave solutions. However, the way in which solitary waves emerge and can become the dominant structure in the long time evolution out of generic initial conditions can take various forms, depending on the structure of the equations, especially in {multiple space  dimensions}. Study of the $N$-particle system for a class of the model PDEs, where the interplay between dispersion and nonlinearity is varied continuously within a one-parameter family, is a convenient way to shed light on this as well as to investigate the role played by traveling waves in the long time evolution from a range of initial data. 

An interesting application for the $N$-particle system is template matching. This is commonly used in problems of image reconstruction and pattern recognition \cite{image1, image2}. Template matching can be formulated as an variational problem, such as finding the shortest or least expensive path of continuous deformation of one geometric object (reference template) into another one (target template). In this context, the time-dependent deformation process produces geodesic evolution equations which falls into the Euler-Poincar\'e theory\cite{bib:hrty04}. 
A practical application for template matching is computational anatomy (CA) \cite{anatomy},
whereby a medical image can be discretized into a set of so-called {\it landmark points}, which in turn can be represented by the $N$-particle system of the model PDEs. The template matching problem, in terms of landmark points, becomes the landmark-matching problem \cite{bib:mumford, bib:hrty04}.  
While the template matching problem is related to the issue of comparing two geometric objects, and thus more concerned with a variational boundary-value problem, the initial-value problem associated to the integration of the model equations and/or their $N$-particle system has important consequences for applications, especially for designing numerical matching procedures \cite{bib:kl14, bib:hrty04, bib:mty06}. As noted above, the $N$-particle algorithms and dynamics play an important role in both the model PDEs and their applications. However, despite some notable efforts \cite{NP, bib:cdm12}, there are aspects of the $N$-particle systems and their dynamics that have not been thoroughly investigated, particularly  when certain smoothness properties are not satisfied. The aim of this paper is to examine some of these aspects, with the brooder goal 
of establishing the foundations of potentially efficient numerical algorithms for the solution of this class of  model PDEs.  

The steps we take towards implementing this goal are as follows. We first introduce the Lagrangian formulation of for the class of PDEs under investigation, which allows us to discretize the resulting integral-differential equations to obtain the $N$-particle systems for the model PDEs. Our approach introduces a mesh size (e.g., $dx\, dy$  in two dimensions) naturally and explicitly, a necessary step for proving the convergence of the particle algorithm. The singular nature of some the particle solutions suggests that a form of regularization might be needed in order to implement numerical algorithms. We examine a possible class of regularizations of the model PDE, and show that this follows simply from assuming non-zero 
constant boundary conditions on the evolving fields. The deformation leads to non-trivial dispersive evolution, and the corresponding dispersion relation explicitly displays  the limitations that this can present when used as regularization for non-smooth solutions (unlike its one-dimensional counterparts, see e.g.~\cite{MZM06}).   We illustrate two-particle dynamics for non-Lipschitz kernels (with particular attention to the example where the power of the associated elliptic operator is equal to $3/2$) via direct numerical simulations and analysis. We analyze the scattering properties under the loss of uniqueness of ODE solutions due to these non-Lipschitz kernels. We also show that when the motion of these particles is confined to a straight line, the dynamics of the associated solitary waves (dubbed as ``conons") coincides with that of the SW equation and is therefore completely integrable, even though the ``single channel" solution retains its dependence on two spatial dimensions.  Finally, we demonstrate that the $N$-particle system can be advantageous for solving the model PDEs with non-smooth solutions, and is also robust enough to capture more regular solutions, by comparing with an operator-splitting pseudospectral method for solving the Eulerian form of the model equation.

\section{Equations of motion}

By using index notation with Einstein convention on sums over repeated indexes for 
the (column) vectors $\bs{m}\equiv \left\{ m_\alpha \right\}_{\alpha=1}^n $ and $\bs{u}\equiv \left\{ u_\alpha \right\}_{\alpha=1}^n $, 
the system of equations can be written as
\begin{equation}\label{eq:EPDIFFi}
\partial_t m_\alpha + u_\beta\partial_\beta m_\alpha+ m_\beta \partial_\alpha u_\beta + m_\alpha \partial_\beta u_\beta=0 \, , 
%\bs{m}_t+ \bs{u}\cdot\nabla\bs{m} + \nabla\bs{u}^{T}\cdot\bs{m}+\bs{m}(\nabla\cdot\bs{u}) = 0,
\end{equation}
or, in short-hand vector notation,
\begin{equation}\label{eq:EPDIFF}
\bs{m}_t  + (\bs{u}\cdot\nabla)\bs{m} +\bs{m} \cdot(\nabla\bs{u})^{T}+\bs{m}(\nabla\cdot\bs{u}) = 0,
\end{equation}
with $t\in\mathbb{R}^+$, $\bs{x}, \bs{u}$ and $\bs{m}\in\mathbb{R}^{n}$, and spatial partial derivatives are labeled by coordinate index. For ease of notation, here and throughout the rest of the paper we will use Greek alphabet indexes to label coordinates, to distinguish them from particle labels (see below) in Latin alphabet, and suppress explicit argument dependences in the notation unless this becomes necessary to avoid confusion. The field $\bs{u}$ and its associated momentum-like variable $\bs{m}$ are formally related by an elliptic  
 operator $\mathcal{L}$
\begin{equation}\label{eq:Yukawa}
\bs{m}=\mathcal{L}\, \bs{u} \, .  
\end{equation}
With boundary vanishing boundary conditions at infinity, 
%$\bs{u} \to 0$ as $\bs{x} \to \infty$ sufficiently fast, 
the operator $\mathcal{L}$ is assumed to 
be invertible, with its inverse being explicitly written in terms of the corresponding Green function $\bs{G}$, so that $\bs{u}$ can also be represented by the convolution
%by convolution of the Green's function of the Yukawa operator $\mathcal{L}$ with the momentum $\bs{m}$. %Let $\bs{G}$ be the Green's function of $\mathcal{L}$. Then
\begin{equation}\label{eq:convol}
\bs{u} =\bs{G}*\bs{m}.
\end{equation}
%\vspace{0.5cm}
%\subsection{Green functions}
In this paper we will restrict our attention to the particular choice of $\mathcal{L}\equiv\mathcal{L}^{ b}$ as the (Yukawa) operator defined by
\begin{equation}\label{eq:elliptic}
\mathcal{L}^{ b}=(\bs{I}-a^2\nabla^2)^{b},
\end{equation}
for  $b >0$. 
%The operator $D$  in the matrix form is
%\beq D\bs{u}=\left( \begin{array}{cccc}
%\partial_{x_1} u_1 & \partial_{x_2} u_1 &\cdots & \partial_{x_n} u_1\\
%\partial_{x_1} u_2 & \partial_{x_2} u_2  &\cdots & \partial_{x_n} u_2 \\
%\vdots & \vdots & \ddots &\vdots\\
%\partial_{x_1} u_n & \partial_{x_2} u_n &\cdots & \partial_{x_n} u_n
%\end{array} \right).\eeq
%$(D\bs{u})^T$ is the transpose of $D\bs{u}$. 
Further, for the domain of $\mathcal{L}^{ b}$ we will take the Schwartz space of rapidly 
decaying functions in $\mathbb{R}^{n}$. 
For any $b>0$, including non-integer values, equation (\ref{eq:Yukawa}) can be defined in Fourier space, 
\beq\label{eq:m-Fourier}
\bs{\hat{u}} = (\hat{\mathcal{L}^{ b}})^{-1}\bs{\hat{m}},\quad\text{where}\,\,\,(\hat{\mathcal{L}^{ b}})^{-1} = \frac{1}{(1+a^2 |\bs{k}|^2)^{ b}},\quad |\bs{k}|=\sqrt{k_1^2+k_2^2\cdots+k_n^2},
\eeq
where $k_\alpha$ is the $\alpha^{th}$ wavenumber. Since $\mathcal{L}^{ b}$ is rotationally invariant and diagonal, then $\bs{G}(\bs{x}) = G_{ b-n/2}(|\bs{x}|)\bs{I}$ for a scalar function $G_{ b-n/2}$, with $|\bs{x}|=\sqrt{x_1^2+x_2^2+\cdot+x_n^2}$.   The scalar Green function $G_{ b-n/2}$ can be obtained by a combination of Bessel and Gamma functions,  
%For $\bs{x}\in\mathbb{R}^{n}$, the Green's function of $\mathcal{L}=(\bs{I}-a^2\nabla^2)^{ b}$ is 
\beq\label{eq:Green-nD}
G_{ b-n/2}(|\bs{x}|) = \frac{2^{n/2- b}}{(2\pi a)^{n/2}a^{ b}\Gamma( b)}|\bs{x}|^{ b -n/2}K_{ b -n/2}\left(\frac{|\bs{x}|}{a}\right),
\eeq
where $K_{ b -n/2}$ is the modified Bessel function of the second kind of order $ b -n/2$ and $\Gamma( b)$ is the usual notation for the Gamma function~\cite{bib:mumford}.

\subsection{Lagrangian formulation}
Equation~(\ref{eq:EPDIFFi}) is the Eulerian counterpart of a
Lagrangian formulation obtained from the characteristics 
$\bs{x}=\bs{q}(\bs{\xi},t)$
\begin{equation}
{d \bs{q}\over d t}\equiv 
\bs{u}(\bs{q}(\bs{\xi},t))\, , \qquad \bs{q}(\bs{\xi},0)=\bs{\xi}\, , 
\label{eq:chars}
\end{equation}
by defining the conjugate field $\bs{p}(\bs{\xi},t)$
\begin{equation}
\bs{m}\big(\bs{q}(\bs{\xi},t),t\big) \equiv \,{\bs{p}(\bs{\xi},t)\over  J(\bs{\xi},t) }  \, , 
\label{eq:charsp}
\end{equation}
where $J(\bs{\xi},t)$ is the Jacobian determinant of the diffeomorphism $\bs{x}=\bs{q}(\bs{\xi},t)$ parametrized by time $t$, 
$$
J(\bs{\xi},t) \equiv \det \left( {\partial x_i \over \partial \xi_j }\right) \, , 
$$
with $J(\bs{\xi},0) =1$. For as long as $J(\bs{\xi},t) \neq 0$ the definition~(\ref{eq:charsp}) is well posed, and the evolution equation preserves 
the smoothness of the initial data. Thus, from the characteristic formulation of equation~(\ref{eq:EPDIFFi}), local well posedness and existence of solutions can 
be readily established. 
The well known property of determinant differentiation
\begin{equation}
{d J \over d t} = J \,\,  \nabla\cdot\bs{u} 
\label{eq:det_t}
\end{equation}
shows  that the $\bs{m}$ evolution in equation~(\ref{eq:EPDIFFi}), with 
our choice of symmetric Green functions,
is defined by 
%\begin{equation}
%{d \bs{p} \over d t} = -(\nabla \bs{G}* {\bs p}) \cdot \bs{p}
%\label{eq:charspdot}
%\end{equation}
%or, with index notation,
%(with index notation) 
\begin{equation}
{d \bs{p} \over d t} = - \int_{\mathbb{R}^{n}} 
G'_{ b-n/2}\big(|\bs{q}(\bs{\xi},t) -\bs{q}(\bs{\eta},t)|\big) 
{\bs{q}(\bs{\xi},t) -\bs{q}(\bs{\eta},t) \over |\bs{q}(\bs{\xi},t) -\bs{q}(\bs{\eta},t)|}
\, \, \bs{p}(\bs{\xi},t) \cdot \bs{p} (\bs{\eta},t) \, d V_\eta \, ,
\label{eq:charspdoti}
\end{equation}
where 
%we have assumed that the Green function depends only on the distance between the 
%observation and evaluation points $\bs{x}=\bs{q}(\bs{\xi},t))$ and $\bs{y}=\bs{q}(\bs{\eta},t))$, respectively, and 
the integration is taken with the measure $d V_\eta$
of $\mathbb{R}^{n}$. In terms of these characteristic variables,  the system formed by equation~(\ref{eq:chars}), rewritten as 
\begin{equation}
{d \bs{q} \over d t} = \int_{\mathbb{R}^{n}} 
G_{ b-n/2}\big(|\bs{q}(\bs{\xi},t) -\bs{q}(\bs{\eta},t)|\big) 
\, \bs{p} (\bs{\eta},t) \, d V_\eta \, ,
\label{eq:charspqot}
\end{equation}
and equation~(\ref{eq:charspdoti}) constitutes the Lagrangian formulation of~equation~(\ref{eq:EPDIFFi}).  In this form, the equations of motion are a canonical 
Hamiltonian system with respect to variational derivatives $\delta/\delta\bs{q}$ and 
$\delta/\delta \bs{p}$ 
\begin{equation}
\dot{\bs{q}} (\bs{\xi},t) ={\delta H \over \delta \bs{p}}\, , 
\qquad \dot{\bs{p}} (\bs{\xi},t) =-{\delta H \over \delta \bs{q}}\,,
\label{eq:charsham}
\end{equation}
of the Hamiltonian functional 
\begin{equation}
H\equiv {1\over 2}\int_{\mathbb{R}^{n}} \int_{\mathbb{R}^{n}} 
G_{ b-n/2}\big(|\bs{q}(\bs{\xi},t) -\bs{q}(\bs{\eta},t)|\big) \, \bs{p}(\bs{\xi},t) \cdot \bs{p}(\bs{\eta},t) \, 
dV_\xi \, dV_\eta  \, .
\label{eq:ham}
\end{equation}
It is straightforward to check that substituting Eq. (\ref{eq:ham}) into Eq. (\ref{eq:charsham}) yields Eqs. (\ref{eq:charspqot}) and (\ref{eq:charspdoti}), which is equivalent to the model PDEs (\ref{eq:EPDIFF}). Hence the canonical Hamiltonian system forms our model equations. 

The Lagrangian version of equation~(\ref{eq:EPDIFFi}) shows that along 
characteristics~$\bs {q}(\bs{\xi},t)$ the evolution of 
the momentum-like variables $\bs {p}(\bs{\xi},t)$ is tied to that of the Jacobian matrix~${\partial_\beta q_\alpha }(\bs{\xi},t)$ by the initial conditions $\bs {p}(\bs{\xi},0)$,
\begin{equation}
p_\alpha(\bs{\xi},t)\, {\partial  q_\alpha \over \partial \xi_\beta}(\bs{\xi},t) = p_\beta(\bs{\xi},0) \,   
\label{eq:invp}
\end{equation}
(sum over repeated index), as it can readily be verified by system~(\ref{eq:charspdoti}),(\ref{eq:charspqot}) and the 
initial condition for characteristics~${\partial_\beta q_\alpha }(\bs{\xi},0)=\delta_{\alpha\beta}$. This is the analog of the constraint evolution for the one-dimensional  
SW equation~\cite{bib:chl06}, and can be used similarly to monitor the error of Lagrangian
numerical schemes to solve system~(\ref{eq:charsham}).  In Appendix \ref{sec:derivation}, we provide details on the connection of the Lagrangian formulation with the Eulerian form of system~(\ref{eq:EPDIFF}). 

\subsection{Dispersive deformation}\label{sec:dispersive}
One of the simplest settings removing the assumption of homogeneous boundary conditions is that of an infinite domain with $\bs{u}(\bs{x},\cdot) \to \bs{\kappa}$ as $\bs{x}\to \infty$ sufficiently fast, for some constant vector~$\bs{\kappa}$. This is most conveniently analyzed by defining the shifted field
$$
\bs{u}\equiv \tilde{\bs{u}}+\bs{\kappa}\, , 
$$
where $\uti$ is assumed to decay rapidly  at infinity. With the ``Galilean boost" 
\beq\label{eq:shift}
\tilde{\bs{x}}=\bs{x}-\bs{\kappa}t, \qquad \tilde{t}=t\, , 
\eeq
system~{eq:EPDIFF} maintains its form as the contributions from the boost and the $\bs{u}$ shift cancel out,
$$
%\beq\label{eq:epdiff_shifted}
\bs{m}_\ti 
-\bs{\kappa}\cdot \nati \bs{m} 
+ ((\uti+\bs{\kappa})\cdot\nati)\bs{m} 
+\bs{m} \cdot(\nati\bs{\uti})^{T}+\bs{m}(\nati\cdot\uti) = 0,
%\eeq
$$
with obvious meaning of the operator $\nati$. The formalism developed for homogeneous boundary conditions in free space can be applied by modifying the link between  $\bs{m}$ and $\bs{u}$ by the corresponding shift
$$
\mti \equiv \bs{m}+\bs{\kappa}
$$
so that domain of the operator $\mathcal{L}$ can remain the same (e.g., the Schwartz space for $\bs{\uti}$ initial data), and 
$$
\mti \equiv \mathcal{L} \uti
$$
so that  
\begin{equation}\label{eq:EPDIFFti}
\mti_\ti 
+ (\uti\cdot\nati)\mti 
+(\mti+\bs{\kappa}) \cdot(\nati\bs{\uti})^{T}+(\mti+\bs{\kappa})(\nati\cdot\uti) = 0.
\end{equation}

Dropping tildes from now on, this deformation can be cast in terms of characteristics as done for system~(\ref{eq:charspdoti}),(\ref{eq:charspqot}), by changing the boundary conditions for the momentum vector $\bs{p}$. 
If we let 
$$
\bs{p}(\bs{\xi},t) \to \bs{\kappa} \equiv {\rm const.} \quad {\rm as}  \quad |\bs{\xi}|\to \infty \, , 
$$ 
with
\begin{equation}
\bs{m}\big(\bs{q}(\bs{\xi},t),t\big)+\bs{\kappa} \equiv \,{\bs{p}(\bs{\xi},t)\over  J(\bs{\xi,t}) }  \, , 
\label{eq:charspk}
\end{equation}
and 
\begin{equation}
{d \bs{q} \over d t} = -\bs{\kappa}+\int_{\mathbb{R}^{n}} 
G_{ b-n/2}\big(|\bs{q}(\bs{\xi},t) -\bs{q}(\bs{\eta},t)|\big) 
\, \bs{p} (\bs{\eta},t) \, d V_\eta \, ,
\label{eq:charspqotk}
\end{equation}
the resulting system is the dispersive counterpart of the evolution equation in Eulerian form 
\begin{equation}\label{eq:EPDIFFk}
\bs{m}_t  + (\bs{u}\cdot\nabla)\bs{m} +(\bs{m}+\bs{\kappa})  \cdot(\nabla\bs{u})^{T}+(\bs{m}+\bs{\kappa})(\nabla\cdot\bs{u}) = 0. 
\end{equation}
The corresponding Hamiltonian is 
\begin{equation}
H\equiv {1\over 2}\int_{\mathbb{R}^{n}} \int_{\mathbb{R}^{n}} \big[
G_{ b-n/2}\big(|\bs{q}(\bs{\xi},t) -\bs{q}(\bs{\eta},t)|\big) \bs{p}(\bs{\xi},t)\cdot\bs{p}(\bs{\eta},t)-\bs{\kappa} \cdot (\bs{p}(\bs{\xi},t)+\bs{p}(\bs{\eta},t)) \big]\, 
dV_\xi \, dV_\eta  \, .
\label{eq:hamk}
\end{equation}

In this form, the system of the model PDEs (\ref{eq:EPDIFFk}) develops a non-trivial dispersion relation for 
the infinitesimal solutions $\bs{u} \to 0$. Linearizing around $\bs{u} = 0$ with 
$\bs{u}=\bs{U}\exp[i(\bs{k} \cdot \bs{x} - \omega t)]$ yields 
\begin{equation}
\omega=\frac{1}{(1+a^2 |\bs{k}|^2)^{ b}}\big(\bs{k} \cdot \bs{\kappa}\pm |\bs{k}||\bs{\kappa}|\big).
\label{eq:disper}
\end{equation}
The dispersion relation shows that when $\bs{k}$ is collinear with $\bs{\kappa}$
the corresponding phase speed $\bs{c}= \omega \, \bs{k}/|\bs{k}|$ can vanish. 
Thus, for unidimensional initial data $\bs{k}=k \, \bs{\kappa}/|\bs{\kappa}|$ 
the linear wave propagation is in fact unidirectional along the direction singled out 
by $\bs{\kappa}$. Note that, in general, the dispersion relation leads to non-trivial group velocity 
$\bs{C} \equiv \bs{\nabla}_{\bs{k}} \, \omega$, thus providing a dispersive mechanism for propagation of ``energy" 
away from localized initial conditions.

For the dispersive case in Lagrangian form, an equivalent formulation, more 
convenient for numerical purposes, can be provided in analogy with that for the 
one-dimensional SW equation presented in~\cite{bib:chl06}. Appending the Lagrangian form 
of the evolution equation~(\ref{eq:det_t}) for the determinant $J(\bs{\xi},t)$,  
\begin{equation}
{d J \over dt}=
J(\bs{\xi},t) \int_{\mathbb{R}^{n}} 
G'_{ b-n/2}\big(|\bs{q}(\bs{\xi},t) -\bs{q}(\bs{\eta},t)|\big) 
{(\bs{q}(\bs{\xi},t) -\bs{q}(\bs{\eta},t) ) \cdot (\bs{p} (\bs{\eta},t)-\bs{\kappa}J(\bs{\eta},t)) \over |\bs{q}(\bs{\xi},t) -\bs{q}(\bs{\eta},t)|}
 \, d V_\eta \, ,
 \label{eq:jchars}
\end{equation}
to the $(\bs{q},\bs{p})$ system~(\ref{eq:charsham}) allows the dispersive time 
evolution for $\bs{\kappa} \neq 0$ to be written equivalently as the system
\begin{equation}
\begin{split}
{d \bs{q} \over d t} &= \int_{\mathbb{R}^{n}} 
G_{ b-n/2}\big(|\bs{q}(\bs{\xi},t) -\bs{q}(\bs{\eta},t)|\big) 
\, (\bs{p} (\bs{\eta},t) -\bs{\kappa} J(\bs{\eta},t))\, d V_\eta \, ,
\\
{d \bs{p} \over d t} &= - \int_{\mathbb{R}^{n}} 
G'_{ b-n/2}\big(|\bs{q}(\bs{\xi},t) -\bs{q}(\bs{\eta},t)|\big) 
{\bs{q}(\bs{\xi},t) -\bs{q}(\bs{\eta},t) \over |\bs{q}(\bs{\xi},t) -\bs{q}(\bs{\eta},t)|}
\, \, \bs{p}(\bs{\xi},t) \cdot (\bs{p} (\bs{\eta},t) -\bs{\kappa} J(\bs{\eta},t)) \, d V_\eta %\, ,
%\\
%{d J \over dt}&=J (\bs{\xi},t)\int_{\mathbb{R}^{n}} 
%G'_{ b-n/2}\big(|\bs{q}(\bs{\xi},t) -\bs{q}(\bs{\eta},t)|\big) 
%{(\bs{q}(\bs{\xi},t) -\bs{q}(\bs{\eta},t) ) \cdot \bs{p} (\bs{\eta},t) \over |\bs{q}(\bs{\xi},t) -\bs{q}(\bs{\eta},t)|}
% \, d V_\eta 
 \, .
 \end{split}
 \label{eq:systkappa}
\end{equation}
(Details of the derivation of system~(\ref{eq:jchars})-(\ref{eq:systkappa}) are reported in Appendix~\ref{sec:derivation}.) Together with its initial conditions, $\bs{q}(\bs{\xi},0), \bs{p}(\bs{\xi},0)$ and 
$J(\bs{\xi},0)=1$,  the evolution system in the form~~(\ref{eq:jchars}) and~(\ref{eq:systkappa}) allows 
for a consistent treatment of the error associated with the numerical 
evaluation of the integrals, which is the foundation for the particle algorithm of Section \ref{sec:ps}. Note that the structure of the original system, 
~(\ref{eq:charspdoti})-(\ref{eq:charspqot}) with Hamiltonian~(\ref{eq:hamk}), is no
longer shared by the modified system (\ref{eq:jchars}) and (\ref{eq:systkappa}), as 
the appended Jacobian variable $J$ does not have a conjugate counterpart in this system.

\subsection{Green functions}
Unless mentioned otherwise, for this paper we will focus on the two-dimensional case, i.e. $n=2$,  for which the Green function reduces to
\beq\label{eq:2D_Green}
G_{ b-1}(|\bs{x}|) = \frac{2^{1- b}}{2\pi a^{1+ b}\, \Gamma( b)}|\bs{x}|^{ b -1}K_{ b -1}\left(\frac{|\bs{x}|}{a}\right).
\eeq
A notable special parametric choice is the two-dimensional Green function for $a=1$ and $ b=3/2$, for which it takes the particularly simple form 
\beq\label{eq:nu=3/2}
G_{1/2}(|\bs{x}|) = \frac{1}{2\pi}e^{-|\bs{x}|} \, . 
%\quad \text{where}\quad G_{1/2}(0) = \frac{1}{2\pi}.
\eeq
%since 
%\[
%\Gamma(\frac{3}{2}) = \frac{\sqrt{\pi}}{2},\quad\text{and}\,\,\, |\bs{x}|^{\frac{1}{2}}K_{\frac{1}{2}}(|\bs{x}|) = \sqrt{\frac{\pi}{2}}e^{-|\bs{x}|}.
%\]
The Green function in equation (\ref{eq:nu=3/2}) is continuous and radially symmetric around the origin~$\bs{x}=\bs{0}$, with a finite jump in radial derivative at the origin. A plot of the function resembles a cone whose peak is located at the origin. In fact, this function is a two-dimensional analog of the peakon solution of the SW equation (similarly to the one-dimensional peakon, this function is also a weak solution of equation~(\ref{eq:EPDIFF}), as further discussed below).

For other values of the parameter $ b$, the Green functions are expressed in terms of the Bessel function $K$. For instance, $ b=2$ and $a =1$, the Green function is 
\beq\label{eq:nu=2}
G_{1} (|\bs{x}|)=  \frac{1}{4\pi} |\bs{x}| K_{1}(|\bs{x}|), \quad \text{where}\quad G_{1}(0) = \frac{1}{4\pi}.
\eeq
 
\begin{figure}[tbh]
\begin{center}
\includegraphics[width=3.5in]{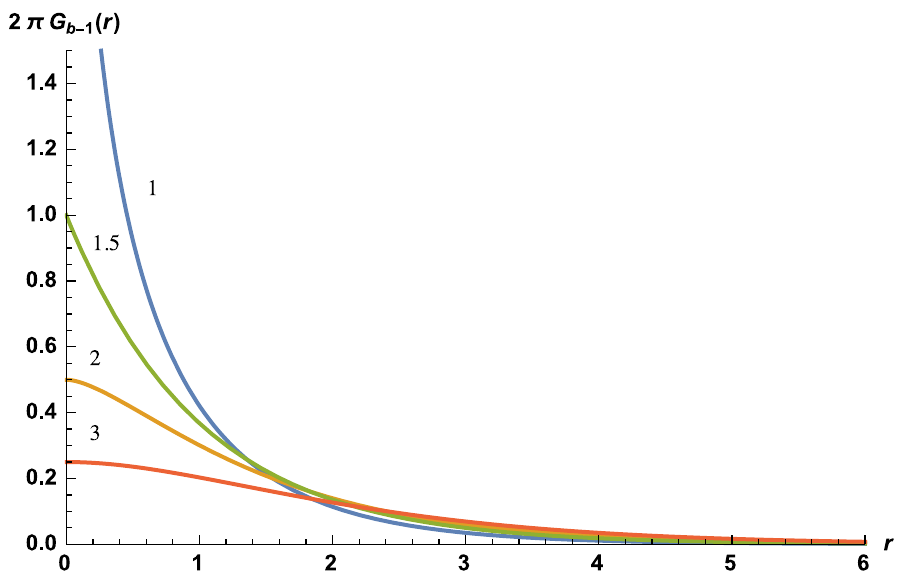}
\end{center}
\caption {Plots of $2\pi G_{ b-1}(r)$ for $ b=1, 1.5, 2$ and $3$, where $G_{ b-1}(r)$ is the two-dimensional Green functions of the Yukawa operator $\mathcal{L}^{ b}$. $a=1$ in the plots.}
\label{fig:Greens}
\end{figure}

The property of the Green function for various ranges of $ b$ is described as follows. For the range $1/4< b\le 1$ the Green function $G_{ b-1}(|\bs{x}|)$ is unbounded. For the range $1< b <3/2$ the function is bounded but non-differentiable at the peak, with the radial derivative suffering an infinite jump there (cusp). At $ b = 3/2$, the jump in radial derivative becomes finite. For the range $3/2 <  b \le 2$ the derivative of the function is continuous, but with an infinite second derivative at the peak. Similar intervals can be defined for higher smoothness properties of the solution. In particular, for $2 <  b<\infty$ the second derivative of the function is continuous. Figure \ref{fig:Greens} plots the function $2\pi G_{ b-1}(r)$ for the critical values $ b=1, 1.5, 2,$ and 3, respectively.

\section{$N$-particle system}
\label{sec:ps}
Replacing the integrals by the truncated Riemann sums in equations (\ref{eq:systkappa}) and (\ref{eq:jchars}) immediately yields a finite-dimensional $N$-particle system 
\beq\label{eq:Riemann-sum-pq}
\begin{split}
\ds\frac{d J_i}{d t} & = dxdyJ_i \sum_{\substack{j=1\\ j\ne i}}^NG_{ b-1}'(|\bs{q}_i-\bs{q}_j|)\frac{(\bs{q}_i-\bs{q}_j)\cdot(\bs{p}_j-\bs{\kappa}J_j)}{|\bs{q}_i-\bs{q}_j|}\\
\ds\frac{d \bs{q}_i}{dt}& = dxdy \sum \limits_{j=1}^N G_{ b-1}(|\bs{q}_i-\bs{q}_j|)(\bs{p}_j-\bs{\kappa}J_j),\\
\ds\frac{d\bs{p}_i}{dt}& = - dxdy\sum_{\substack{j=1\\ j\ne i}}^N(\bs{p}_i\cdot(\bs{p}_j-\bs{\kappa} J_j))G_{ b-1}'(|\bs{q}_i-\bs{q}_j|)\frac{\bs{q}_i-\bs{q}_j}{|\bs{q}_i -\bs{q}_j|}.
\end{split}
\eeq
The field $\bs{u}$ can be recovered by 
\beq\label{eq:Riemann-sum-u}
\bs{u}(\bs{x}, t) = dxdy \sum_{j=1}^{N}G_{ b-1}(|\bs{x} - \bs{q}_j|)(\bs{p}_j-\bs{\kappa}J_j).
\eeq

An alternative viewpoint, proposed by Mumford \& Desolneux \cite{bib:mumford} for the nondispersive case $\bs{\kappa}=0$,  is to obtain equations (\ref{eq:Riemann-sum-pq})  by the {ansatz} 
%for system~(\ref{eq:EPDIFF})
\beq\label{eq:u}
\bs{u}(\bs{x}, t) = \sum_{j=1}^{N}G_{ b-1}(|\bs{x} - \bs{q}_j|)\bs{p}_j,
\eeq
where 
\begin{eqnarray}
\ds\frac{d \bs{q}_i}{dt}=u(\bs{q}_i, t)=\sum \limits_{j=1}^N G_{ b-1}(|\bs{q}_i-\bs{q}_j|)\bs{p}_j,\label{eq:Q}\\
\bs{m}(\bs{x}, t) =\mathcal{L}^{b}u(\bs{x}, t) = \sum_{j=1}^{N} \bs{p}_j\delta(\bs{x}-\bs{q}_j).\label{eq:m}
\end{eqnarray}
Substituting this {ansatz} into the {\it weak formulation} of system~(\ref{eq:EPDIFF}) with respect to 
an appropriate test-function space~\cite{bib:mumford} 
 yields an equation for the $\bs{p}_i$'s that closes the $\{\bs{q},\bs{p}\}$~system, i.e., the 
%\beq\label{eq:P}
%\ds\frac{d\bs{p}_i}{dt} = - \sum_{\substack{j=1\\ j\ne i}}^N(\bs{p}_i\cdot\bs{p}_j)\ds\frac{\partial G(|\bs{q}_i-\bs{q}_j|)}{\partial \bs{q}_i} = -\sum_{\substack{j=1\\ j\ne i}}^N(\bs{p}_i\cdot\bs{p}_j)G'(|\bs{q}_i-\bs{q}_j|)\frac{\bs{q}_i-\bs{q}_j}{|\bs{q}_i -\bs{q}_j|}.
%\eeq
finite-dimensional $N$-particle system 
\beq\label{eq:N-particle}
\begin{split}
\ds\frac{d \bs{q}_i}{dt}& =\sum \limits_{j=1}^N G_{ b-1}(|\bs{q}_i-\bs{q}_j|)\bs{p}_j,\\
\ds\frac{d\bs{p}_i}{dt}& = -\sum_{\substack{j=1\\ j\ne i}}^N(\bs{p}_i\cdot\bs{p}_j)G_{ b-1}'(|\bs{q}_i-\bs{q}_j|)\frac{\bs{q}_i-\bs{q}_j}{|\bs{q}_i -\bs{q}_j|},
\end{split}
\eeq
where $i=1,\cdots, N$ for equations~(\ref{eq:EPDIFF}). Note that equations~(\ref{eq:Riemann-sum-pq}) for  $\bs{\kappa}=0$ and (\ref{eq:N-particle}) are equivalent, since $dxdy$ can be scaled into the momentum variable $\bs{p}_j$, although the interpretation of the system is somewhat different in the two approaches, as the Lagrangian derivation bypasses the weak formulation of the evolution equation. 

By denoting with $\bs{q}_i$ and $\bs{p}_i$ the 2-vectors 
\[
\bs{q}_i = \left[\begin{array}{cc}
q^1_i\\ q^2_i\end{array} \right],\quad \bs{p}_i = \left[\begin{array}{cc}
p^1_i\\ p^2_i\end{array} \right],
\]
the inner product in equation (\ref{eq:N-particle}) is $\bs{p}_i\cdot\bs{p}_j = \bs{p}_i^{T}\bs{p}_j = p^1_i p^1_j + p^2_i p^2_j $.
For the special case $ b=3/2$ and $a=1$, as mentioned previously, we have $G_{1/2}(|\bs{x}|)= e^{-|\bs{x}|}/{2\pi}$ and $G_{1/2}'(|\bs{x}|)=- e^{-|\bs{x}|}/{2\pi}$. For other values of $ b$, we recall the recursive formula for the modified Bessel function of second kind 
for real~$b$~\cite{bib:as65},
\beq
\frac{d}{dr}\left[r^{ b}K_{ b}\right] = - r^{ b}K_{ b-1}.
\eeq
Thus, in two dimensions, we have
%\beq
%G'(r)=\frac{d}{dr}G_{ b-1}(r) = \frac{1}{2\pi\alpha^{1+ b}\Gamma( b)}\left[r^{ b -2}K_{ b -1}\left(\frac{r}{\alpha}\right) +r^{ b-1}\alpha^{-1}K'_{ b-1}\left(\frac{r}{\alpha}\right)\right],\,\,r\ne 0.
%\eeq 
\beq
\begin{split}
G'_{ b-1}(r)  =\frac{d}{dr}G_{ b-1}(r)  & = \frac{2^{1- b}}{2\pi a^{1+ b}\Gamma( b)}\frac{d}{dr}\left[r^{ b-1}K_{ b-1}(\frac{r}{ a})\right] \\ &= -\frac{2^{1- b}}{2\pi a^{2+ b}\Gamma( b)} r^{ b-1}K_{ b-2}(\frac{r}{a})\\
& = -\frac{r}{2( b-1)a^2}\left[\frac{2^{2- b}}{2\pi a^{ b}\Gamma( b-1)} r^{ b-2}K_{ b-2}(\frac{r}{a})\right] \\
& = -\frac{r}{2( b-1)a^2}G_{ b-2}(r),\,\,\, r\ne 0.                         
%&=-\left(\fac{r}{ a b}\right)\frac{1}{2\pi a^{ b}\Gamma( b-1)}r^{ b-2}K_{ b-2} \\
%& = -\left(\frac{r}{ a b}\right) G_{ b-2}(r),\,\,\ r\ne 0.
\end{split}
\eeq

\begin{note}\label{note:Lip}
For $ b<3/2$,  $G_{ b-1}$ is not differentiable at zero and $G'_{ b-1}(0)\rightarrow\infty$. For $b=3/2$ the radial derivative $G'_{ b-1}$ is discontinuous at zero , which is a bounded discontinuity for the $\bs{p}$ equation in the particle system.  For $b=2$,  $G'_{ b-1}$ is continuous, but is not Lipschitz continuous at zero. In general (see, e.g., \cite{bib:perko}), if $F$ in the ODE system $\dot{\bs Y} = F(\bs{Y})$ is not continuous, the existence of the solution of the ODE is not guaranteed. Furthermore, if $F$ is not Lipschitz continuous, the uniqueness of the solution of the ODE is not guaranteed. Hence for $b=3/2$, the existence of the solution of the particle system for particle collision is not guaranteed, and likewise for $b=2$, solution uniqueness may fail.
For $ b>2$,  $G'_{ b-1}$  is differentiable at zero, and hence existence and uniqueness of solutions hold.  
%Moreover, for $ b>2$, since the Green function $G_{ b-1}$ is symmetric and continuously differentiable at zero, it is natural to set $G_{ b-1}'(0)=0$ in this range.
%We refer to the Green's kernel with $G'_{ b-1}$ that is not Lipschitz continuous at zero as a locally non-Lipschitz particle. 

\end{note}

\begin{note}
Without further specification, for the rest of the paper, we will only consider the case $ a=1$ for our analysis and numerical examples.
\end{note}

\section{Traveling wave solutions}

The (nondispersive) system~(\ref{eq:EPDIFF}) admits the traveling wave solution, 
\beq\label{eq:traveling}
 \bs{u}(\bs{x}, t)= \left(\bs{p}/G_{ b-1}(0)\right)G_{ b-1}(|\bs{x}-\bs{x}_0 - t\, \bs{p}|),
\eeq
for some constant vector $\bs{p}$. $G_{ b-1}(0)$ is the Green function evaluated at the origin. At $t=0$ the wave is centered at $\bs{x}_0$, and the initial condition of $\bs{u}$ is
\beq\label{eq:traveling_u_0}
\bs{u}_0 =\bs{u}(\bs{x}, 0)= (\bs{p}/G_{ b-1}(0))G_{ b-1}(|\bs{x}-\bs{x}_0|).
\eeq
The behavior of the traveling wave depends on the Green function of the elliptic operator. For $ b$ in the range of $1/4< b<1$, the traveling wave solution moves along the vector $\bs{p}$ with a moving unbounded point $\bs{x}= \bs{x}_0+\bs{p}t$ at the center. For the range $1< b<3/2$, the center is bounded but its radial derivative is unbounded. At $ b=3/2$, the center becomes continuous, but its radial derivative has a finite jump, i.e., a two-dimensional {peakon}, which, because of its conical shape, we will henceforth refer to as a ``conon."

The traveling-wave solution can be easily verified by placing only one particle at $\bs{x}_0$ initially in the $N$ particle system, i.e. $N=1$ and $\bs{q}_1(0) =\bs{x}_0$, with an unknown initial momentum $\bs{p}_1(0)$. Then, by using the initial data of the traveling wave (\ref{eq:traveling_u_0}), one can find this initial momentum. Recall the definition of $\bs{m}$,
\beq\label{eq:m_0-travel}
\bs{m}_0=\bs{m}(\bs{x}, 0) = \mathcal{L}^{b}\bs{u}_0 = \left(\bs{p}/G_{ b-1}(0)\right)\mathcal{L}^{b}G_{ b-1}(|\bs{x}-\bs{x}_0|) = \left(\bs{p}/G_{ b-1}(0)\right)\delta(\bs{x}-\bs{x}_0);
\eeq
by comparing equations (\ref{eq:m}) and (\ref{eq:m_0-travel}), we obtain
\beq
\bs{m}(\bs{x}, 0) = \left(\bs{p}/G_{ b-1}(0)\right)\delta(\bs{x}-\bs{x}_0) = \bs{p}_1(0)\delta(\bs{x}-\bs{q}_1(0))= \bs{p}_1(0)\delta(\bs{x}-\bs{x}_0), 
\eeq
and thus $\bs{p}_1(0) = \bs{p} /G_{ b-1}(0)$. Given $\bs{q}_1(0)$ and $\bs{p}_1(0)$, 
the one-particle system is simply 
\beq\label{eq:one-particle}
\begin{split}
\ds\frac{d \bs{q}_1(t)}{dt}& =G_{ b-1}(0)\bs{p}_1(t),\\
\ds\frac{d\bs{p}_1(t)}{dt}& = \bs{0}.
\end{split}
\eeq
Integration of the first system of ODE gives $\bs{q}_1(t) = G_{ b-1}(0)\bs{p}_1(0) t + \bs{q}_1(0)= \bs{p}t + \bs{x}_0$. From equation (\ref{eq:u}), the field $\bs{u}$ is 
then reconstructed by 
\beq
\bs{u}(\bs{x}, t) = G_{ b-1}(|\bs{x} - \bs{q}_1(t)|)\bs{p}_1(t)=  (\bs{p}/G(0))G_{ b-1}(|\bs{x} - \bs{x}_0 - t\, \bs{p}|),
\eeq
since $\bs{p}_1(t) =\bs{p}_1(0) = \bs{p}/G_{ b-1}(0)$. Thus,  the solution obtained by the $N$-particle system using the initial data of the traveling wave is consistent with the exact traveling-wave solution at later times.

%As an example, we consider the special case $ b=3/2$ and $\bs{x}_0=\bs{0}$. The constant vector $\bs{p}$ is chosen to be $(1, 1)$. Figure \ref{fig:traveling_wave} shows the waveforms of the first component of the traveling wave at $t=0$ and $t=10$, respectively. It can be seen that the conon solution behaves as expected, 
%%, resembling a cone whose radial derivative at the peak has a finite jump, 
%traveling from the initial location $(0, 0)$ to the location $(10, 10)$ at $t=10$. The simulations are constructed on a $41\times 41$ mesh grid in the domain $[-5, 15]\times [-5,15]$. ???need to discuss this, it is simply an illustration obtained by matlab plots or is there something computed here????
%\begin{figure}[tbh]
%\begin{center}
%(a)\includegraphics[width=2.75in]{t0_traveling-eps-converted-to.pdf}
%(b)\includegraphics[width=2.75in]{t10_traveling-eps-converted-to.pdf}
%\end{center}
%\caption{(a) Initial waveform of the first component of the traveling wave. (b) The first component of the traveling wave at $t=10$. We use $ b=3/2$, $\bs{x}_0=\bs{0}$, and $\bs{p} = (1, 1)^T$ in the simulations.}
%\label{fig:traveling_wave}
%\end{figure}
%
%\vskip 12pt

%For the case of traveling wave, we have 
%\beq
%\bs{u}(\bs{x}_i, 0) = (\bs{p}/G(\bs{0}))G(|\bs{x}_i-\bs{x}_0|).
%\eeq

% Equation (\ref{eq:traveling}) becomes
%\beq\label{eq:traveling_scaled}
%\bs{u}(\bs{x}, t)= \bs{p}\hat{G}(|\bs{x}-\bs{x}_0 - t\bs{p}|),
%\eeq

\subsection{Normalization of the Green functions}
It is easy to check that the constant in front of the Green function in equation (\ref{eq:Green-nD}) can be absorbed into a time rescaling. For our computational purpose, it may be convenient to normalize the Green function as it were an element of a basis system. If we normalize the Green function by $G_{ b-1}(0)$ and introduce the pair of scaled functions  
\beq
\tilde{G}_{ b-1}(r) = \frac{G_{ b-1}(r)}{G_{ b-1}(0)},\quad \tilde{\bs{p}}_j = G_{ b-1}(0)\bs{p}_j
\eeq
then the {ansatz} for system~(\ref{eq:EPDIFF}) becomes
\begin{eqnarray}
\bs{u}(\bs{x}, t) = \sum_{j=1}^{N}\tilde{G}_{ b-1}(|\bs{x} - \bs{q}_j|)\tilde{\bs{p}}_j,\label{eq:new-anstz}\\
\ds\frac{d \bs{q}_i}{dt}=\bs{u}(\bs{q}_i, t)=\sum \limits_{j=1}^N\tilde{G}_{ b-1}(|\bs{q}_i-\bs{q}_j|)\tilde{\bs{p}}_j,\label{eq:new-Q}\\
\bs{m}(\bs{x}, t) =\mathcal{L}^{ b}\bs{u}(\bs{x}, t) = \frac{1}{G_{ b-1}(0)}\sum_{j=1}^{N} \tilde{\bs{p}}_j\delta(\bs{x}-\bs{q}_j).\label{eq:new-m}
\end{eqnarray}
and the equation for $\tilde{\bs{p}}_j$ is
\beq\label{eq:new-P}
\ds\frac{d\tilde{\bs{p}}_i}{dt} = -\sum_{\substack{j=1\\ j\ne i}}^N(\tilde{\bs{p}}_i\cdot\tilde{\bs{p}}_j)\tilde{G}_{ b-1}'(|\bs{q}_i-\bs{q}_j|)\frac{\bs{q}_i-\bs{q}_j}{|\bs{q}_i -\bs{q}_j|}.
\eeq
%equations (\ref{eq:new-anstz}) to (\ref{eq:new-P}) are equivalent to equations (\ref{eq:u}) to (\ref{eq:P}). 
The scaled Green function and momentum give rise to the traveling wave solution 
\beq\label{eq:traveling-new}
 \bs{u}(\bs{x}, t)= \bs{p}\tilde{G}(|\bs{x}-\bs{x}_0 - t\bs{p}|).
\eeq
The above solution can be verified by the scaled one-particle system
\beq\label{eq:one-particle-scaled}
\begin{split}
\ds\frac{d \bs{q}_1(t)}{dt}& =\tilde{G}_{ b-1}(0)\tilde{\bs{p}}_1(t),\\
\ds\frac{d\tilde{\bs{p}}_1(t)}{dt}& = \bs{0},
\end{split}
\eeq
where $\tilde{\bs{p}}_1(0)=\bs{p}$.  It is worth noting that for $ b=3/2$, $G_{1/2}(0) =\ds\frac{1}{2\pi}$, and for $ b=2$, $G_{1}(0)=\ds\frac{1}{4\pi}$.

\subsection{An example of traveling wave}
A numerical test of the particle algorithms is offered by the traveling solution of system~(\ref{eq:EPDIFF}).  For $N>1$, one way to obtain the momenta on a mesh from a given $\bs{u}$ is to use equation (\ref{eq:u}) instead of equation (\ref{eq:m}). Suppose that $N$ particles are placed on a mesh initially. The initial locations of the particles are at the mesh grid, i.e. $\{\bs{q}_1,\dots,\bs{q}_N\}=\{\bs{x}_1,\dots,\bs{x}_N\}$, and hence $\bs{u}$ for the $i^{th}$ particle is
\beq
\bs{u}(\bs{x}_i, 0) =\sum_{j=1}^N G_{ b-1}(|\bs{x}_i-\bs{x}_j|)\bs{p}_j(0).
\eeq
The above equation in matrix-vector form is the linear system 
%\beq\label{eq:linear} \left(\begin{array}{cccc}
%\bs{u}(\bs{x}_1,0)\\ \bs{u}(\bs{x}_2,0)\\ \vdots\\ \bs{u}(\bs{x}_N,0)\end{array} \right)=\left( \begin{array}{cccc}
%G_{ b-1}(|\bs{x}_1-\bs{x}_1|) & G_{ b-1}(|\bs{x}_1-\bs{x}_2|) &\cdots & G_{ b-1}(|\bs{x}_1-\bs{x}_N|)\\
%G_{ b-1}(|\bs{x}_2-\bs{x}_1|) & G_{ b-1}(|\bs{x}_2-\bs{x}_2|) &\cdots & G_{ b-1}(|\bs{x}_2-\bs{x}_N|) \\
%\vdots & \vdots & \vdots &\vdots\\
%G_{ b-1}(|\bs{x}_N-\bs{x}_1|) & G_{ b-1}(|\bs{x}_N-\bs{x}_1|) &\cdots & G_{ b-1}(|\bs{x}_N-\bs{x}_N|)
%\end{array} \right)
%\left(\begin{array}{cccc}
%\bs{p}_1(0)\\ \bs{p}_2(0)\\ \vdots\\ \bs{p}_N(0)\end{array} \right).\eeq

\beq\label{eq:linear} \left(\begin{array}{ccc}
\bs{u}(\bs{x}_1,0)\\ \vdots\\ \bs{u}(\bs{x}_N,0)\end{array} \right)=\left( \begin{array}{ccc}
G_{ b-1}(|\bs{x}_1-\bs{x}_1|) &\cdots & G_{ b-1}(|\bs{x}_1-\bs{x}_N|)\\
\vdots & \vdots &\vdots\\
G_{ b-1}(|\bs{x}_N-\bs{x}_1|) &\cdots & G_{ b-1}(|\bs{x}_N-\bs{x}_N|)
\end{array} \right)
\left(\begin{array}{ccc}
\bs{p}_1(0)\\  \vdots\\ \bs{p}_N(0)\end{array} \right).
\eeq
Inverting the system, we obtain the initial momenta $\bs{p}_j(0)$, $j=1\dots N$ for the $N$-particle system.

We consider the scaled traveling waves (\ref{eq:traveling-new}) and the scaled $N$-particle system (\ref{eq:new-Q})-(\ref{eq:new-P}). We first use equation (\ref{eq:linear}) with scaled Green functions to find the initial particle momenta $\tilde{\bs{p}}_i(0)$, $i=1,\dots, N$. Then we evolve the $N$-particle system to some finite time. Finally, we use the particle locations and momenta to reconstruct the field $\bs{u}$. 

\begin{figure}[tbh]
\begin{center}
%(a) \includegraphics[width=2.75in]{initial_N21-eps-converted-to.pdf}
(a)\includegraphics[width=2.75in]{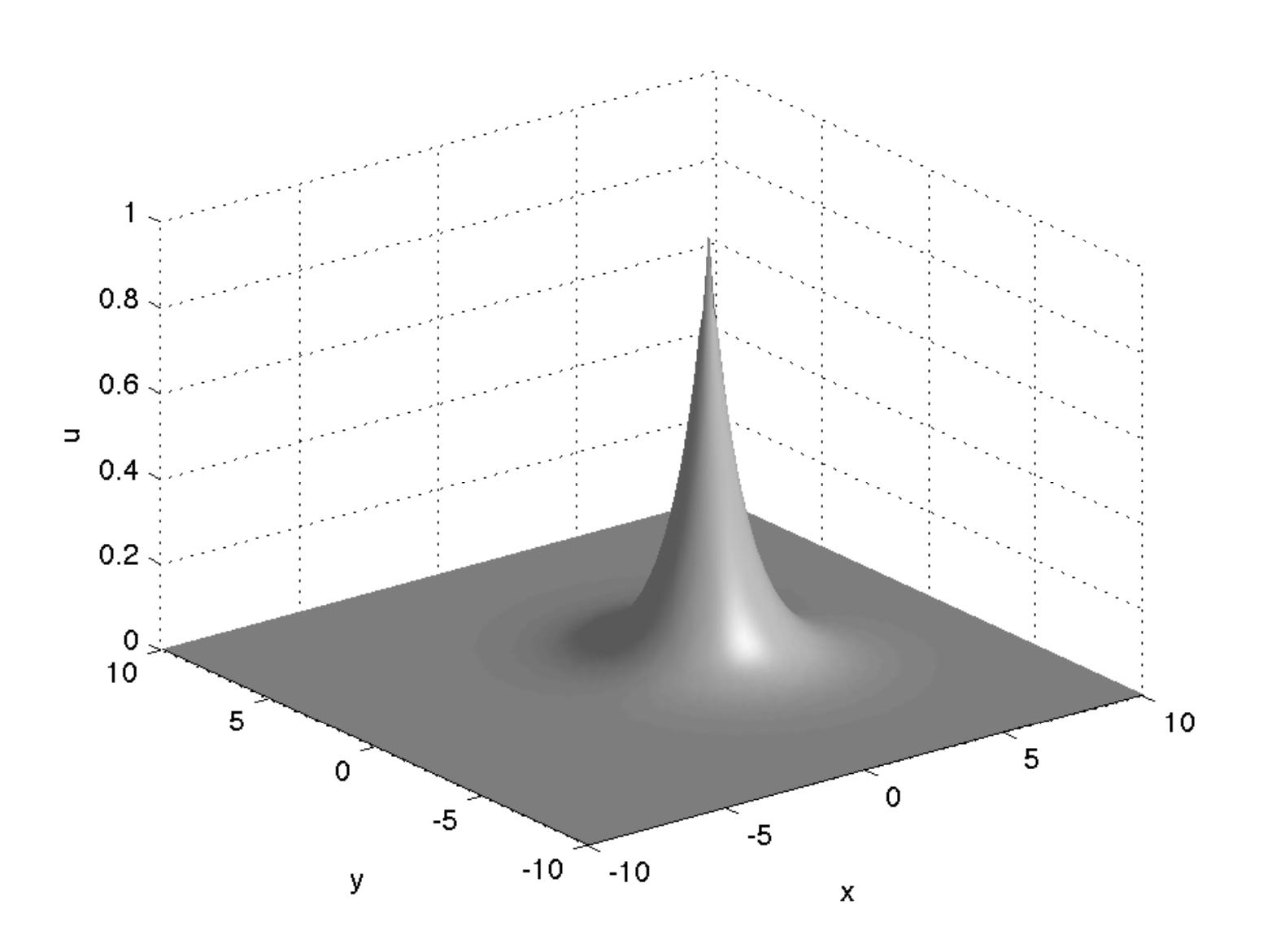}
(b)\includegraphics[width=2.75in]{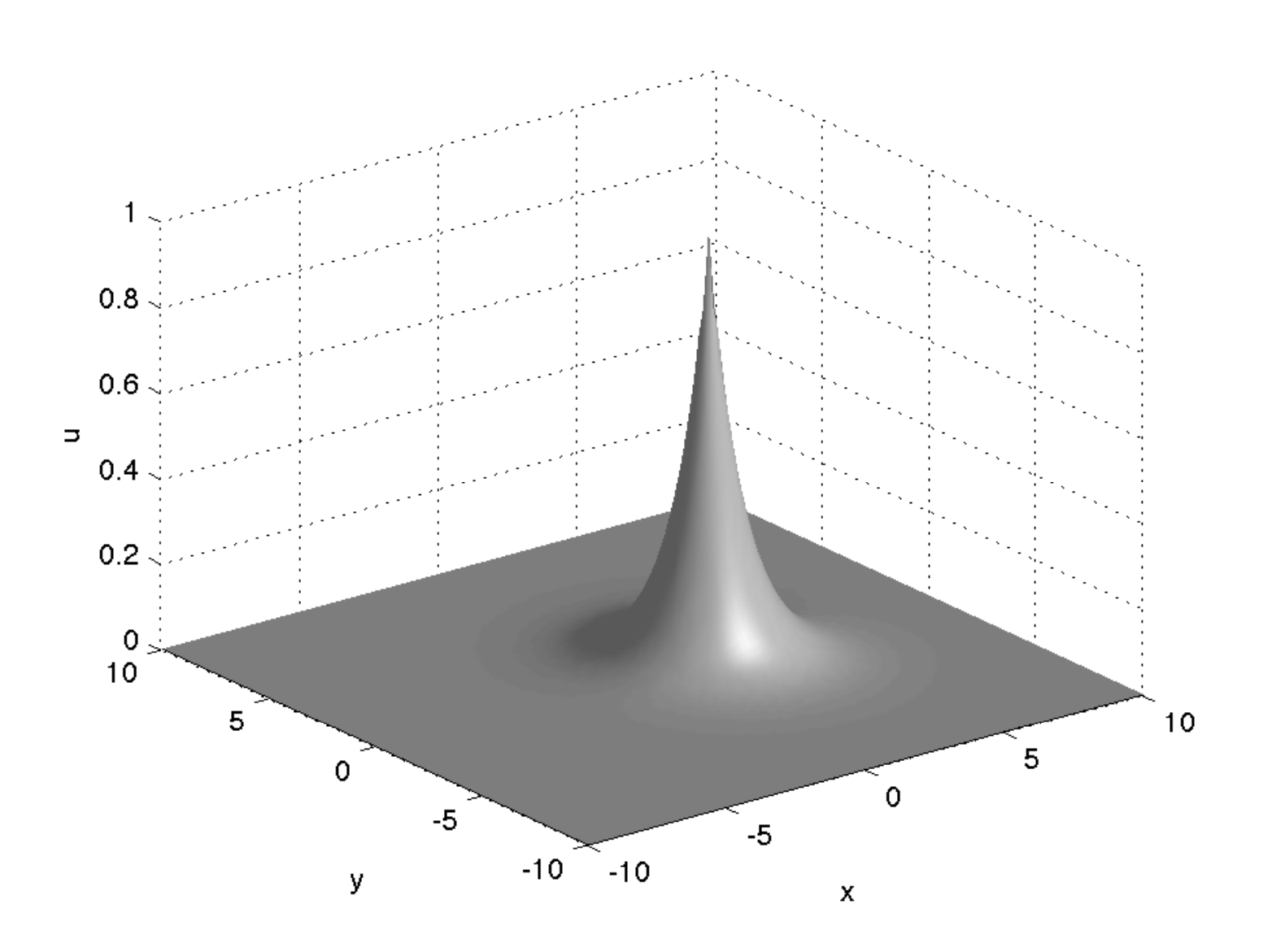}
\end{center}
\caption{$N$ particles, $N=1681$, are placed initially on a $41\times 41$ mesh in the domain of $[-10, 10]\times [-10, 10]$. The solution $\bs{u}$ is reconstructed on a $101\times 101$ mesh from the solutions of the $N$-particle system. The constant vector is $\bs{p} = (1, 0)$ and $ b=3/2$. (a) The first component of the exact traveling-wave solution at $t=2$. (b) The computed solution for (a). The 2-norm error is $2.0461\times 10^{-15}$.}
\label{fig:traveling_n41}
\end{figure}

%For the numerical experiments in this section, the initial velocity field of a traveling wave (\ref{eq:traveling_u_0}) is place onto a mesh based on the following two characteristics. Let the mesh points be $\{\bs{x}_1,\dots,\bs{x}_N\}$ and $\bs{q}_i(0)=\bs{x}_i$.
%\begin{enumerate}
%\item {\bf Case I:}  $\bs{x}_0 \in \{\bs{q}_1(0),\dots,\bs{q}_N(0)\}$. One of the particles is placed at the center of the traveling wave $\bs{x}_0$. i. e. $\bs{x}_0 = \bs{q}_j$(0) for some $j$.
%\item {\bf Case II:}  $\bs{x}_0 \not\in \{\bs{q}_1(0),\dots,\bs{q}_N(0)\}$. No particle is placed at the center of the traveling wave $\bs{x}_0$.
%\end{enumerate}
A traveling wave, $\bs{p}\tilde{G}_{ b-1}(|\bs{x}-\bs{x}_0-\bs{p}t|)$, where $\bs{p}=(1, 0)$ and $\bs{x}_0=\bs{0}$, is placed on a two-dimensional mesh in the domain  $D=[-10, 10]\times [-10, 10]$. We consider the case that  the operator $\mathcal{L}^{ b}$ has power $ b=3/2$, and its Green function $G_{ b-1}$ is described as in equation (\ref{eq:nu=3/2}), divided 
by ${2\pi}$. The initial data $\bs{q}_i(0)$ and $\tilde{\bs{p}}_i(0)$, $i=1,\cdots N$ for the $N$-particle system (\ref{eq:N-particle}) are obtained as follows. We initially place $N$ particles on a $41\times 41$ mesh over the domain $D$  ($N=1681$). We solve the linear system (\ref{eq:linear}) to obtained $\tilde{\bs{p}}_{i}$ for the $N$-particles. We remark that a single particle of a given amplitude would yield a traveling wave solution of the PDE with trivial evolution. This cannot in general be seen by assigning this as an initial condition to the $\bs{u}$-field for the numerical particle algorithm. Instead, the discretization of the initial data $\bs{u}$ would yield a particle system with as many particles as the initial grid points. 
%One of the particles is placed at the center of the traveling wave, i.e.,  $\bs{q}_{841}(0)=\bs{x}_0$. 
%To obtain $\tilde{\bs{p}}_i(0)$, we solve the linear system (\ref{eq:linear}), and discover that $\tilde{\bs{p}}_{841}(0) =\bs{p}$ (counting in the row direction from south-west corner) , and $\tilde{\bs{p}}_{i}(0)=0$ for $i$ otherwise. Using these $\bs{q}_i(0)$ and $\tilde{\bs{p}}_i(0)$, we solve the $N$-particle system (\ref{eq:new-Q})-(\ref{eq:new-P}) by using an explicit second-order Runge-Kutta method with two-stages 
The $N$-particle system (\ref{eq:new-Q})-(\ref{eq:new-P}) solved by using an explicit second-order Runge-Kutta method with two-stages
\beq\label{eq:2ndRK}
y^{n+1}=y^{n} + \Delta t f\left(t^n+\frac{1}{2}\Delta t, y^{n}+\frac{1}{2}\Delta t f(t^{n}, y^{n})\right).
\eeq
Figure \ref{fig:traveling_n41}(a) shows the first component of the exact traveling wave solution at $t=2$, and  Figure \ref{fig:traveling_n41}(b) is the computed counterpart. The 2-norm error for the computed solution is $2.0461\times 10^{-15}$, with $2$-norm defined as 
\beq\label{eq:2norm}
||e||_{2} = \sqrt{dx dy\sum_{i=1}^{n}\sum_{j=1}^{n} e^2_{i,j} }.
\eeq
The solution $\bs{u}$ is reconstructed on a $101\times 101$ mesh points in the domain of $[-10, 10]\times [-10, 10]$ from the solutions of the $N$-particle system.

%\vskip 12pt
%
%\noindent
%{\bf Remark:}
%We note that since only one particle has non-zero momentum initially, solving the $N$-particle system is equivalent to solving the coupled linear one-particle system (\ref{eq:one-particle}) by the Runge-Kutta method. Moreover, our numerical experiments confirm that as long as there is one particle placed at the location of $\bs{x}_0$, the momenta $\tilde{\bs{p}}_i(0)$, $i=1,\dots, N$ are always ????has $\tilde{\bs{p}}_{\bs{x}_0} = \bs{p}$, and zeros (or very small numbers) everywhere else, no matter what the size of $N$ we choose.

\section{Two-particle dynamics}

\subsection{Phase portrait for $b=3/2$}\label{sec:nu=1} As remarked in Note \ref{note:Lip}, the existence and/or uniqueness of the solution of the two-particle system at zero are not guaranteed for $b=3/2$ and $b=2$, for which the Green kernels have bounded discontinuity or non-Lipschitz continuity, respectively. 
In this section, we investigate the two particle system for these two special cases, to illustrate these existence issues.  In particular, we focus on the solution of particle collisions. As we will see, while exact solutions by quadrature are possible, the issue of how to continue past a collision can arise, and this can be overcome by imposing a conservation a law such as that of the Hamiltonian. However, when solving the two particle ODE system numerically, such conservation would depend on the algorithm, and it will be seen that the way continuation past collision is selected (if at all possible) can in fact depend on the details of the numerical scheme and on its parameters.

\begin{figure}[tbh]
\begin{center}
(a)\includegraphics[width=2.75in]{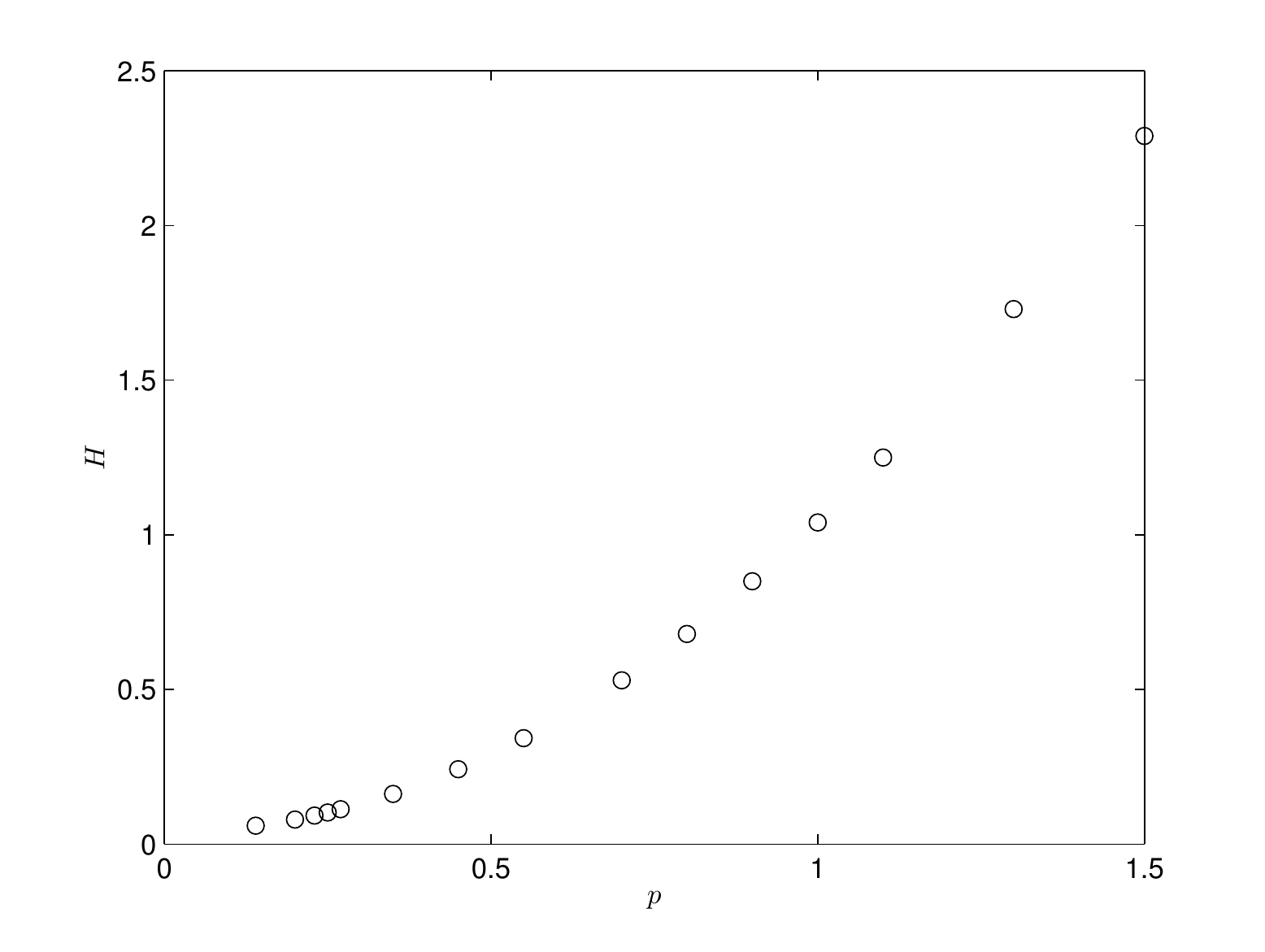}
(b)\includegraphics[width=2.75in]{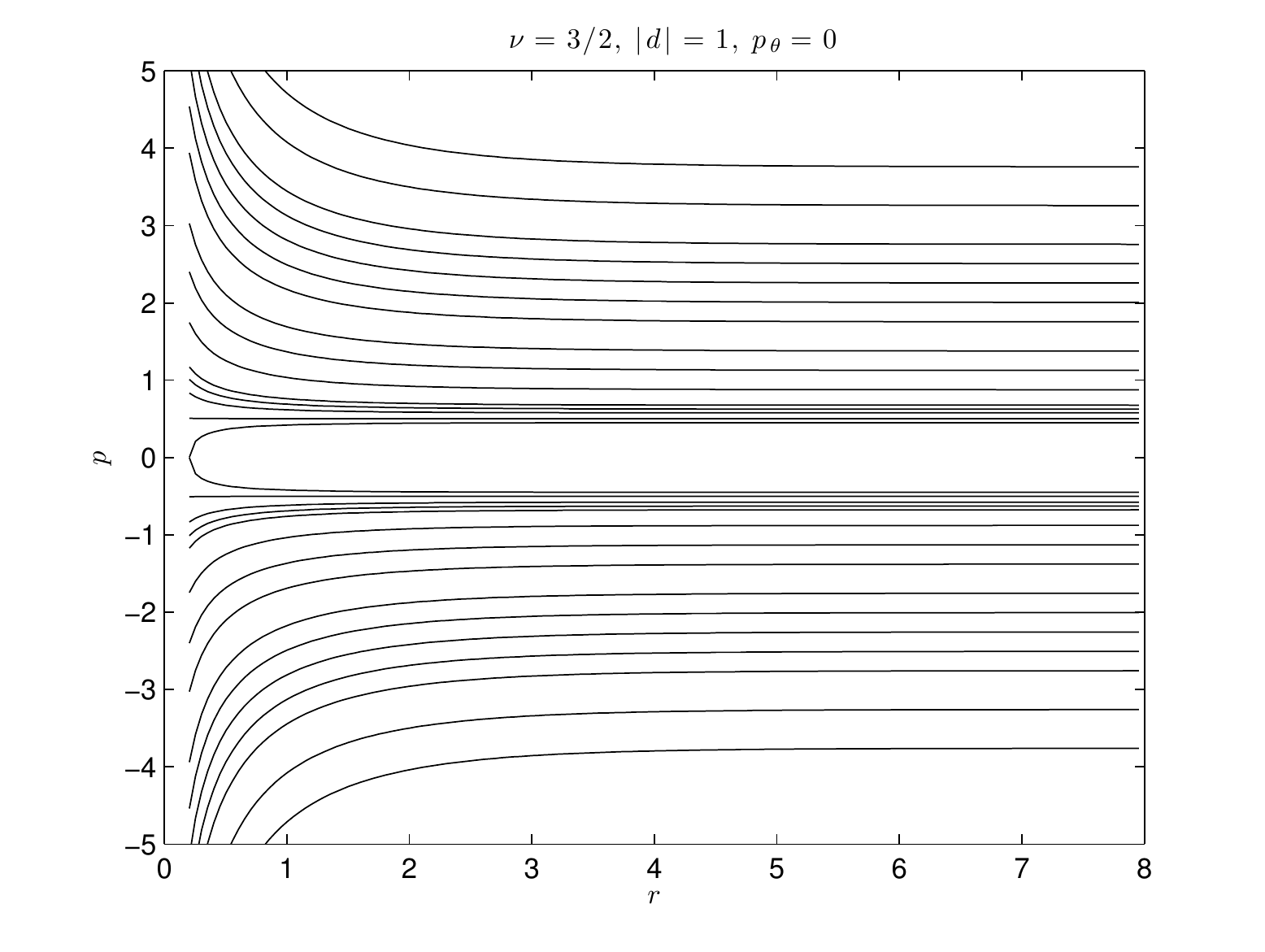}
\end{center}
\caption{Phase portrait for $|\bs{d}|=1$, $p_{\theta}=0$ for $ b=3/2$. (a) The Hamiltonian values computed with the given $|\bs{d}|$, $p_{\theta}$, and $p$ by using equation (\ref{eq:Hamiltonian}) at $r=8$. (b) Three principle behaviors are exhibited in the phase portrait for the two-particle interaction. The ejection orbits are in the upper-half plane, whereas the capturing orbits are in the lower-half plane. The scattering orbits are in the middle.}
\label{fig:phase_portrait1}
\end{figure}

For the phase-portrait analysis of two-particle interaction, we adopt the approaches in \cite{NP, Landau} and define the Hamiltonian 
\beq\label{eq:general_Hamiltonian}
H= \frac{1}{2}\sum_{i, j=1}^{N}\bs{p}_i\cdot\bs{p}_jG_{ b-1}(|\bs{q}_i-\bs{q}_j|).
\eeq
The Hamiltonian $H$ is conserved \cite{NP}. 
%The other generic conserved quantities for the $N$-particle system (\ref{eq:N-particle}) are the linear momentum $\sum_{i=1}^{N}\bs{p}_i$ and the angular momentum $\sum_{i=1}^{N}\bs{q}_i\times \bs{p}_i$ \cite{NP}. 
If $N=2$,
\beq  
H = \bs{p}_1\cdot\bs{p}_2G_{ b-1}(|\bs{q}_1-\bs{q}_2|) + \frac{1}{2}\left(|\bs{p}_1|^2 +  |\bs{p}_2|^2 \right)G_{ b-1}(0).
\eeq
Let 
\beq
\bs{d} =\bs{p}_1 + \bs{p}_2,\quad \bs{p}=(\bs{p}_1-\bs{p}_2)/2,\quad{\bs s} = (\bs{q}_1+\bs{q}_2)/2,\quad\bs{q} = \bs{q}_1 -\bs{q}_2,
\eeq
then
\beq\label{eq:N=2-H}
H=\left(\frac{1}{4}|\bs{d}|^2-|\bs{p}|^2\right)G_{ b-1}(|\bs{q}|) + \left(\frac{1}{4}|\bs{d}|^2+|\bs{p}|^2\right)G_{ b-1}(0).
\eeq
Parameterizing $\bs{q}$ and $\bs{p}$ in the polar coordinates yields
\beq
\bs{q} = (r\cos\theta, r\sin\theta),\quad\bs{p} = (p\cos\theta-p_{\theta}\sin\theta/r, p\sin\theta-p_{\theta}\cos\theta/r),
\eeq
where $r=|\bs{q}|$, $\theta$ is the angle between $\bs{q}$ and the $x$-axis,  $p$ is the linear momentum, and $p_{\theta}$ is the angular momentum.
With the new coordinate variables, the Hamiltonian reduces to
\beq\label{eq:Hamiltonian}
H=\frac{1}{4}|\bs{d}|^2\left(G_{ b-1}(0) + G_{ b-1}(r)\right) +\left(p^2+\frac{p_{\theta}^2}{r^2}\right)\left(G_{ b-1}(0)-G_{ b-1}(r)\right).
\eeq
If we treat $p$ as a function of $r$ and every other variables as parameters, then 
\beq\label{eq:p(r)}
p(r) =\pm \sqrt{\frac{H-\ds\frac{1}{4}|\bs{d}|^2\left(G_{ b-1}(0) + G_{ b-1}(r)\right)}{G_{ b-1}(0)-G_{ b-1}(r)} -\ds\frac{p_{\theta}^2}{r^2}} \, . 
\eeq
One can plot the linear momentum $p(r)$ versus $r$ for some fixed values of $|\bs{d}|$ and $p_{\theta}$ as the phase portraits for two-particle dynamics.  We consider these steps for the special case $ b=3/2$. We first compute the Hamiltonians at $r=8$ for various $p$'s for some fixed values of $|\bs{d}|$ and $p_{\theta}$ by using equation (\ref{eq:Hamiltonian}). With these values, we then compute the function $p(r)$ through equation (\ref{eq:p(r)}) for $r_0\le r\le 8$, where $r_0$ is chosen so that the second component of the vector $\bs{p}$ is $p^2>0$. Finally, we plot $p(r)$ versus $r$ as the phase portrait for the fixed values of $|\bs{d}|$ and $p_{\theta}$. Figure \ref{fig:phase_portrait1} is the phase portrait for $|\bs{d}|=1$, $p_{\theta}=0$. Three main behaviors are exhibited in the graph. The ejection and capture orbits are in the upper
and lower-half plane, respectively. The scattering orbits are in the middle. These
orbits correspond to particle collisions when a particle with larger momentum collides with and overcomes one of smaller momentum. 

%This phase portrait is similar to the case for $ b=2$ shown in the paper by McLachlan and Marsland \cite{NP}. 

\begin{figure}[tbh]
\begin{center}
(a)\includegraphics[width=2.75in]{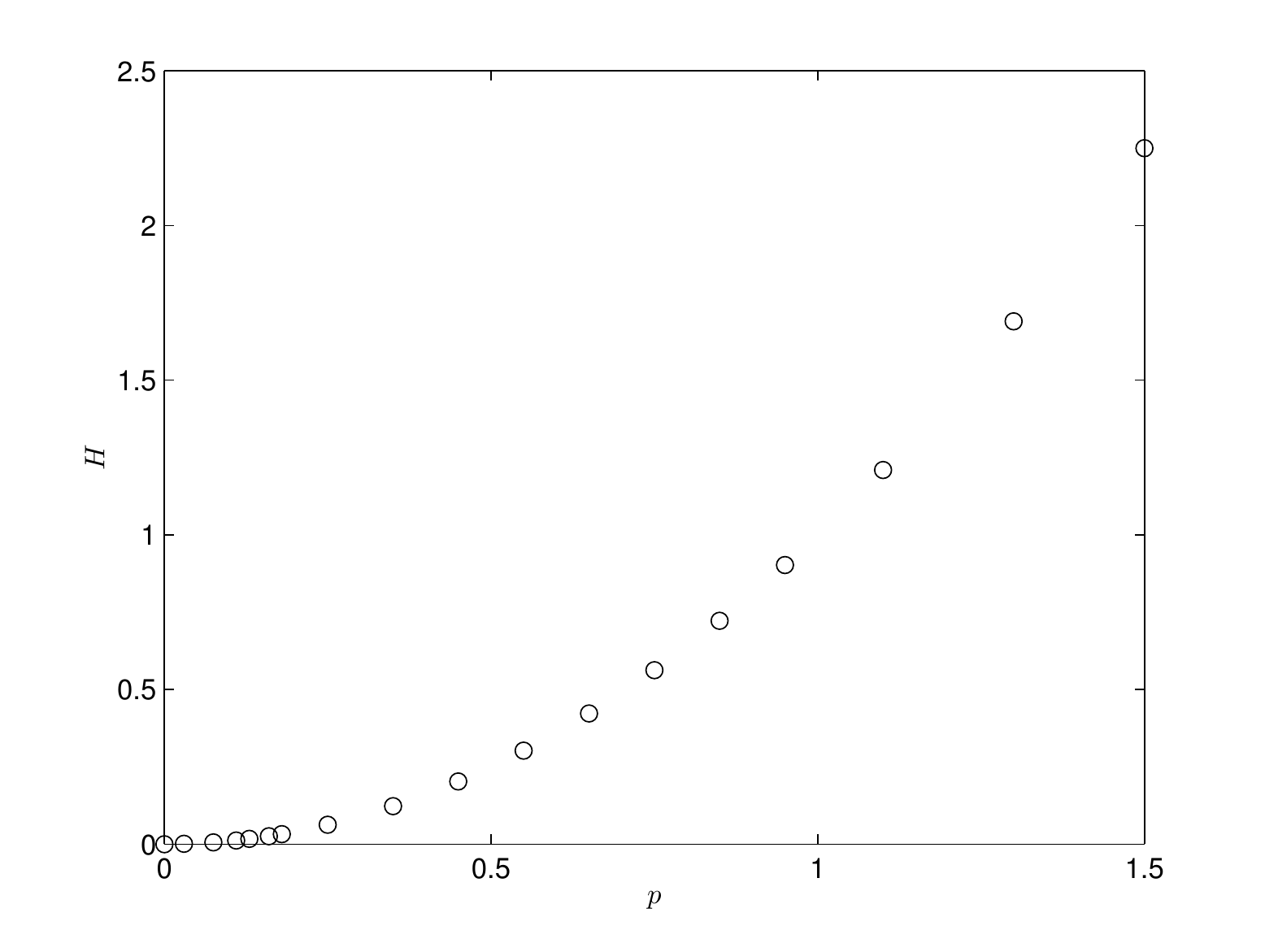}
(b)\includegraphics[width=2.75in]{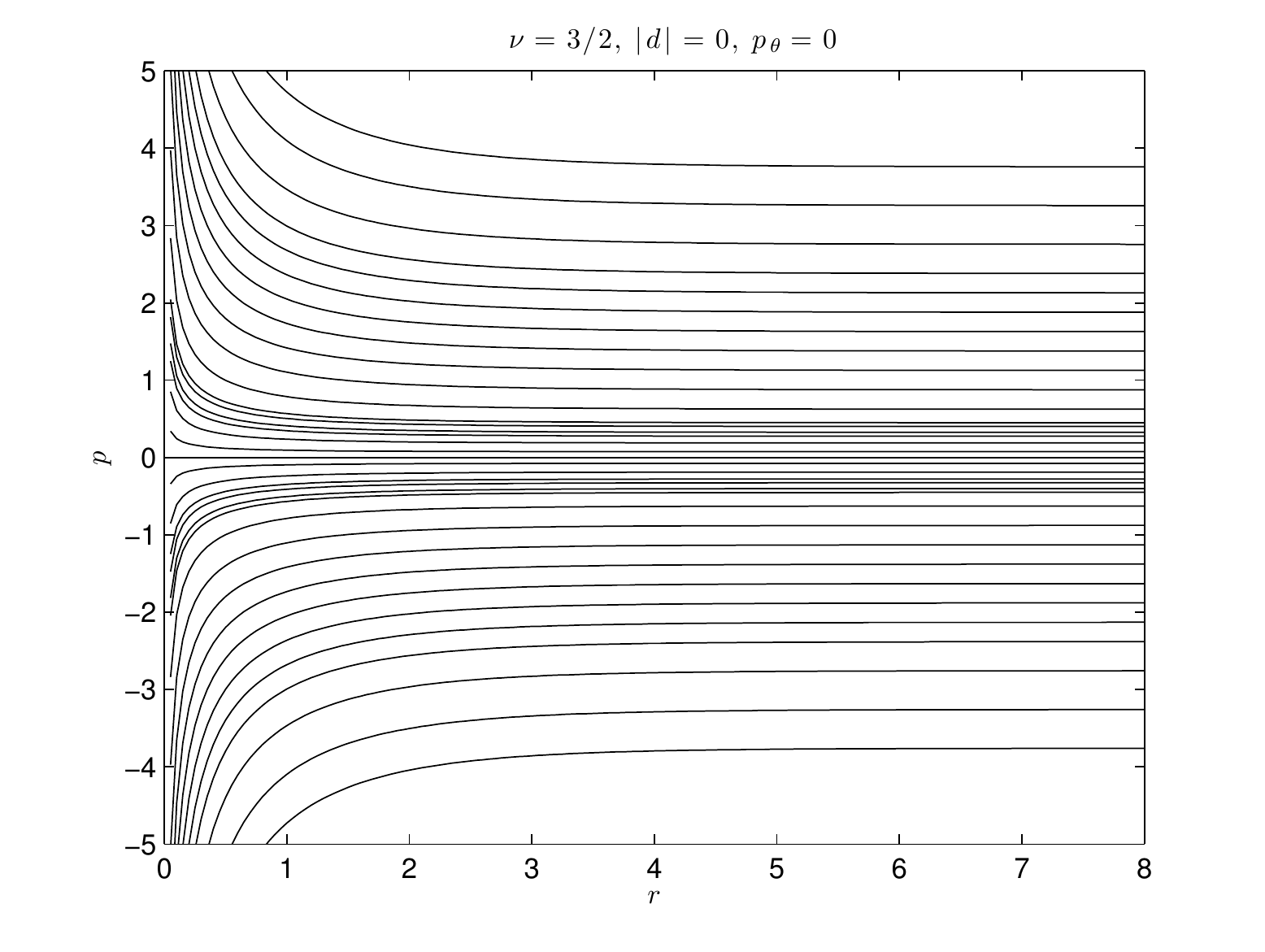}
\end{center}
\caption{Phase portrait of particle-antiparticle head-on collision ($|\bs{d}|=0$ and $p_{\theta}=0$) for the conon case $b=3/2$. (a) The Hamiltonian values computed with the given $|\bs{d}|$, $p_{\theta}$, and $p$ by using equation (\ref{eq:Hamiltonian}) at $r=8$. (b) The graph shows that when particles get closer, their relative linear momenta increase dramatically.}
\label{fig:phase_portrait2}
\end{figure}

Figure \ref{fig:phase_portrait2} is the phase portrait for the particle-antiparticle head-on collision ($|\bs{d}|=0$ and $p_{\theta}=0$). The graph shows that when particles get closer their relative linear momentum increases dramatically. There is, however, no information revealed in the phase portrait about what happens to the linear momentum when $r\ge 0$. We note that the particle motion in Figure \ref{fig:phase_portrait1}-\ref{fig:phase_portrait2} is confined to a line due to the zero angular momentum. Moreover, the scattering orbits in Figure \ref{fig:phase_portrait1} suggest that the relative linear momentum $p$ changes sign at $p=0$ for the case when the sum of linear momenta is non-zero, whereas the lack of scattering orbits in Figure \ref{fig:phase_portrait2} indicates that the relative linear momentum $p$ can only change sign passing through infinity in the particle-antiparticle head-on collision case (the sum of linear momenta is zero).

This behavior is similar to that exhibited by the one-dimensional SW equation, where the Hamiltonian is not conserved when the support of a peakon and an antipeakon coincide in a head-on collision. This leads to divergence of the momenta in the limit to the collision time~\cite{chh94}. Similarly, for the $N$-particle system in this paper,  equation (\ref{eq:N=2-H}) suggests that in a particle-antiparticle head-on collision, when the peaks overlap, the Hamiltonian becomes
\beq
H=\frac{1}{2}|\bs{d}|G_{ b-1}(0).
\eeq
Since $\bs{\dot{d}} =0$, we have $H=0$ if $\bs{d}$ is zero initially.  This would lead to blow-up of the linear momentum  as $r\rightarrow 0$.  As mentioned in the beginning of this section, the continuation of the solutions can be achieved by imposing a conservation a law such as that of the Hamiltonian. The phase portrait analysis for $b=2$, we refer readers to the results in reference \cite{NP}.

\subsection{Two-particle dynamics for the reduced systems}

We continue our study of two-particle dynamics, but focus on the reduced systems (the motion is restricted in the $x$-axis) in this section. For the particle system (see Eqs. (\ref{eq:Q}) \& (\ref{eq:P})), if we let $\bs{Q}^2=0$ and $\bs{P}^2=0$ (the particles are restricted in the $x$-axis and the initial momenta in the $y$-direction is zero), then we obtain
\beq\label{eq:integrable}
\begin{split}
\frac{d\bs{Q}^1}{dt} & = \bs{A}\bs{P}^1\\
\frac{d\bs{Q}^2}{dt} & = \bs{0}\\
\frac{d\bs{P}^1}{dt} & = -\bs{I}^{P^1}\bs{B}^{1}\bs{P}^{1}\\
\frac{d\bs{P}^2}{dt} &= \bs{0}.
\end{split}
\eeq
If $b=3/2$, the above reduced system is the SW equation in two-dimensional space for an arbitrary number of $N$ particles, and hence is completely integrable. (More discussion in Section \ref{sec:6}).

For the rest of this section, we consider $N=2$ for two-particle dynamics.  We will discuss the $N>2$ case in Section \ref{sec:6}.  In this section, we mainly investigate the cases, $b=3/2$ and $b=2$, and we will comment about the case when $b>2$.  

\begin{figure}[tbh]
\begin{center}
(a)\includegraphics[width=2.35in]{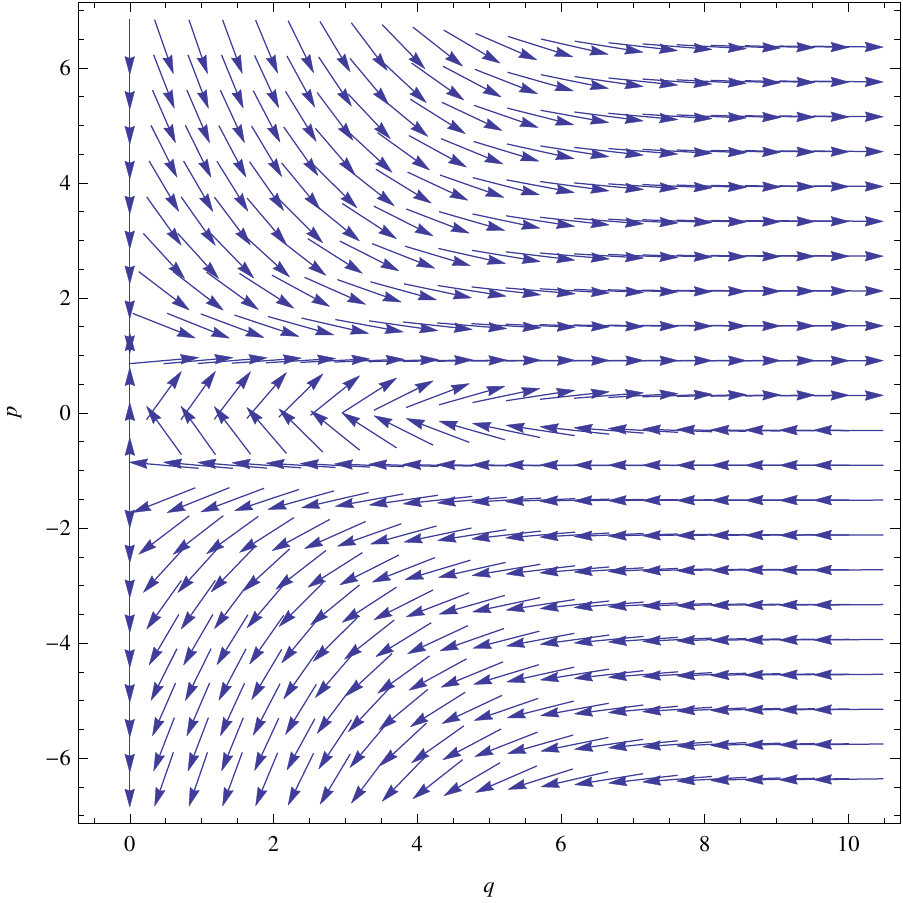}\hskip 24pt
(b)\includegraphics[width=2.35in]{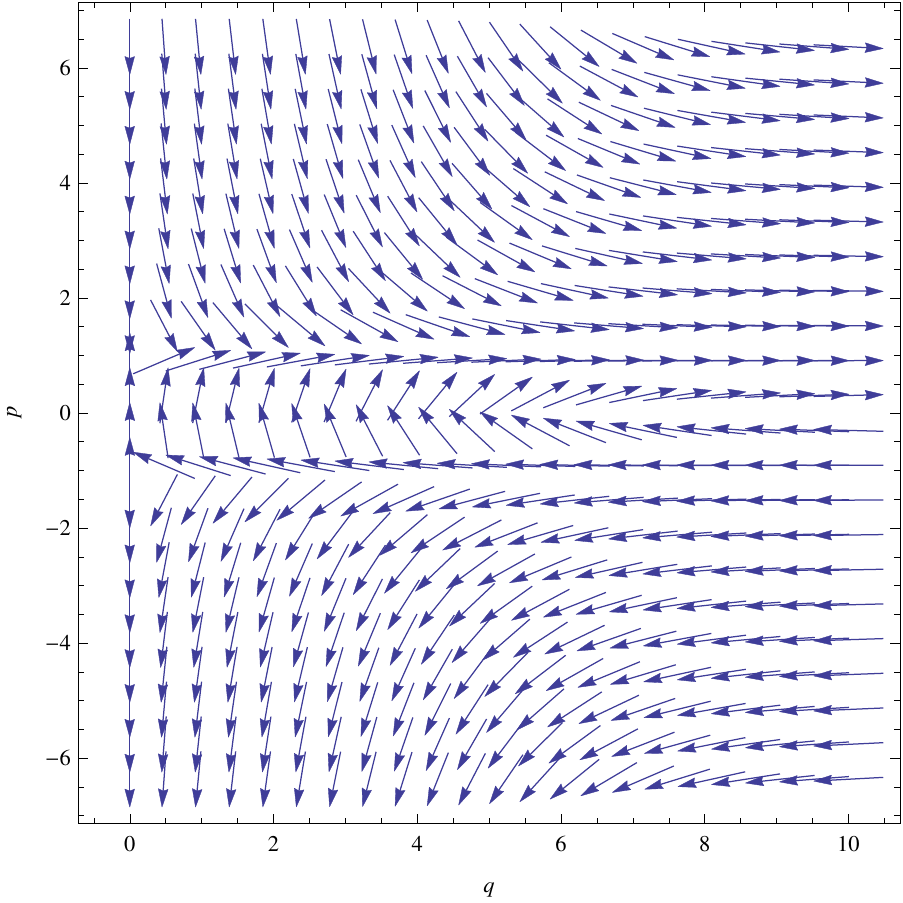}
\end{center}
\caption{The direction field of equation (\ref{eq:ODE-diff}) with $P=1$.  (a) $ b=3/2$. (b) $ b=2$. The graphs clearly show the ejection, capture, and scattering orbits for both $ b=3/2$ and $ b=2$.}
\label{fig:phase_portrait_reduced1}
\end{figure}

Suppose that for $N=2$, two particles are well separated initially (e.g. the distance between the particles approaches infinity), and travel at speeds $c_1$ and $c_2$, respectively, along the $x$-axis. The corresponding traveling waves can be represented by a reduced normalized $N$-particle system, using equations (\ref{eq:new-Q}) and  (\ref{eq:new-P}) (dropping tilde~$\tilde{\cdot}$ notation), for which the second component of the momentum and   position variables is zero, i.e.,
\beq
\bs{p}_1 = [p_1, 0]^{T},\quad  \bs{p}_2=[p_2,0]^{T},\quad\text{and}\quad \bs{q}_1=[q_1, 0]^{T},\quad \bs{q}_2=[q_2,0]^{T}.
\eeq 
Because the motion is confined to a line ($x$-axis), the problem reduces to  one-dimensional dynamics. Only the first component of the two-particle system governs the motion. The system of ODEs for the first component of the two-particle system is 
\beq\label{eq:ODE-1st}
\begin{split}
\frac{d q_1}{dt} &= G_{ b-1}(0)p_1+G_{ b-1}(|q_1-q_2|)p_2,\\
\frac{d q_2}{dt} &= G_{ b-1}(|q_2-q_1|)p_1+G_{ b-1}(0)p_2,\\
\frac{d p_1}{dt} &= -p_1p_2G'_{ b-1}(|q_1-q_2|)\frac{q_1-q_2}{|q_1-q_2|},\\
\frac{d p_2}{dt} &= -p_1p_2G'_{ b-1}(|q_2-q_1|)\frac{q_2-q_1}{|q_2-q_1|}.\
\end{split}
\eeq
Introducing the sum and difference variables 
\beq\label{eq:sum-diff}
\begin{split}
P=p_1+p_2,\quad Q=q_1+q_2\\
p=p_1-p_2,\quad  q=q_1-q_2,
\end{split}
\eeq
we obtain a system of ODEs for $q$ and $p$
\begin{equation}\label{eq:ODE-diff}
\begin{split}
\dot{q}& =\left(G_{ b-1}(0)-G_{ b-1}(|q|)\right)p;\\
\dot{p} & =\frac{p^2-P^2}{2}G'_{ b-1}(|q|).
\end{split}
\end{equation}
The direction field of equation (\ref{eq:ODE-diff}) with $P=1$ is shown in Figure \ref{fig:phase_portrait_reduced1}, where (a) is the case $b=3/2$, and (b) corresponds to $b=2$. It can be seen that the phase dynamics described in Figure \ref{fig:phase_portrait_reduced1}(a) is the same as that of Figure \ref{fig:phase_portrait2}(b) for the same setup and parameters. Figure \ref{fig:phase_portrait_reduced1}, however, clearly shows the ejection, capture, and scattering behaviors indicated in \cite{NP}.

\begin{figure}[tbh]
\begin{center}
(a)\includegraphics[width=2.35in]{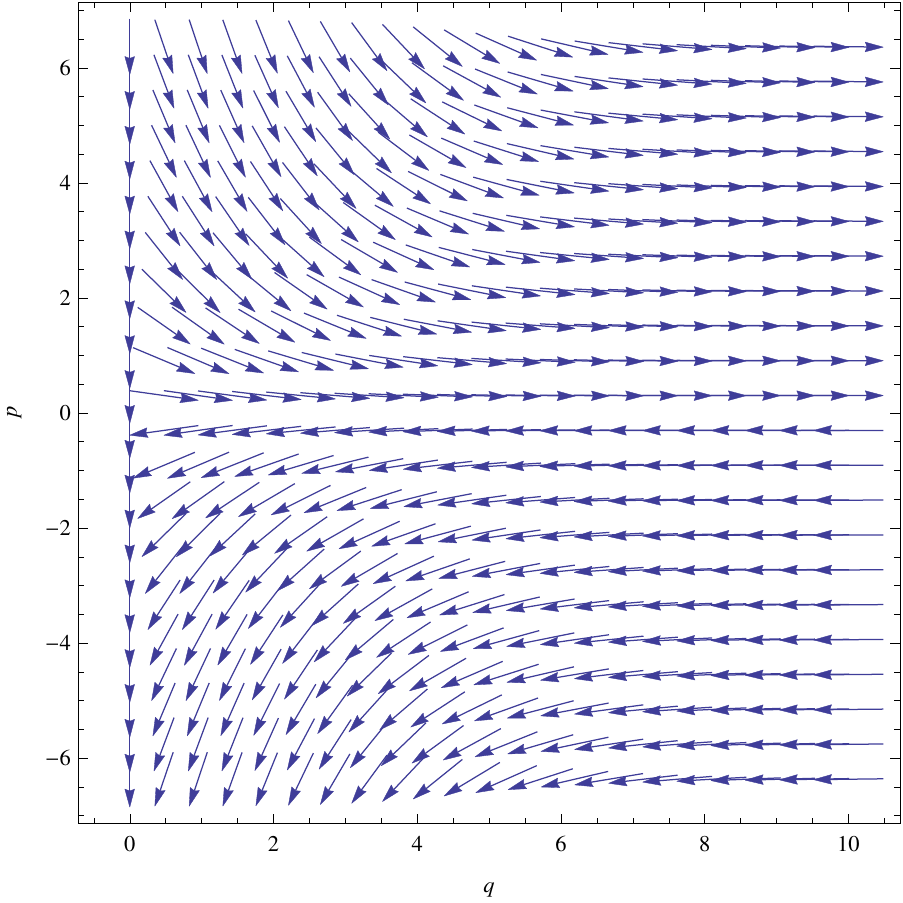}\hskip 24pt
(b)\includegraphics[width=2.35in]{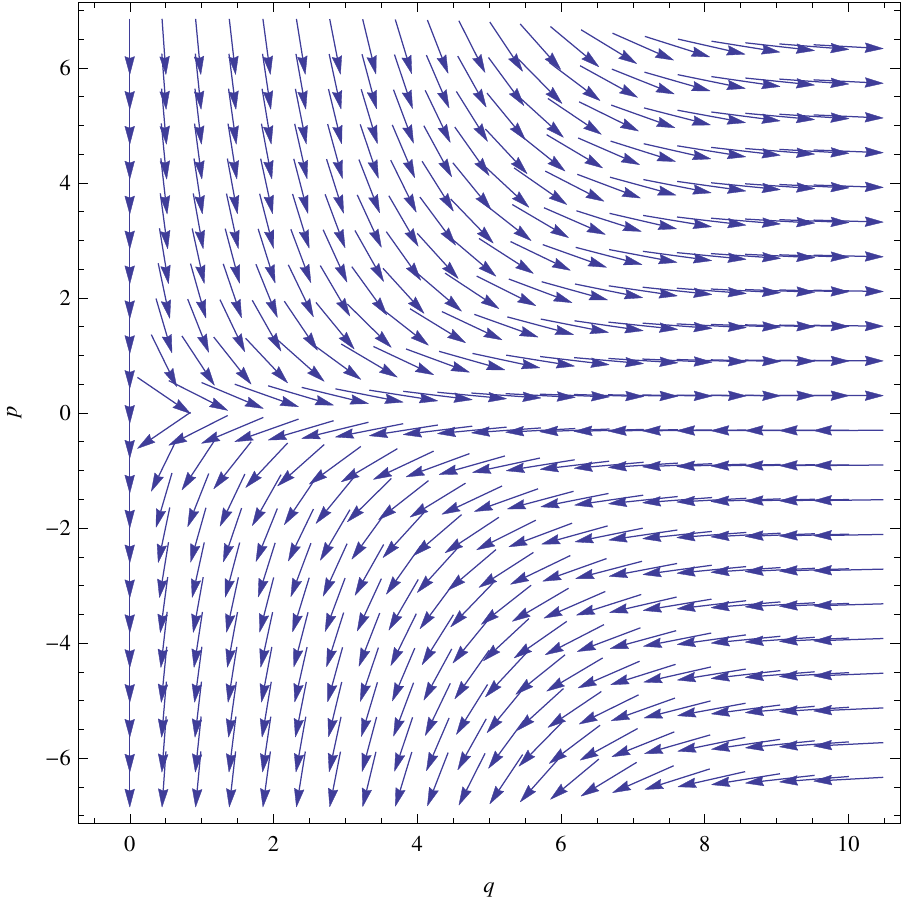}
\end{center}
\caption{The direction field of equation (\ref{eq:ODE-diff}) with $P=0$ (head-on collision).  (a) $ b=3/2$. (b) $ b=2$. The graphs show that there are no scattering orbits, only  ejection and capture orbits for both $ b=3/2$ and $ b=2$.}
\label{fig:phase_portrait_reduced2}
\end{figure}

Next, we consider the head-on collision case, for which $P=0$. The direction fields shown in Figure \ref{fig:phase_portrait_reduced2} indicate that there are no scattering orbits, only ejection and capture orbits exist for head-on collision for both $b=3/2$ and $b=2$ cases. This is consistent with the phase portrait in Figure \ref{fig:phase_portrait1}(b) and those illustrated in \cite{NP}.

As mentioned earlier, the lack of of scattering orbits implies that $p$ can only change sign through infinity in the case of particle-antiparticle head-on collision. 
%The blow-up of $p$ at $q=0$ prevents us from gaining information {\it after} the collision by using  Equation (\ref{eq:ODE-diff}).  
To investigate further the dynamics of head-on collisions, we recall the Hamiltonian (\ref{eq:N=2-H}) for the motion of two particles confined to a line
\beq\label{eq:H_collision}
H=\frac{1}{4}\left(P^2-p^2\right)G_{ b-1}(q) + \frac{1}{4}\left(P^2+p^2\right)G_{ b-1}(0).
\eeq
For a particle-antiparticle head-on collision, $P=0$ and we have
\beq\label{eq:stable-1}
\left(G_{ b-1}(0)-G_{ b-1}(q)\right)= \frac{4H}{p^2}.
\eeq
Using the above relation allows to rewrite equation (\ref{eq:ODE-diff}) as
\begin{eqnarray}
\dot{q}&= &4Hz,\label{eq:q-z-1}\\
\dot{z}&= &-\frac{1}{2}G_{ b-1}'(q)\label{eq:q-z-2},
\end{eqnarray}
where $z=\ds\frac{1}{p}$. Let $\bs{Y}=[q, z]^{T}$. The above equations represent a nonlinear autonomous system
%\beq\label{eq:ODE-Y}
$\dot{\bs Y} = F(\bs{Y})$.
%\eeq
%For $ b\ge 3$, $G_{ b-1}\in C^2$ and the function $F\in C^{1}$ everywhere. Hence by the fundamental existence-uniqueness theorem, there exist a unique solution $\bs{Y}(t)$ on some interval $[-a , a]$ \cite{bib:perko}. 
Since the radial derivative of $G_{ b-1}$ vanishes at $q=0$, due to the symmetry, $\bs{Y}(\bs{0})$ is a fixed-point. 
%The Jacobian matrix for the system (\ref{eq:q-z-1}) \& (\ref{eq:q-z-2}) at the fix-point is 
%\beq
%J=\left(\begin{array}{cccc}
%0 & 4H\\
%-\ds\frac{1}{2}G_{ b-1}''(0) & 0
%\end{array}\right).
%\eeq
%Because $G''_{ b-1}(0) < 0$ for $ b \ge 3$ and $H>0$, 
For smooth enough particles ($G_{ b-1}\in C^2$, $ b\ge 3$), it is straightforward to show that the eigenvalues of the Jacobian matrix for the system (\ref{eq:q-z-1})-(\ref{eq:q-z-2}), linearized around the fixed-point, are real, 
\beq
\lambda_{1, 2} = \pm \sqrt{-2HG''_{ b-1}(0)}.
\eeq
Hence the fixed-point $\bs{Y}(\bs{0})$ is locally a saddle, since $G''_{ b-1}(0) < 0$ for $ b \ge 3$ and $H>0$.

From the Lyapunov function computed in Appendix \ref{app:L}, we know that for particle-antiparticle head-on collisions with $b\ge 3$,  once the motion of the particles is confined to the $x$-axis, the solution stays on the stable manifold. Hence there are no scattering orbits and the particles capture each other.  We conclude that if the solitary waves are smooth enough ($b\ge 3$), for particle-antiparticle collision, scattering orbits can only exist when the motion of the particles is not confined to a line, or the relative angular momentum is non-zero.

%This conclusion does not hold for non-smooth particles. For example, for $ b=2$, $G'_{ b-1}$ in equation (\ref{eq:q-z-2}) is continuous, but its radial derivative ($G''_{ b-1}$) is discontinuous at $q=0$. This means that the function $F$ in equation (\ref{eq:ODE-Y}) is $C^0$, but not $C^1$ in a neighborhood $E$, containing the fixed-point $Y(\bs{0})$, i.e., in such a neighborhood solutions of equation (\ref{eq:ODE-Y}) exist but need not be unique \cite{bib:perko}. 

\begin{figure}[tbh]
\begin{center}
(a)\includegraphics[width=2.75in]{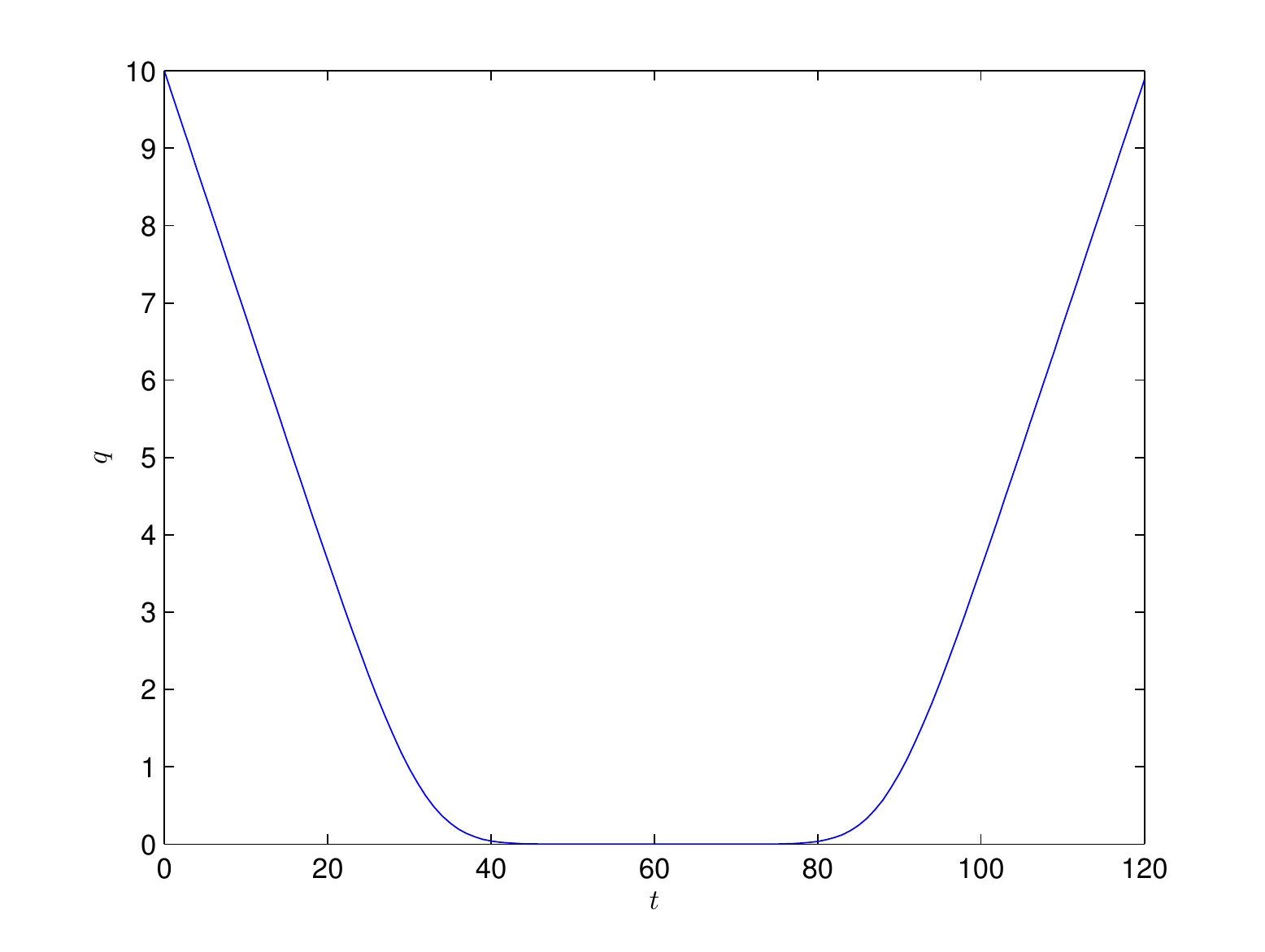}
(b)\includegraphics[width=2.75in]{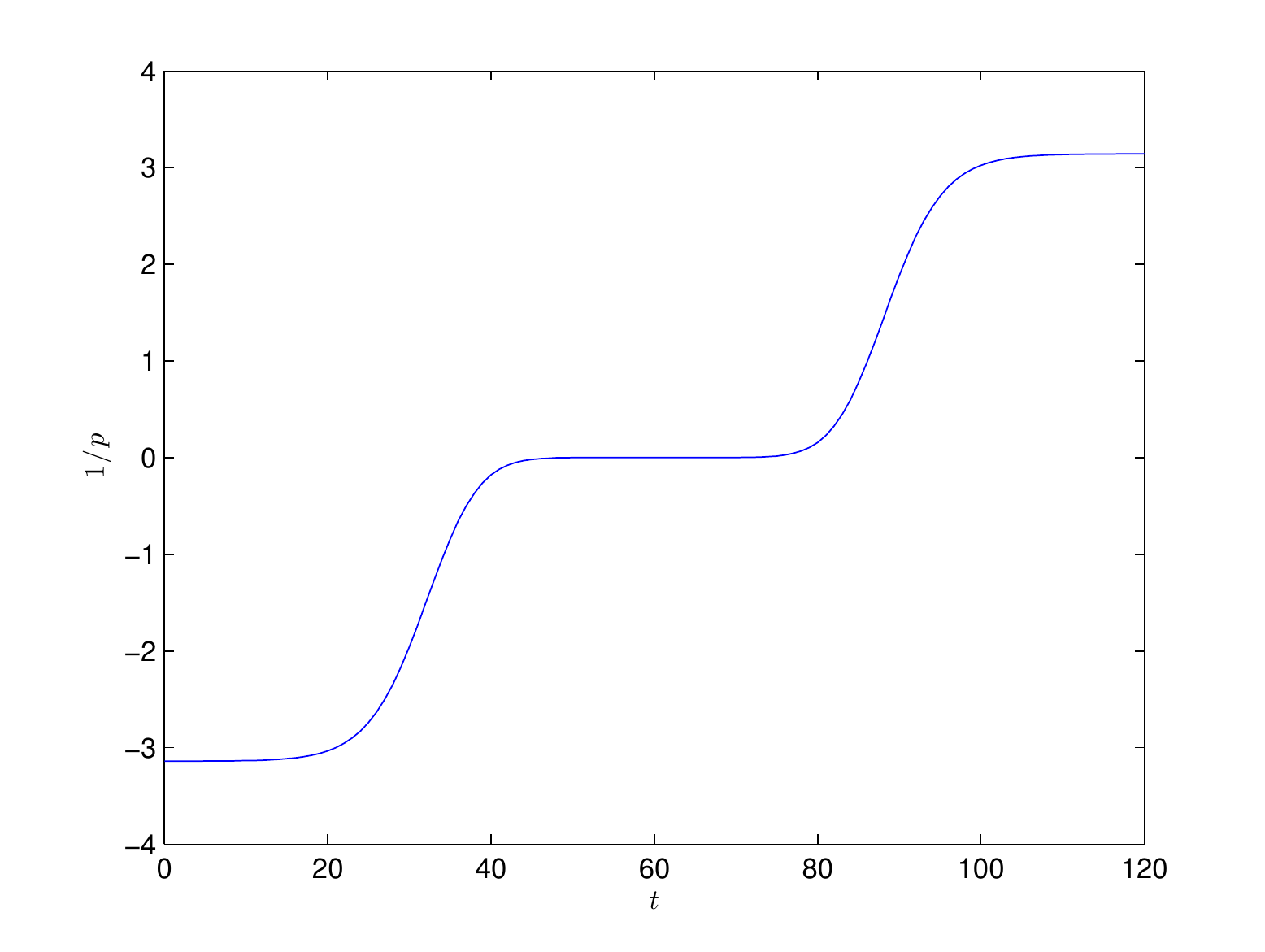}
\end{center}
\caption{A typical scattering solution of two-particle collision for $ b=2$. (a) The plot of $q$ vs $t$. (b) The plot of $z$ vs $t$.}
\label{fig:q-z-reduced-nu-2}
\end{figure}

The property of non-uniqueness may allow scattering solutions for particle-antiparticle collision, even when the motion of particles is restricted to a line. A typical scattering solution is shown in Figure \ref{fig:q-z-reduced-nu-2}.
In the figure, the resonance period ($q=0$) is between $t=40$ and $t=80$. In general, our numerical experiments show that the length of resonance can be arbitrary (due to the non-uniqueness of solutions). Figure~\ref{fig:q-z-reduced-nu-2}  is the numerical integration for the two-particle collision. The figure shows that after the resonance period the particles could exchange momenta as in collisions of two elastic bodies, and move away from each other. Nevertheless, particles are also allowed to keep their momenta, and these solutions allow $q$ to become negative after a resonance period. Figure \ref{fig:q-z-reduced-nu-2} is generated by solving equations (\ref{eq:q-z-1}) and (\ref{eq:q-z-2}) using the sixth-order Runge-Kutta method developed in \cite{bib:butcher64}. The initial conditions are $z(0)=\pi$ and $q(0) = 10$. The time step is $\Delta t=3.125$e-5. 

Finally, we focus on the choice $b=3/2$, for which the derivative of the Green function $G_{ b-1}$ is not continuous at $q=0$. Hence, in a neighborhood of $E$  containing $Y(\bs{0})$, the solution may or may not exist.  One can enforce a continuation rule for solutions of the particle-antiparticle head-on collision. In particular, this rule can be assigned to correspond to elastic collisions. i.e., the particles exchange momenta and scatter after the collision. We take a closer look at such solutions next.

\subsection{Exact solution for $ b=3/2$}

For the special case $ b=3/2$, we write the system of equations (\ref{eq:ODE-1st}) in terms of the sum and difference variables (\ref{eq:sum-diff}) as
\beq\label{eq:ODE-transformed}
\begin{split}
&\dot{P}=0,\qquad\qquad\qquad\qquad\qquad\qquad\quad\dot{Q}=\left(G_{ b-1}(0) + G_{ b-1}(|q| \right)P,\\
&\dot{p}= -\frac{1}{2}\left(P^2-p^2\right)\sgn(q) G'_{ b-1}(|q|), \quad \dot{q}=\left(G_{ b-1}(0) - G_{ b-1}(|q| \right)p,
\end{split}
\eeq
where $\sgn(q)$ is the signum function.
%defined as 
%\beq
%\sgn(q)=\begin{cases}
%+ 1,\quad &q > 0,\\
%-1,\quad & q<0.
%\end{cases}
%\eeq
The second pair in the above equations is the same as equations (\ref{eq:ODE-diff}), since the Green functions for the elliptic equations are evenly symmetric. We consider the case of $ b=3/2$, for which the normalized Green function and its derivative are 
\beq\label{eq:3/2-normalized-green}
G_{1/2}(r) = e^{-r},\quad  G_{1/2}(0) =1,\quad \text{and}\quad G'_{1/2}(r) = - e^{-r}.
\eeq
Thus for this special case  equation(\ref{eq:ODE-transformed}) becomes
\beq\label{eq:ODE-nu-3/2}
\begin{split}
&\dot{P}=0,\qquad\qquad\qquad\qquad\qquad\dot{Q}=\left(1 + e^{-|q|} \right)P,\\
&\dot{p}= \frac{1}{2}\left(P^2-p^2\right)\sgn(q) e^{-|q|},\quad\dot{q}=\left(1 - e^{-|q|} \right)p.
\end{split}
\eeq
The above reduced system is a two-dimensional 2-body collision problem restricted to the $x$-axis. We note that these equations coincide with those for the interaction of two solitons of the one-dimensional SW equation \cite{chh94}. (The exact solution of the above system was derived in \cite{bib:ch93} and~\cite{chh94}.) In Appendix \ref{sec:exact}, we present an example of the exact solution and use this to test the numerical solution of equations (\ref{eq:q-z-1})-(\ref{eq:q-z-2}).

\subsection{Examples of particle interaction}

We present numerical integration for the ODEs system to illustrate two-particle interaction. W focus on the special case $b=3/2$.

%based on the particle system and the analysis developed in the previous sections focusing on two examples. 

\vskip 12 pt

\begin{figure}[h]
\begin{center}
(a)\includegraphics[width=2.75in]{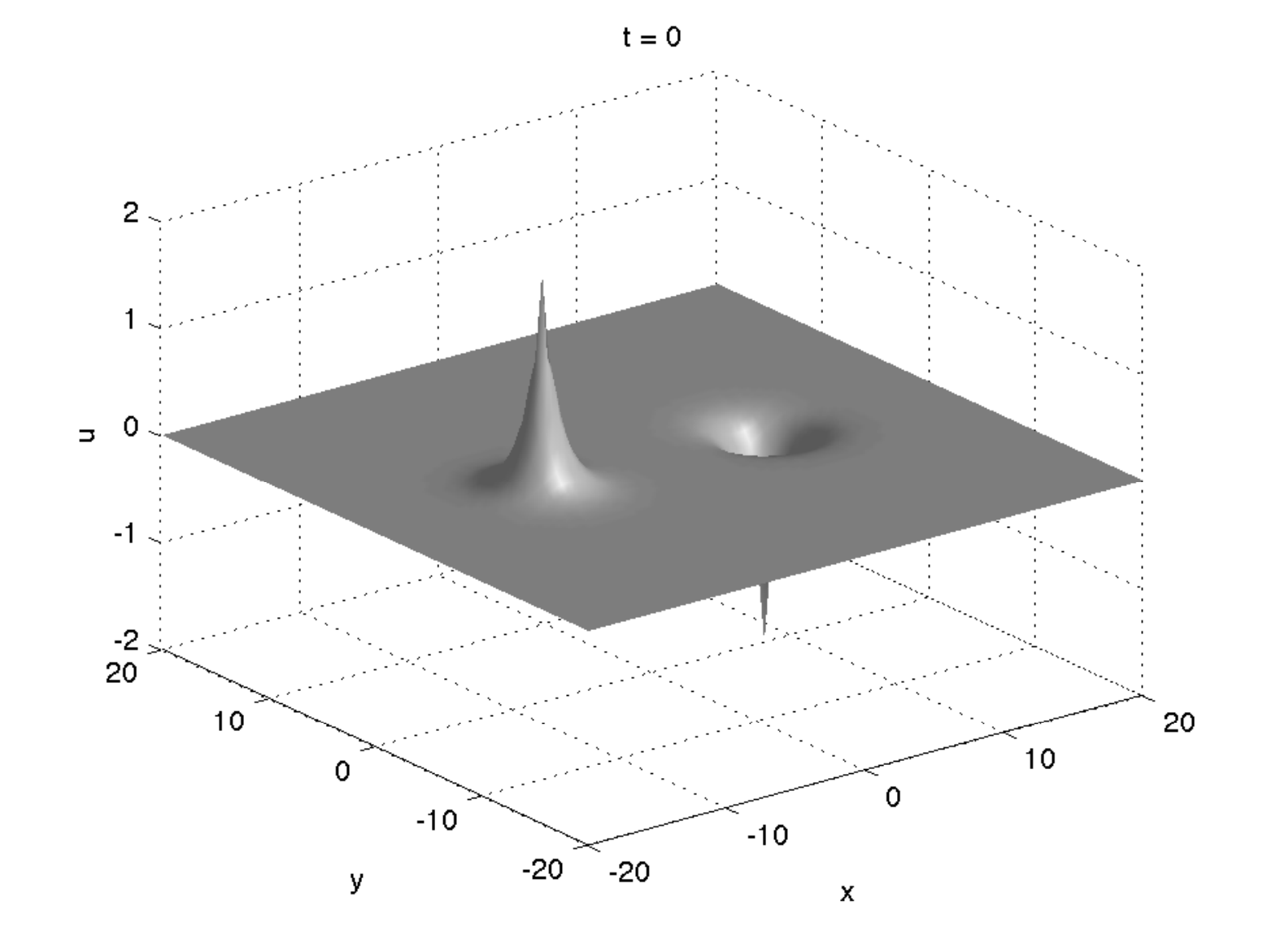}
(b)\includegraphics[width=2.75in]{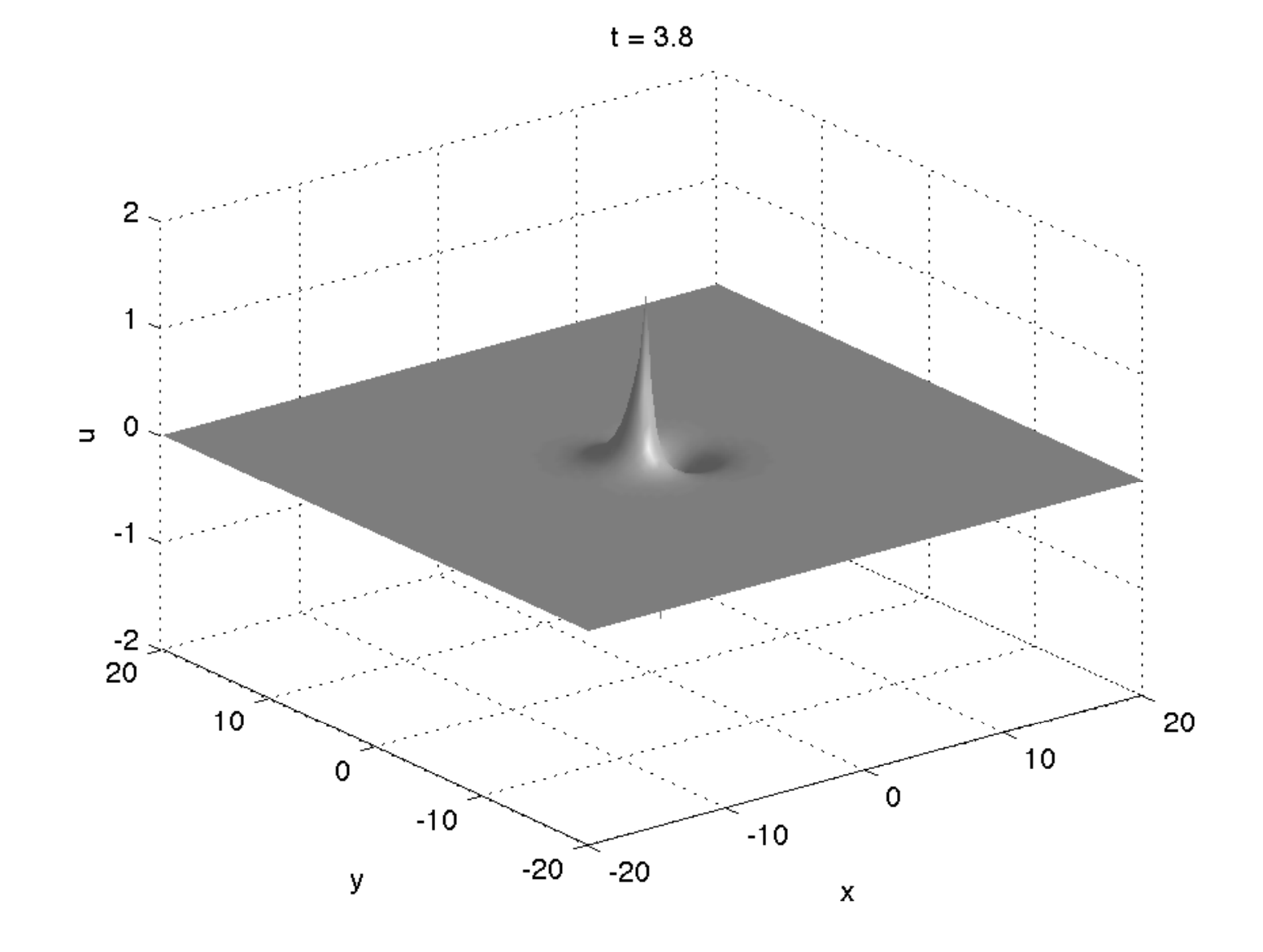}
(c)\includegraphics[width=2.75in]{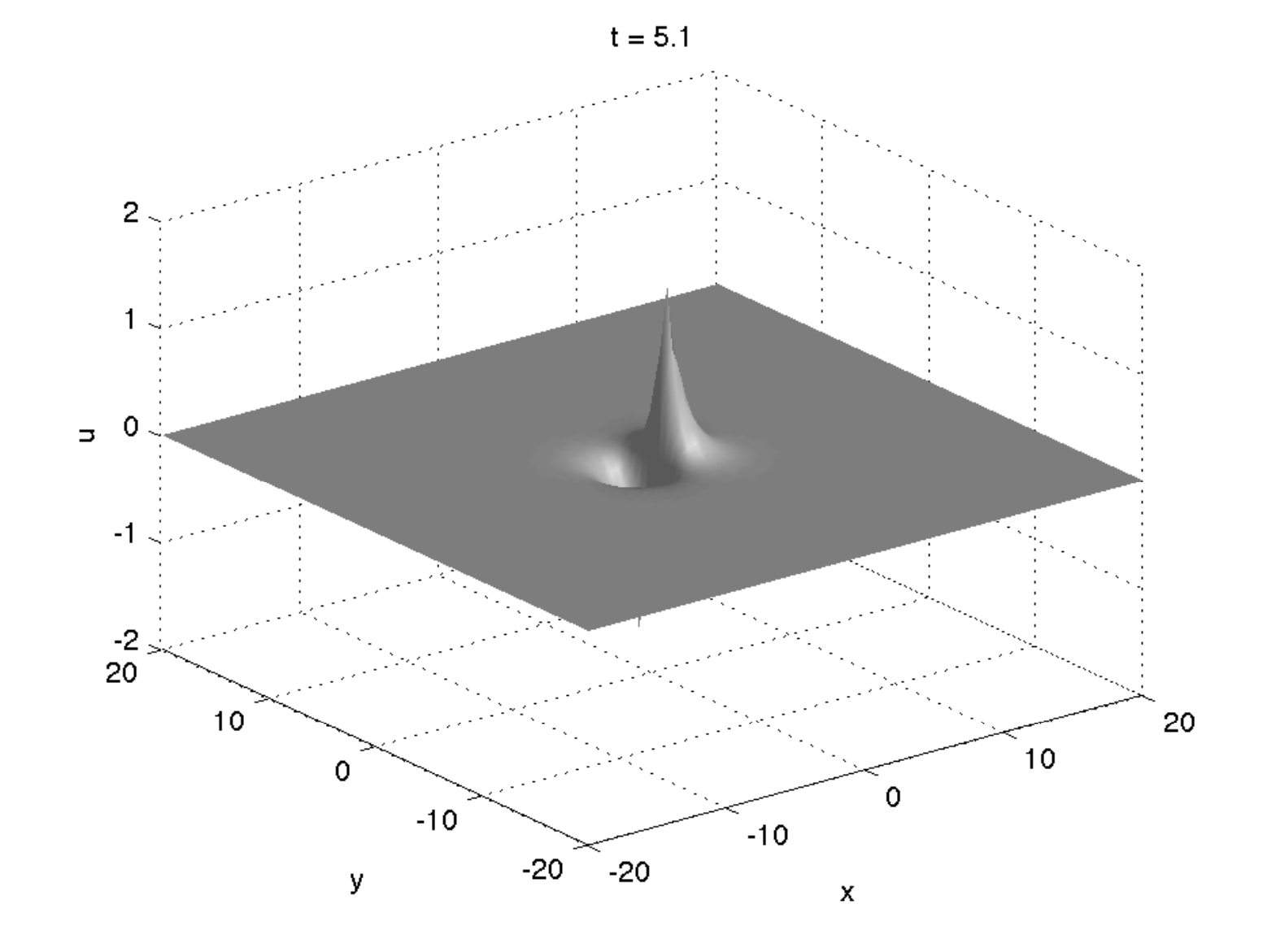}
(d)\includegraphics[width=2.75in]{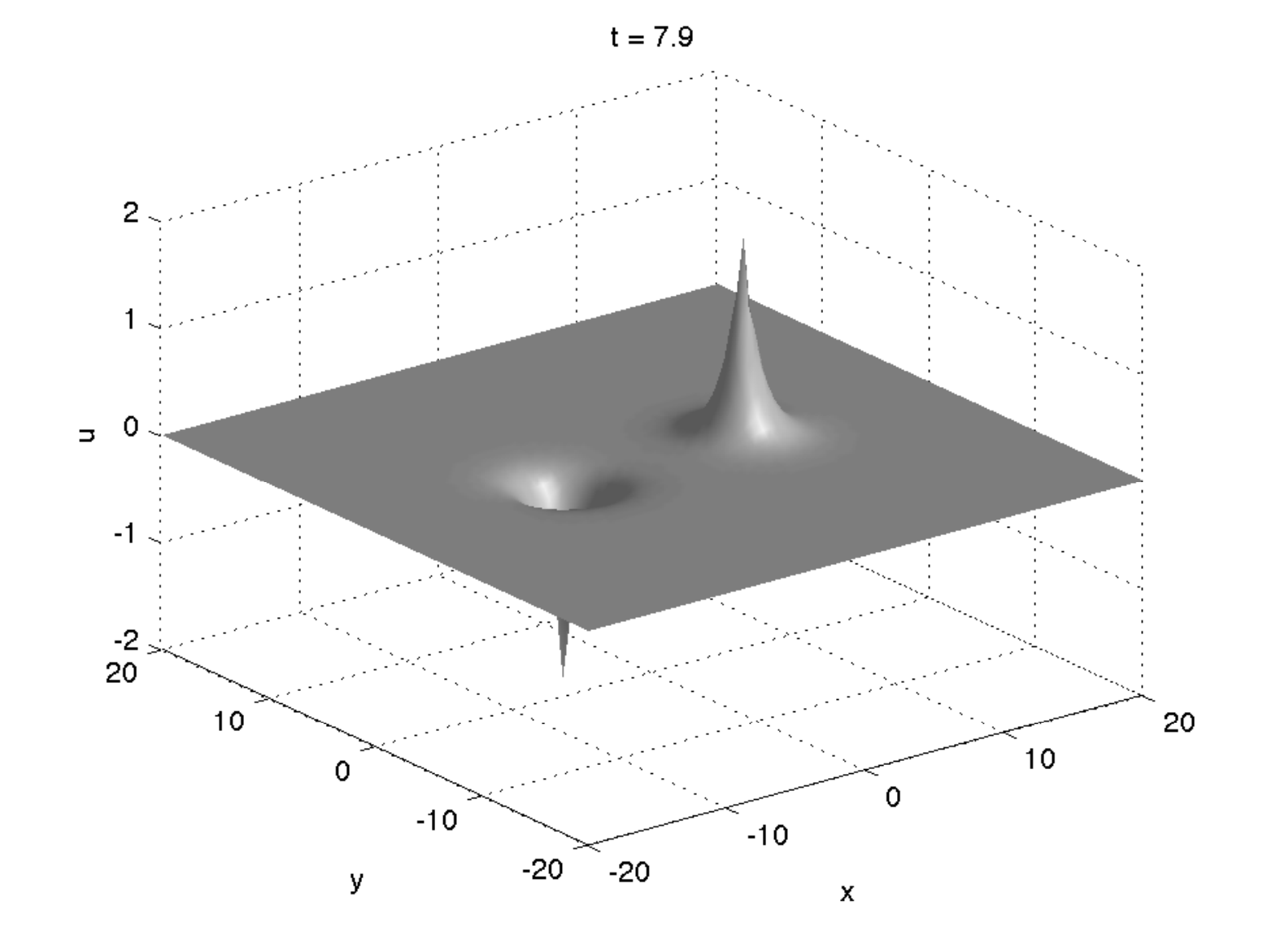}
\end{center}
\caption{Numerical reconstruction of the first-component of $\bs{u}$ for head-on collision between particle and antiparticle for $ b=3/2$. The first component of $\bs{u}$ is shown at (a) $t=0$, (b) $t=3.8$, (c) $t=5.1$, and (d) $t=7.9$, respectively. Initially the particle $\bs{p}_1=[2, 0]^{T}$ is located at $\bs{q}_1=[-8, 0]^{T}$, while the antiparticle  $\bs{p}_2=[-2, 0]^{T}$ is located at $\bs{q}_2=[8, 0]^{T}$. (The particle locations and momenta are computed by using equation (\ref{eq:sol4}), as shown in Figure \ref{fig:comp_exact_RK6}
in Appendix D.}
\label{fig:head-on-32-direct}
\end{figure}

\noindent 
{\bf Example 1:}  We first show the particle-antiparticle head-on collision for $ b=3/2$. The integration of the two-particle system suffers from divergence of the momentum when the two particles collide. Instead, we reconstruct the solution $\bs{u}$ by using the exact solutions of $p$ and $q$ obtained by equation (\ref{eq:sol4}), and the reconstruction formula (\ref{eq:new-anstz}). Suppose that initially the particle $\bs{p}_1=[2, 0]^{T}$ is located at $\bs{q}_1=[-8, 0]^{T}$, while the antiparticle  $\bs{p}_2=[-2, 0]^{T}$ is located at $\bs{q}_2=[8, 0]^{T}$. From equation (\ref{eq:tc}), the two particles collide at $t_c\approx 4.346573576213$. Figure \ref{fig:head-on-32-direct} is the plots of the first component of $\bs{u}$ before and after the collision at (a) $t=0$, (b) $t=3.8$, (c) $t=5.1$ and (d) $t=7.9$, respectively.

%\begin{figure}[h]
%\begin{center}
%(a)\includegraphics[width=2.75in]{nu2_qz_q-eps-converted-to.pdf}
%(b)\includegraphics[width=2.75in]{nu2_qz_z-eps-converted-to.pdf}
%\end{center}
%\caption{The solutions of equations (\ref{eq:q-z-1}) \& (\ref{eq:q-z-2}), with $ b=2$, for the initial data $q(0)=16$ and $z(0)=1/4$ corresponding to the initial conditions $\bs{q}$ and $\bs{p}$ in {\bf Exampl 1}. (a) $q$ vs $t$. (b) $z$ vs $t$.}
%\label{fig:head-on-2-qz}
%\end{figure}

%\begin{figure}[h]
%\begin{center}
%(a)\includegraphics[width=2.75in]{t0_nu_2_head_on-eps-converted-to.pdf}
%(b)\includegraphics[width=2.75in]{t1_nu_2_head_on-eps-converted-to.pdf}
%(c)\includegraphics[width=2.75in]{t2_nu_2_head_on-eps-converted-to.pdf}
%(d)\includegraphics[width=2.75in]{t3_nu_2_head_on-eps-converted-to.pdf}
%\end{center}
%\caption{Direct numerical simulation for head-on collision between soliton and antisoliton for $ b=2$. The first component of the velocity $\bs{u}$ is shown at (a) $t=0$, (b) $t=3.85$, (c) $t=8.7$ and (d) $t=12.2$, respectively.  Initially the soliton $\bs{p}_1=[2, 0]^{T}$ is located at $\bs{q}_1=[-8, 0]^{T}$, while the antisoliton  $\bs{p}_2=[-2, 0]^{T}$ is located at $\bs{q}_2=[8, 0]^{T}$. Later the particle locations and momenta are computed by using equations (\ref{eq:q-z-1}) \& (\ref{eq:q-z-2}), as shown in Figure \ref{fig:head-on-2-qz}. The velocity is reconstructed by  equations (\ref{eq:new-anstz}).}
%\label{fig:head-on-2-direct}
%\end{figure}

\begin{figure}[h]
\begin{center}
(a)\includegraphics[width=2.75in]{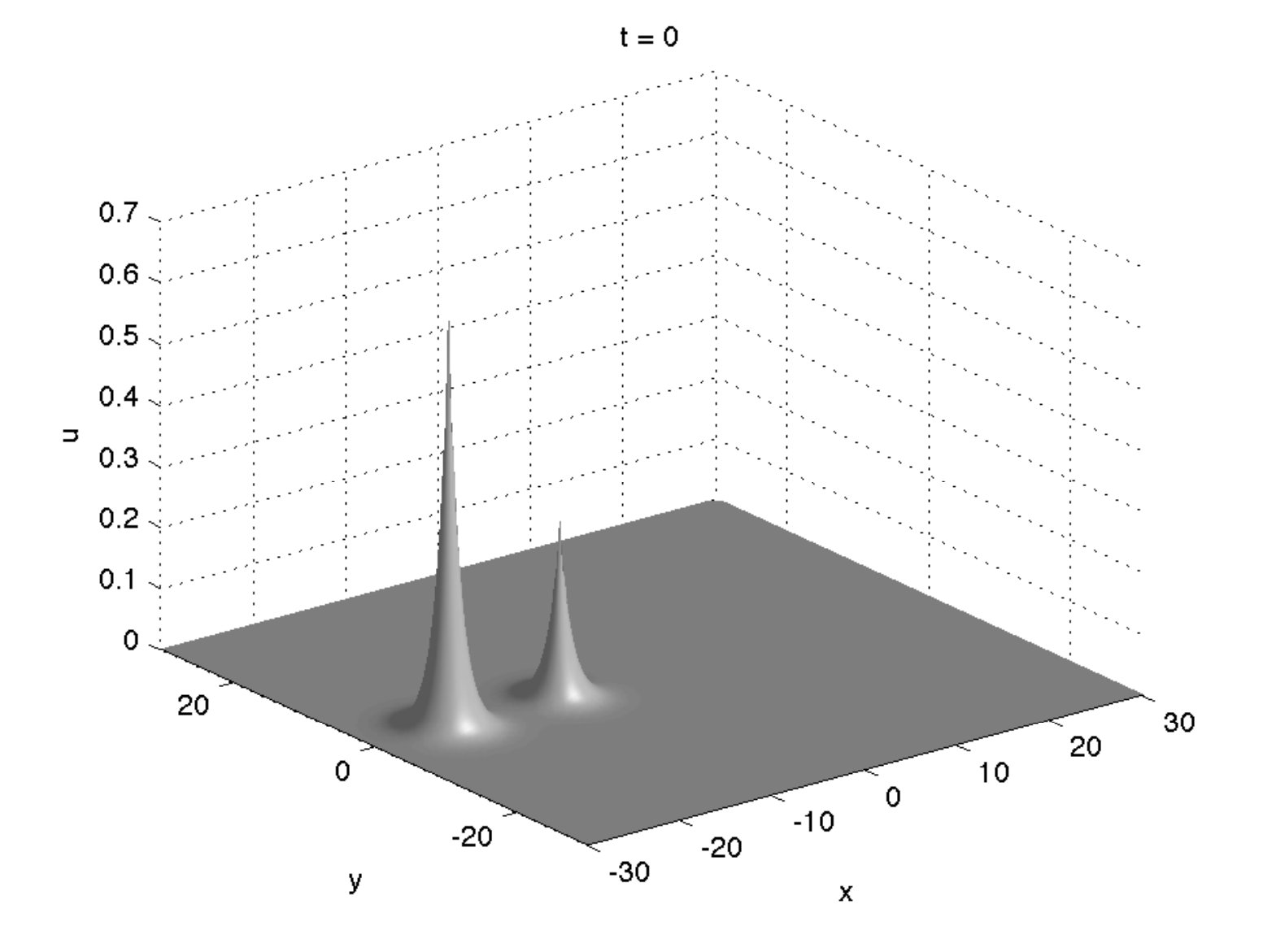}
(b)\includegraphics[width=2.75in]{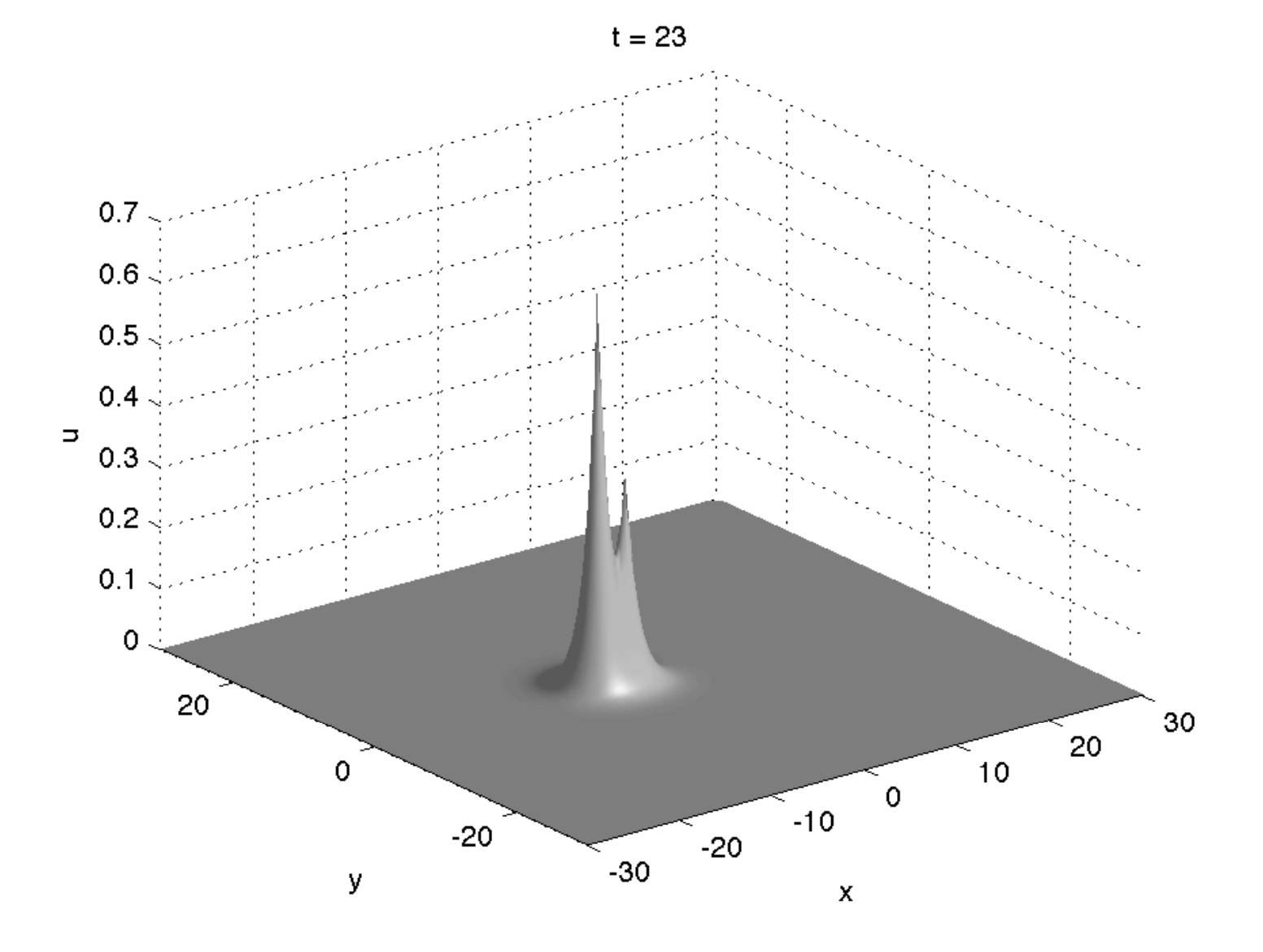}
(c)\includegraphics[width=2.75in]{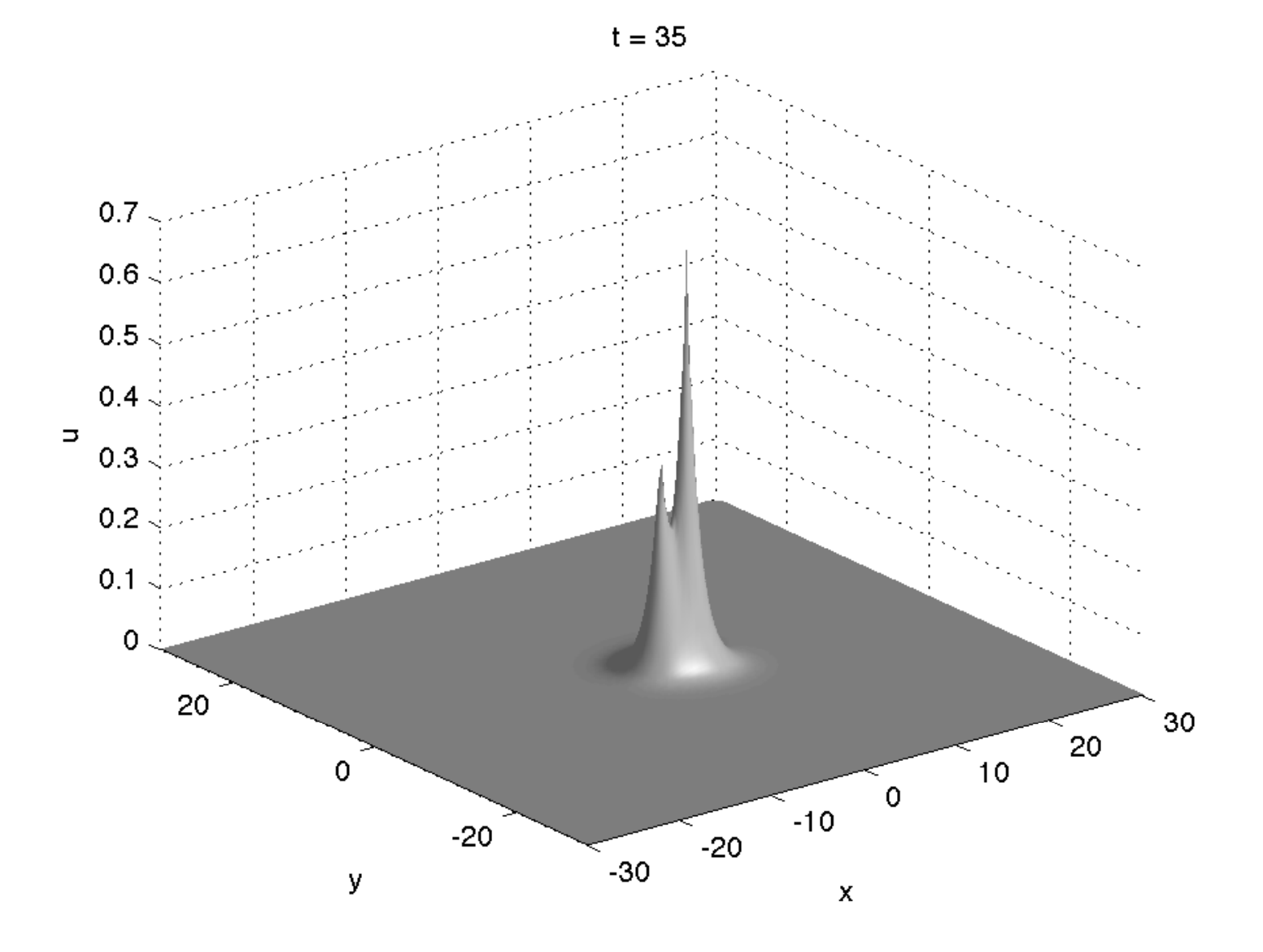}
(d)\includegraphics[width=2.75in]{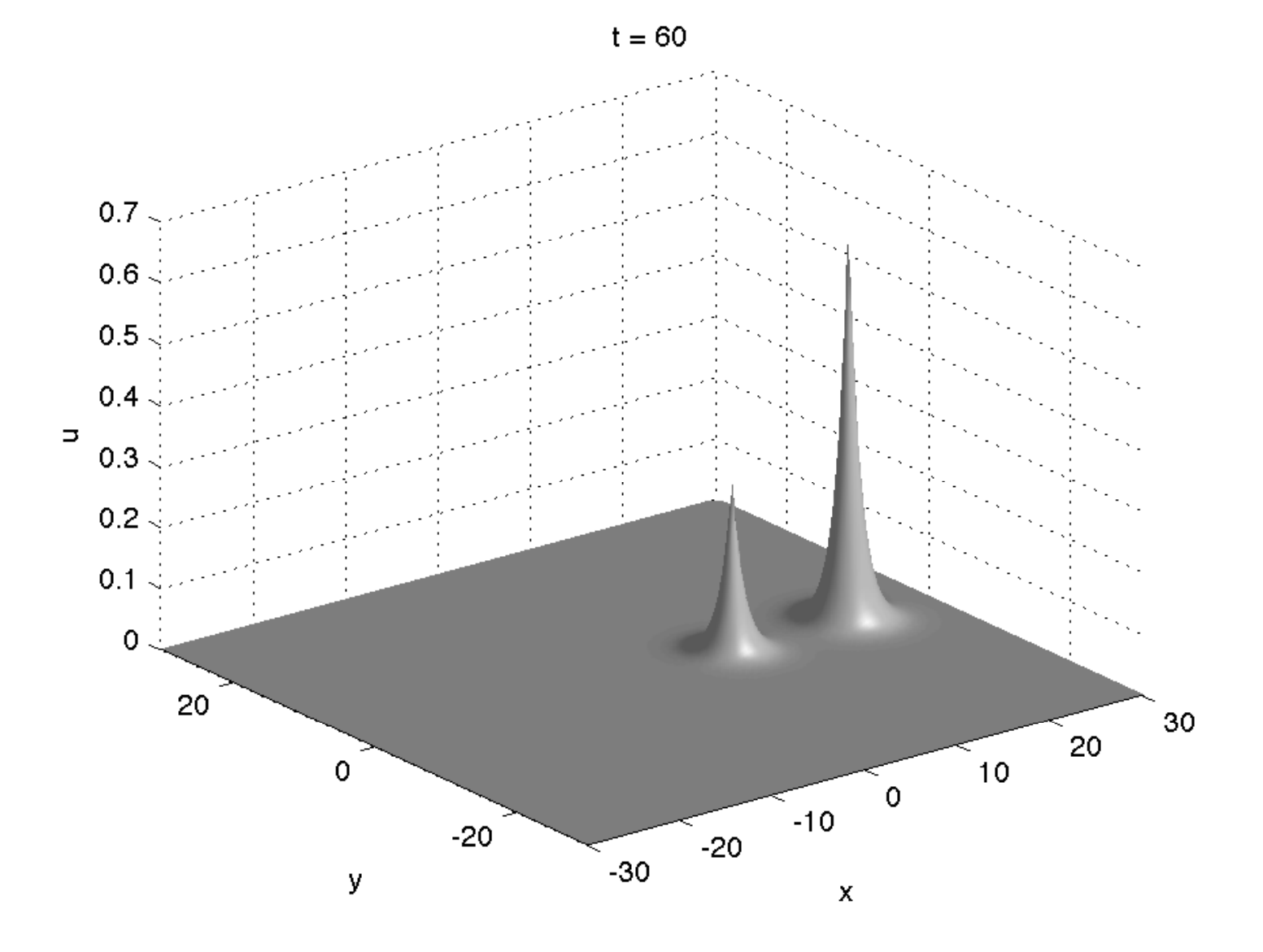}
\end{center} 
\caption{Numerical integration for the two-particle system for $ b=3/2$. Plots show the first component of $\bf{u}$. (a) $t=0$, (b) $t=23$, (c) $t=35$, (d) $t=60$. The fast particle overtakes and then leave the slow one. The initial conditions are $\bs{p}_1=[0.7, 0]^{T}$, $\bs{q}_1=[-22, 0]^{T}$,  $\bs{p}_2=[0.3, 0]^{T}$, and $\bs{q}_2=[-10, 0]^{T}$, respectively.  From Figure \ref{fig:phase_portrait1}, scattering orbits exist for these initial data. Plots are conducted by integrating the $N$-particle algorithm, (\ref{eq:new-Q}) \& (\ref{eq:new-P}). The time step is $\Delta t= 1.0$e-4.}
\label{fig:overtake-32-direct}
\end{figure}

%
%\vskip 12 pt
%
%\noindent 
%{\bf Exampl 2:} In this example, we present the similar simulation as shown in {\bf Exampl 1} for $ b=2$. Since there are no exact solutions for $p$ and $q$ for $ b=2$, we compute $p$ and $q$ by using equations (\ref{eq:q-z-1}) \& (\ref{eq:q-z-2}), for which the derivative of the Green's function is $G_1'(q)=-qK_0(q) $, where $K_0(q)$ is the modified Bessel function of order zero. Figure \ref{fig:head-on-2-qz} is the solutions of equations (\ref{eq:q-z-1}) \& (\ref{eq:q-z-2}) for the initial data $q(0)=16$ and $z(0)=1/4$ corresponding to the initial conditions $\bs{q}$ and $\bs{p}$ in {\bf Exampl 1}. Figure \ref{fig:head-on-2-direct} is the simulation of the first component of the velocity $\bs{u}$ before and after a capture period at (a) $t=0$, (b) $t=3.85$, (c) $t=8.7$ and (d) $t=12.2$, respectively.

\vskip 12pt

\noindent
{\bf Example 2:} We consider the case of  two solitary waves travelling in the same direction. Suppose the solitary wave that has a larger amplitude (momentum) is behind and travels faster than the other one. The fast solitary wave will overcome the slow one, and after separating  two solitary waves will emerge which continue to travel at their original speeds.  From Figure~\ref{fig:phase_portrait1}, for $ b=3/2$, if the sum of momenta is 1, there exist scattering orbits for relative momenta $p\le 0.5$. Because the momentum does not blow up when the two waves collide, the plots are obtained by numerical integration of the $2$-particle sysetm, equations (\ref{eq:new-Q}) \& (\ref{eq:new-P}). We choose the initial data as $\bs{p}_1=[0.7, 0]^{T}$, $\bs{q}_1=[-22, 0]^{T}$,  $\bs{p}_2=[0.3, 0]^{T}$, and $\bs{q}_2=[-10, 0]^{T}$, respectively, so that the sum of momenta is 1 and the relative momentum is $0.4$. The time step for the integration is $\Delta t= 1.0$e-5. Figure \ref{fig:overtake-32-direct} shows the simulation for the waves before, during, and after the overtaking process.

%\section{$N$ particles restricted to the $x$-axis}\label{sec:6}
%
%It is known (see, e.g., \cite{chh94,bib:DCDS03, bib:Camassa05}), that a genetic behavior of the dispersionless SW equation ($\kappa=0$) consists of a train of solitons emerging from the location of an initial hump, whereas for the case ($\kappa\ne 0$), the solution additionally exhibits dispersive behavior with an ensuing wave train. In Section \ref{sec:6.1}, we demonstrate a similar behavior observed for the general model equations with 
% $\bs{\kappa=0}$),  while in Section \ref{sec:6.2}, with the same initial data prescribed along the $x$-axis in the two-dimensional domain,  for the the dispersive equations derived in Section \ref{sec:dispersive}.
%\todo{the preamble to section 6 is not really understandable and is poor English, tried to fix it but seems to be redundant anyway, so commented it out above} 

\section{Initial data $p_1(x) =a\,\sech^2(x)$ for the reduced systems}\label{sec:6}

%\subsection{Dispersionless case ${\bs\kappa}={\bs 0}$}\label{sec:6.1}
\subsection{Dispersionless case ${\bs\kappa}={\bs 0}$}\label{sec:6.1}

In this section, we investigate the reduced system (\ref{eq:integrable}) for $N > 2$ and $ b=3/2$. The normalized Green function is provided by equation (\ref{eq:3/2-normalized-green}). Suppose that we place $N$ particles on the $x$-axis with non-zero first momentum-component while setting the second component of momentum to zero. Suppose all other particles in the domain have zero momenta for both components.  Then the first component of the particle system of equation (\ref{eq:N-particle}) reduces to
\beq\label{eq:2D-qp}
\begin{split}
\dot{q}_i & = \sum_{j=1}^{N}e^{-|q_i-q_j|}p_j,\\
\dot{p}_i & = -\sum_{\substack{j=1\\ j\ne i}}^N\sgn(q_i-q_j) e^{-|q_i-q_j|}p_i p_j,
\end{split}
\eeq
where $i=1,\cdots, N$. Here $q_i$ is the first coordinate of the $i^{th}$ particle on the $x$-axis, while $p_i$ is its momentum. Equation (\ref{eq:2D-qp}) is a completely integrable system. It shares the same form and hence the same properties as that of the completely integrable $N$-particle system for the one-dimensional (1-D) SW equation studied in \cite{bib:DCDS03, bib:chl06,bib:cl08}. To illustrate that in two-dimensional space the solution behaviors, based on the setup and system (\ref{eq:2D-qp}), are the same as those of the 1-D SW equation, we consider the following initial data. A $100\times 100$ particle grid is placed in a domain of size $[-20, 20]\times[-20, 20]$. Along the $x$-axis, the first-component momentum for the particles is given by
%\beq\label{eq:sech2_x-axis}
$p(x)=\frac{1}{2}\sech^2(x)$.
%\eeq
The first momentum-component  is zero outside the $x$-axis, whereas the second momentum-component  is zero everywhere. This initial condition is chosen to emulate the sharp traveling wave solution of the SW equation. The initial wave hump sharpens as it moves to the right, followed by others emerging humps from the initial condition support. 

\begin{figure}[tbh]
\begin{center}
(a)\includegraphics[width=2.75in]{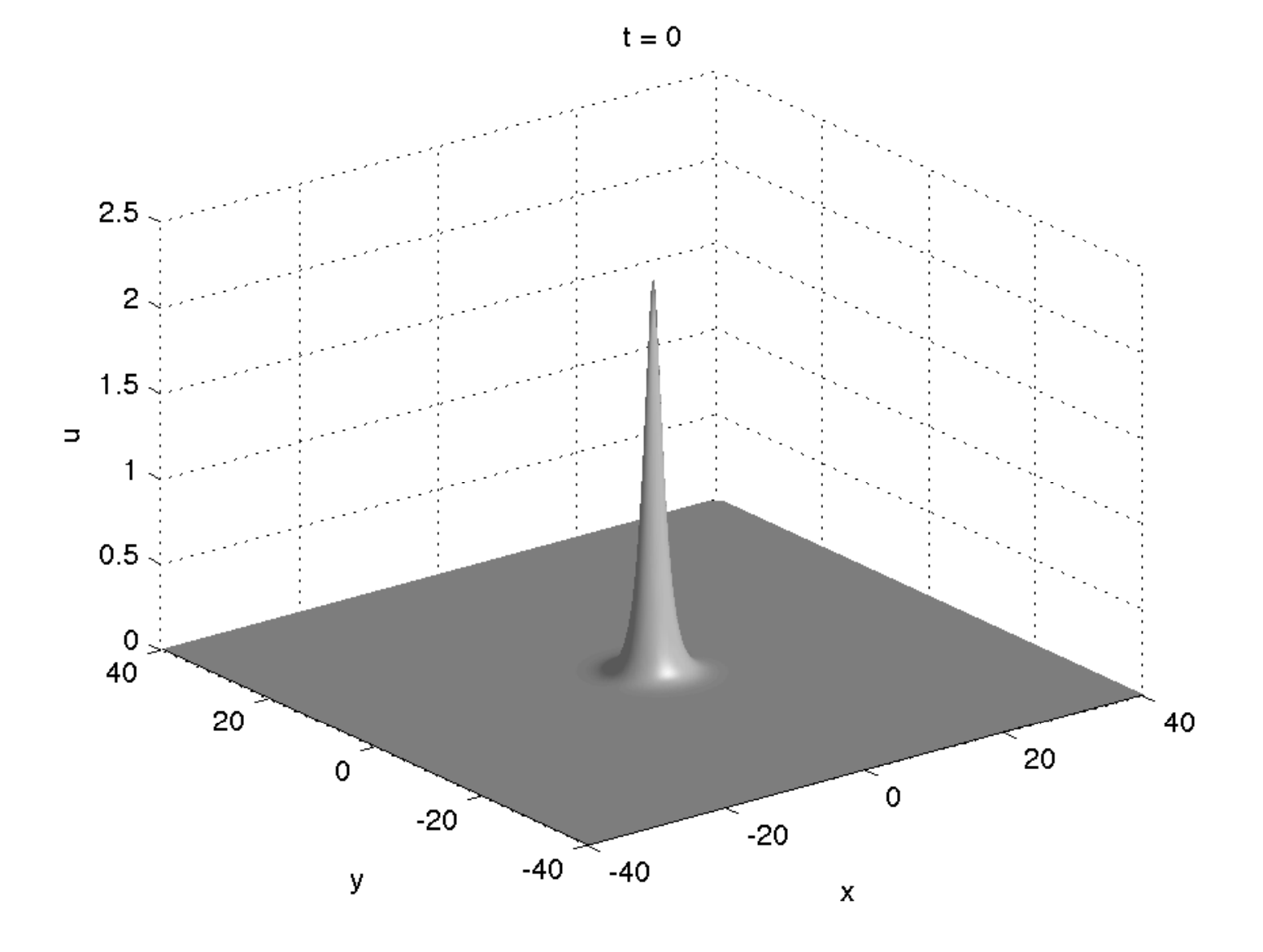}
(b)\includegraphics[width=2.75in]{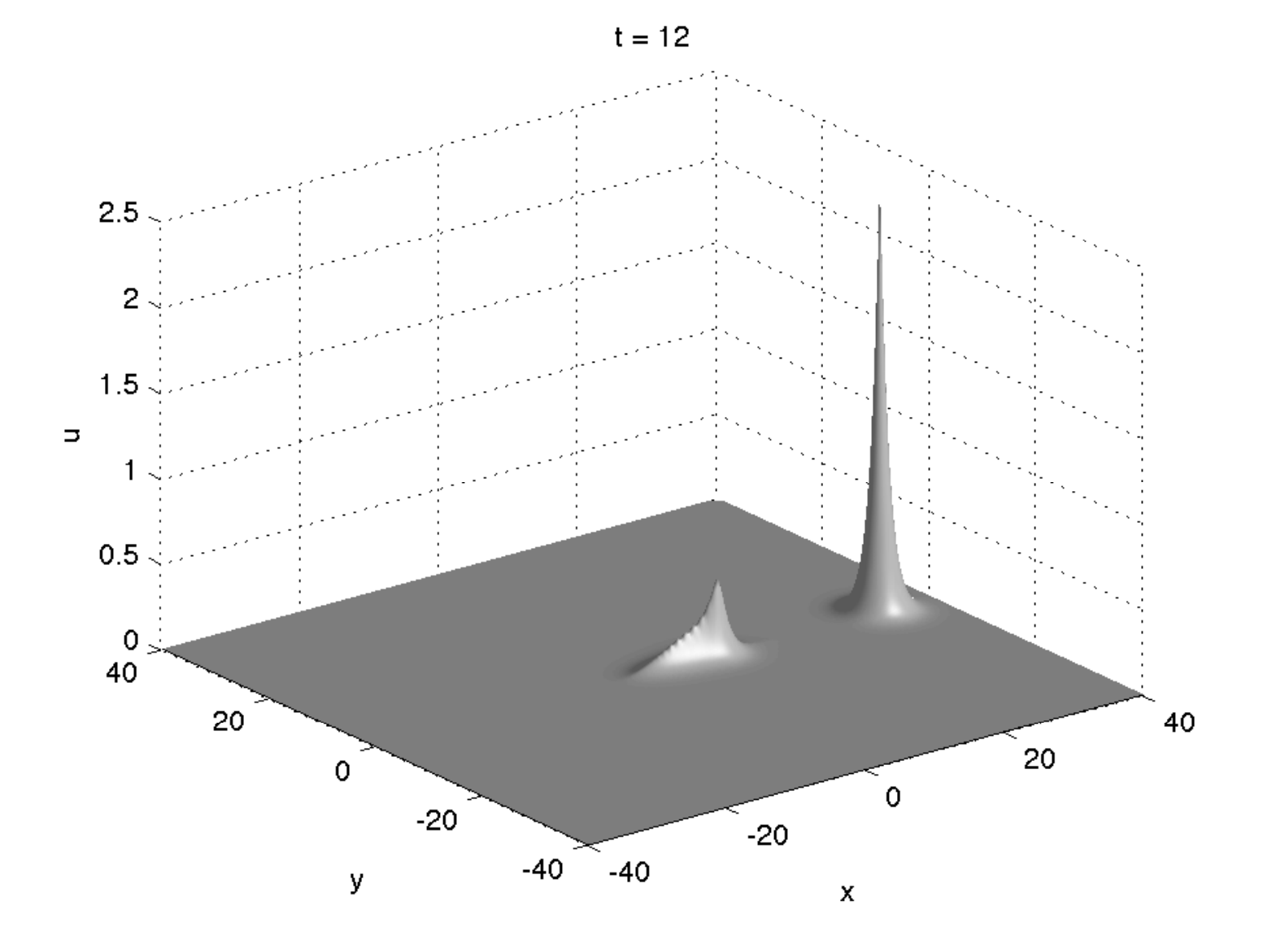}
%(a)\includegraphics[width=2.75in]{P_sech2_init}
%(b)\includegraphics[width=2.75in]{P_sech2_t_12}
\end{center}
\caption{ Simulations for the completely integrable system (\ref{eq:2D-qp}), $\bs{\kappa=0}$. Numerical parameters as in Figure~\ref{fig:overtake-32-direct}. (a) Initial first-component of $\bs{u}$. (b) First-component of $\bs{u}$ at $t=12$. The initial data have non-zero first-component momenta along the $x$-axis and zero second-component momenta everywhere. The non-zero momenta along $x$-axis are distributed by $\frac{1}{2}\sech^2(x)$.}
\label{fig:2D_qp_u}
\end{figure}

\begin{figure}[tbh]
\begin{center}
(a)\includegraphics[width=2.75in]{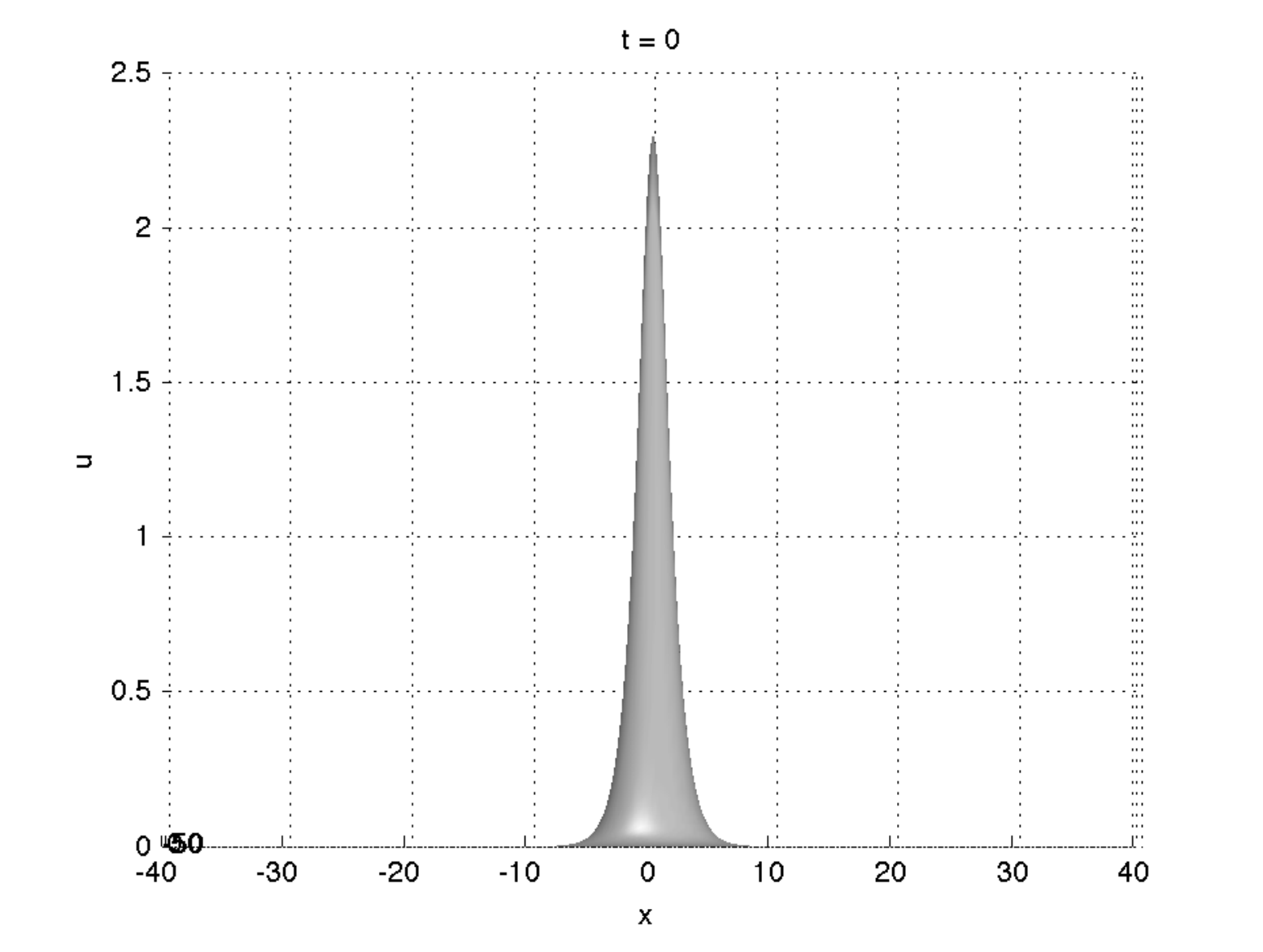}
(b)\includegraphics[width=2.75in]{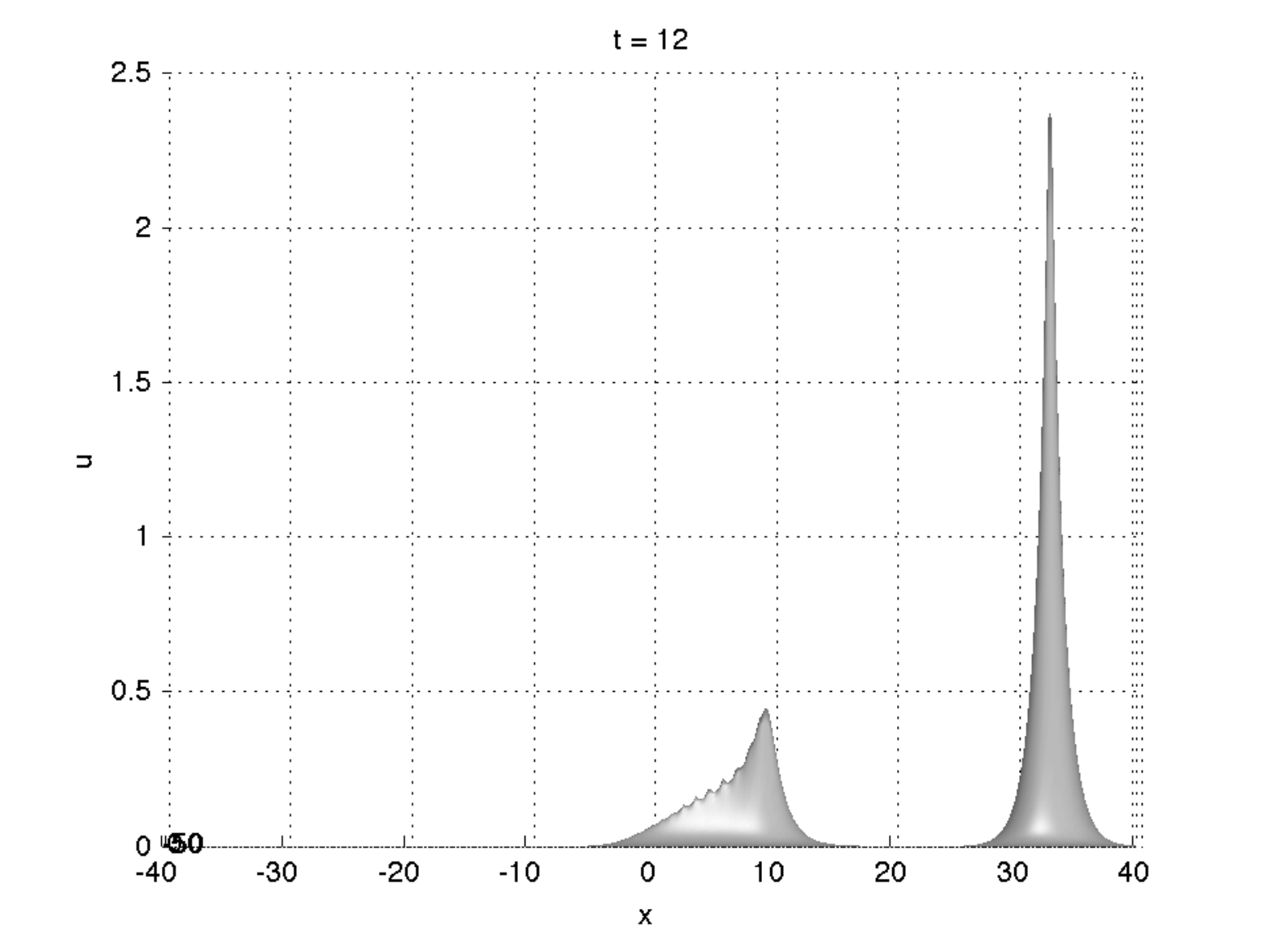}
%(a)\includegraphics[width=2.75in]{P_sech2_init_x_axis}
%(b)\includegraphics[width=2.75in]{P_sech2_t_12_x}
\end{center}
\caption{A frontal view of Figure \ref{fig:2D_qp_u}. The view direction is perpendicular to the $x$-axis. The initial smooth hump is sharpening as it moves to the right, followed by another hump emerging from the location of the initial hump. First-component of~$\bs{u}$ at $t=0$, (a), and at $t=12$, (b).}
\label{fig:2D_qp_u-x}
\end{figure}

\begin{figure}[tbh]
\begin{center}
(a)\includegraphics[width=2.75in]{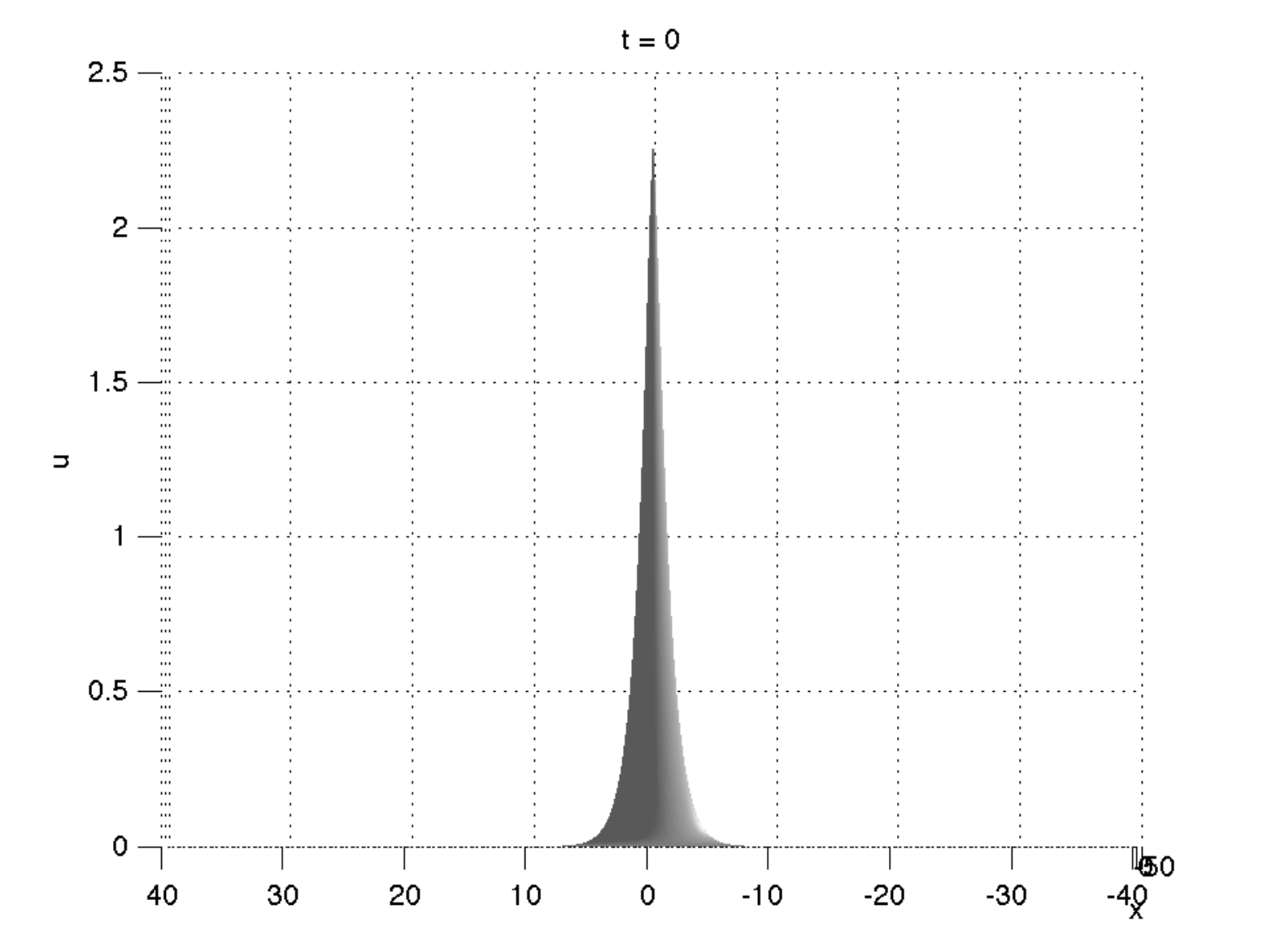}
(b)\includegraphics[width=2.75in]{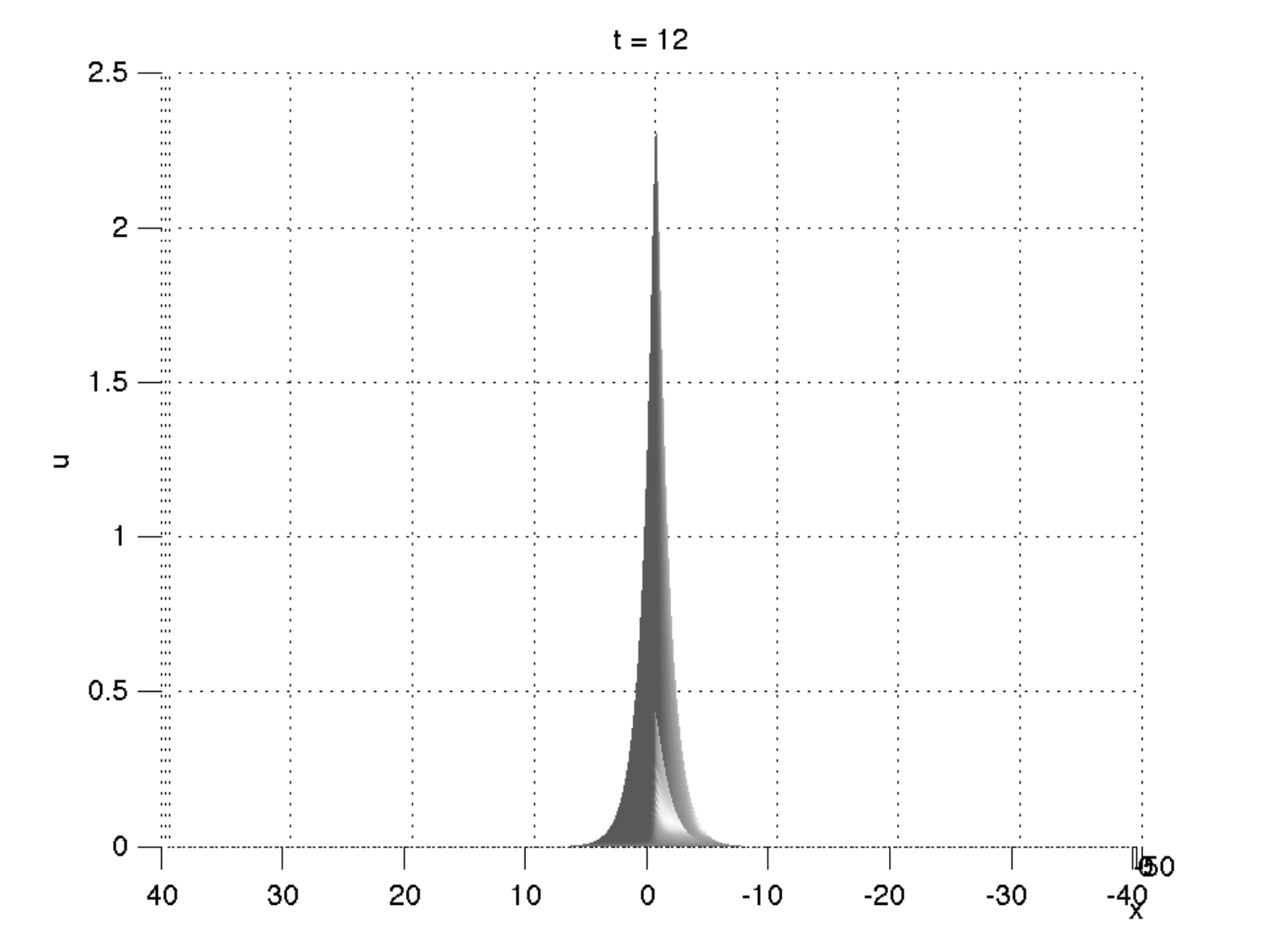}
%(a)\includegraphics[width=2.75in]{P_sech2_init_y_axis}
%(b)\includegraphics[width=2.75in]{P_sech2_t_12_y}
\end{center}
\caption{A frontal view of Figure \ref{fig:2D_qp_u}. The view direction is perpendicular to the $y$-axis. From this direction, the waves look the same as the conical solitary waves, for which the radial derivative has a finite jump. First-component of~$\bs{u}$ at $t=0$, (a), and at $t=12$, (b).}
\label{fig:2D_qp_u-y}
\end{figure}

Figure \ref{fig:2D_qp_u}(a) shows the initial first-component of $\bs{u}$, and (b) shows $\bs{u}$ at $t=12$.
Figure \ref{fig:2D_qp_u-x} is a frontal view of Figure \ref{fig:2D_qp_u}. The view direction is perpendicular to the $x$-axis. The figure shows that the slice along the $x$-axis is a smooth hump initially. Similar to the 1-D case, the initial hump sharpens as it moves to the right, followed by another hump which emerges from the location of the initial condition.
Figure \ref{fig:2D_qp_u-y} is another frontal view of Figure \ref{fig:2D_qp_u}. The view direction is perpendicular to the $y$-axis. From this view direction, the waves look the same as the conons, for which the radial derivative has a finite jump.   

We remark that simulations in this section use the full two-dimensional $N$-particle algorithm (\ref{eq:N-particle}) at the expense of computational cost. The number of particles in the calculation is $N=1280$ ($320$ particle on each slice in the $x$-direction and $4$ particles on the slice in the $y$-direction initially; compared with 1000 particles in the 1-D simulation in \cite{bib:chl06}). That is why the second emerging wave shown in Figures \ref{fig:2D_qp_u}  and \ref{fig:2D_qp_u-x} displays a saw-tooth-like roughness.  A much less expensive way to obtain a high-resolution result would be to compute the particle evolution on the $x$-axis only, i.e.,  evolve the $q$'s and $p$'s in equation (\ref{eq:2D-qp}), and then reconstruct the field $\bs{u}$ onto the whole two-dimensional plane by using equation (\ref{eq:u}) with the evolved $q$'s and $p$'s. Of course, in this manner we would not 
provide a full two-dimensional test of the numerics but rather use the analytical 
reduction to one dimensional settings. 

\subsection{Dispersive case $\bs{\kappa}\ne\bs{0}$}\label{sec:6.2}
In this section we demonstrate the effect of taking $\bs{\kappa}\ne 0$, i.e., of considering the dispersive deformation.

As noted in Section \ref{sec:dispersive}, the dispersion relation is non-trivial in the dispersive deformation. It depends on both amplitude and direction of the constant parameter $\bs{\kappa}$. Unlike its one-dimensional counterpart, this may have consequences in considering the limit $\bs{\kappa}\to 0$ as a possible dispersive regularization mechanism to handle non-smooth kernels. 

\begin{figure}[tbh]
\begin{center}
(a)\includegraphics[width=2.75in]{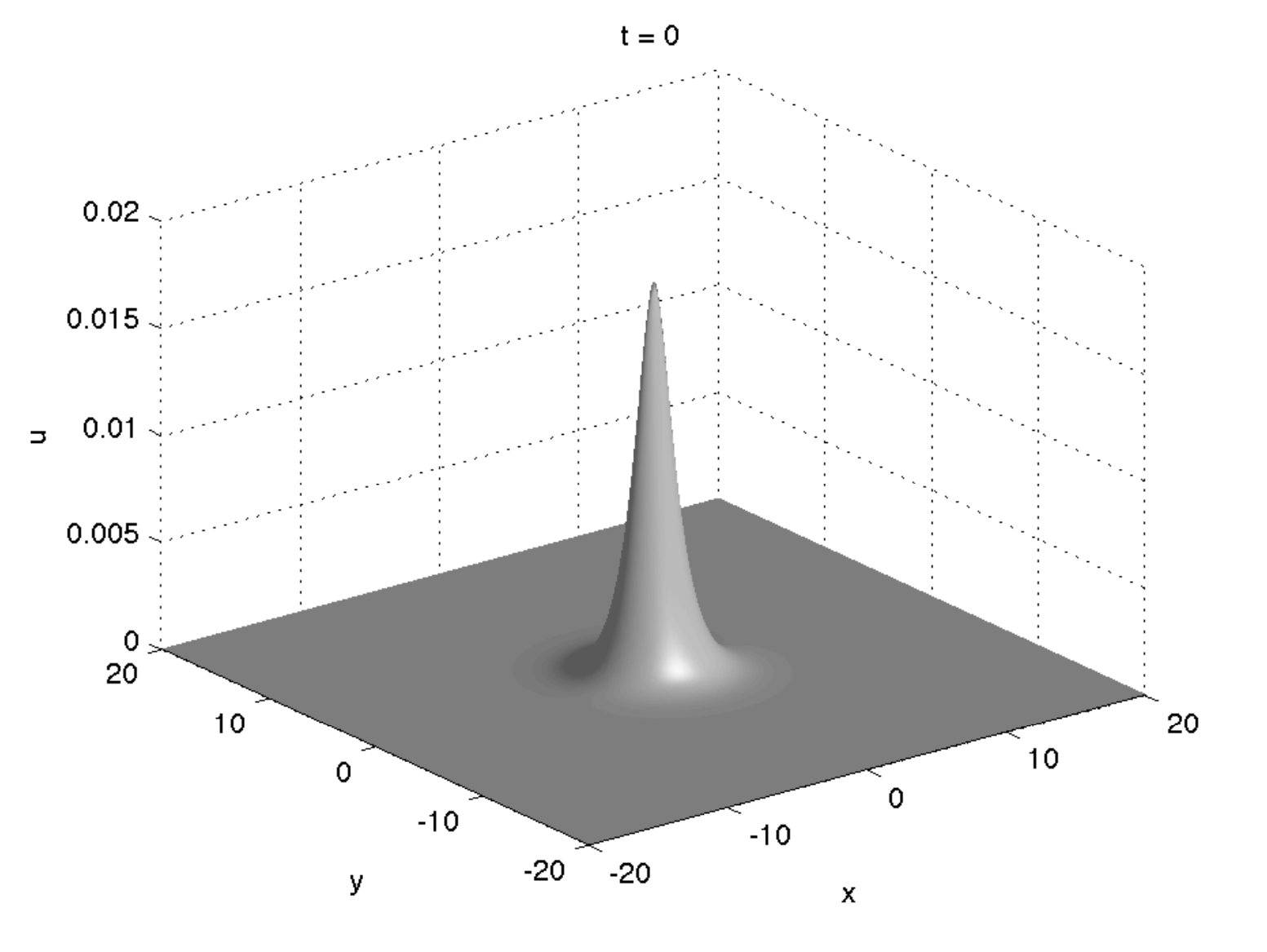}
(b)\includegraphics[width=2.75in]{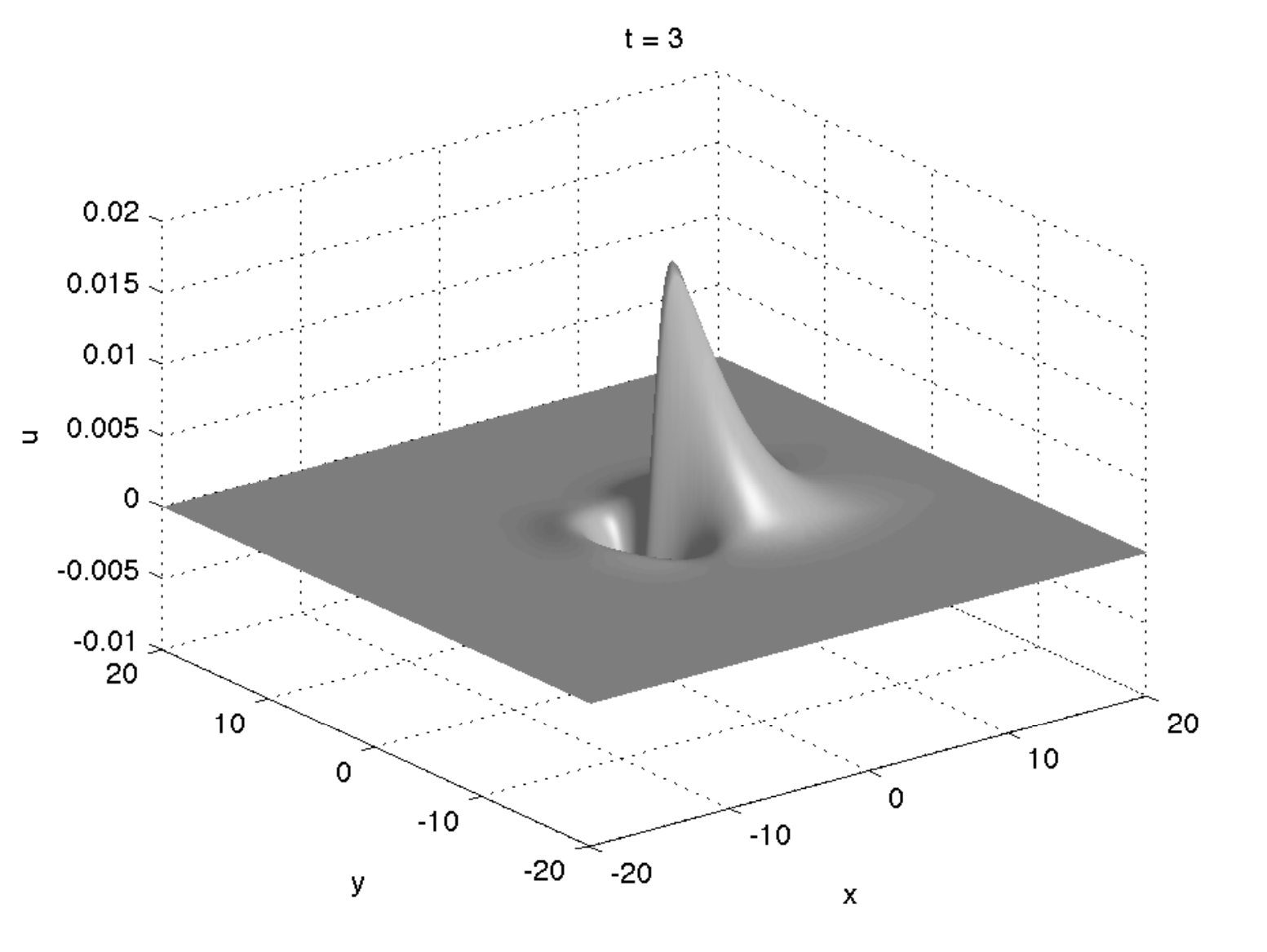}
%(a)\includegraphics[width=2.75in]{dispersion_initial}
%(b)\includegraphics[width=2.75in]{dispersion_kx1_ky0_t3}
\end{center}
\caption{ Same as Figure \ref{fig:2D_qp_u-x}, but $\bs{\kappa}\ne \bs{0}$. (a) Initial first-component of $\bs{u}$. (b) First component of $\bs{u}$ at $t=3$. The initial data have non-zero first-component momenta along the $x$-axis and zero second-component momenta everywhere. The non-zero momenta along $x$-axis are distributed by $\sech^2(x)$. The dispersive vector is ${\bs \kappa}= (1,\,0)$.}
\label{fig:2D_qp_u_k_ne_0}
\end{figure}

In our numerical experiments, the initial condition for the first component of $\bs{p}$ is specified by $p(x)=\sech^2(x)$ along the $x$-axis, and zero everywhere else, while the second component of ${\bs p}$ is zero everywhere, for the finite dimensional dynamical system (\ref{eq:Riemann-sum-pq}), corresponding to the dispersive PDE (\ref{eq:EPDIFFk}). However, unlike the previous non-dispersive  example where ${\bs\kappa} ={\bs 0}$, we carry out our numerical simulations by using equation (\ref{eq:Riemann-sum-pq}) directly, without absorbing $dx\,dy$ into the ${\bs p}$ variable. Furthermore, we use the original Green functions without normalization. Similar to the example for ${\bs\kappa}={\bs 0}$, the special case $ b=3/2$ is considered. Figure \ref{fig:2D_qp_u_k_ne_0}(a) shows the initial data. 
Let the dispersive (constant) vector be denoted by ${\bs \kappa}= (\kappa_1,\,\kappa_2)$.  Figure \ref{fig:2D_qp_u_k_ne_0}(b) shows the first-component of $\bs{u}$, evolving from the initial data in Figure \ref{fig:2D_qp_u_k_ne_0}(a) to the final time $t=3$ with the dispersive vector ${\bs \kappa}= (1,\,0)$.
The computational domain is $[-20, 20]\times[-5, 5]$ for ${\bs p}$ and ${\bs q}$. The mesh size is $dx=dy=0.2$. The temporal step size is $\Delta t=0.1$. The field $u$ is reconstructed on the domain $[-20, 20]\times[-20, 20]$. 
\begin{figure}[tbh]
\begin{center}
(a)\includegraphics[width=2.75in]{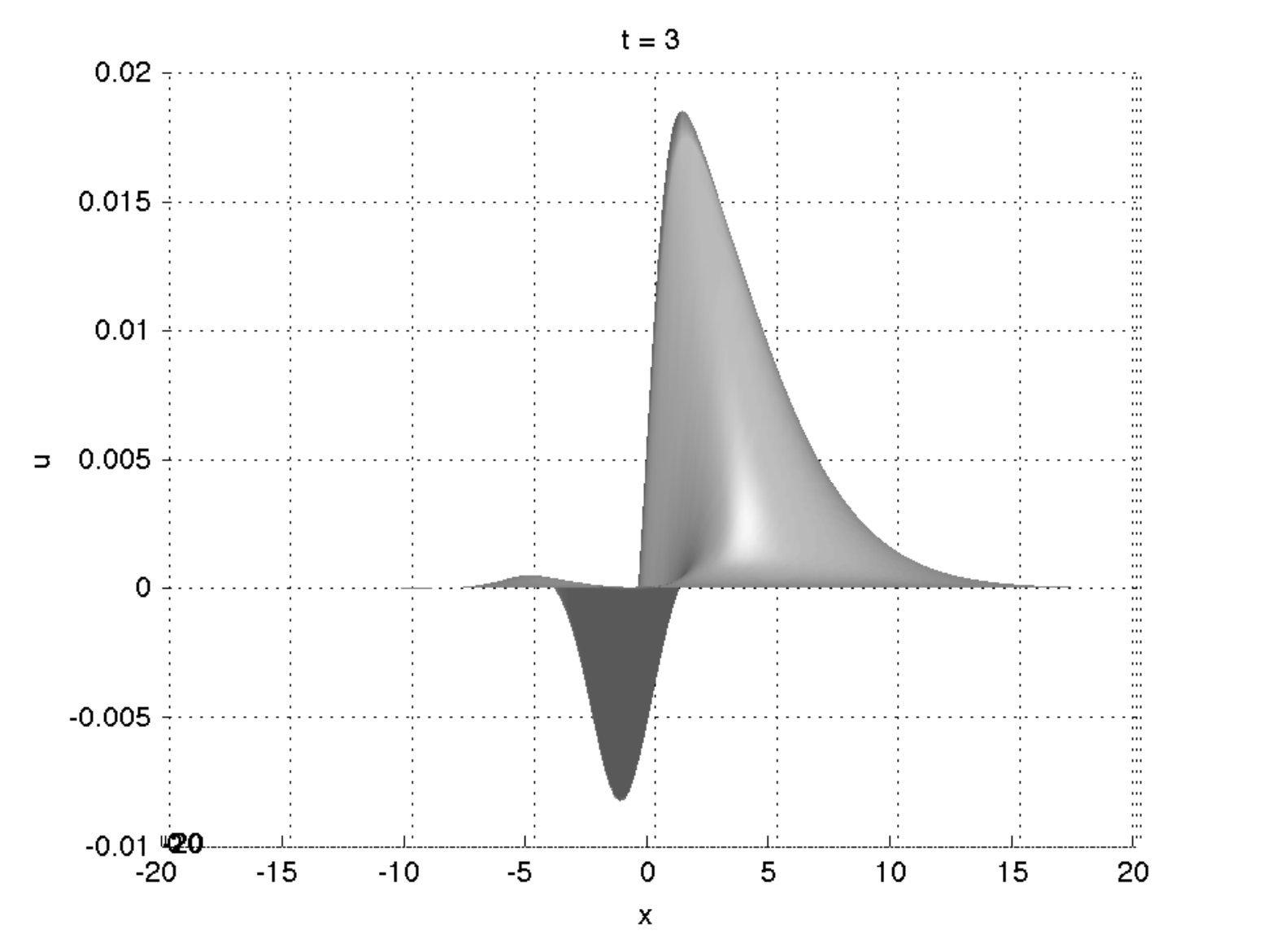}
(b)\includegraphics[width=2.75in]{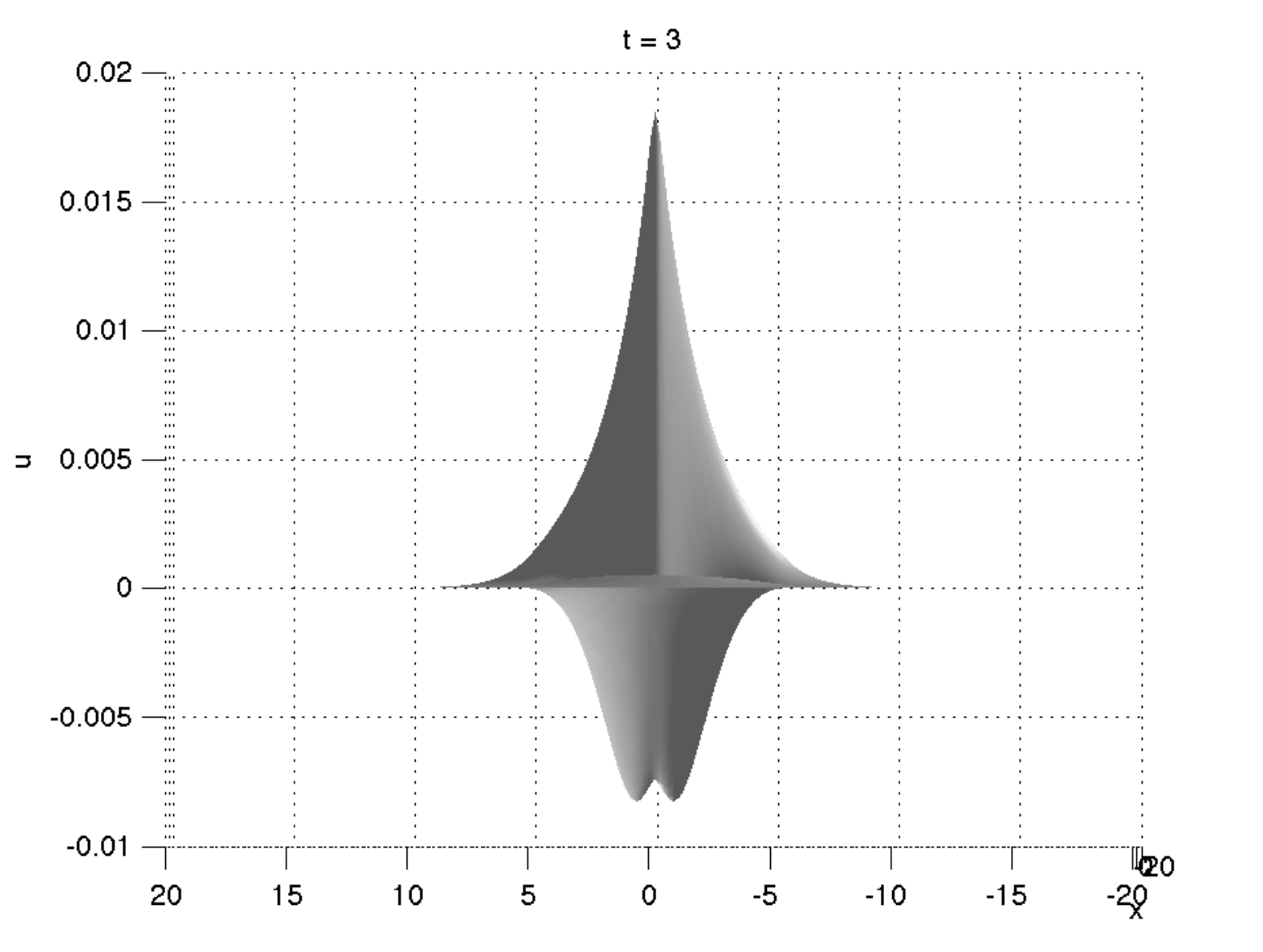}
%(a)\includegraphics[width=2.75in]{dispersion_kx1_ky0_t3_x_axis}
%(b)\includegraphics[width=2.75in]{dispersion_kx1_ky0_t3_y_axis}
\end{center}
\caption{Frontal views of Figure \ref{fig:2D_qp_u_k_ne_0}(b). (a) View direction  perpendicular to the $x$-axis. (b) View direction perpendicular to the $y$-axis.}
\label{fig:2D_qp_u_k_ne_0_frontal}
\end{figure} 
Figure\ref{fig:2D_qp_u_k_ne_0_frontal}(a) shows that along the $x$-axis the initial data evolve into a front advancing from left to right followed by an oscillatory wave train, similar to the example observed in \cite{bib:DCDS03} for the nonlinear SW equation. 
Further, in analogy with Figures \ref{fig:2D_qp_u_k_ne_0} and \ref{fig:2D_qp_u_k_ne_0_frontal}, Figure \ref{fig:2D_qp_u_k_ne_0_more} shows the numerical experiments for ${\bs \kappa}= (1,\,1)$ (left panel) and  ${\bs \kappa}= (0,\,1)$ (right panel), respectively. The development of oscillatory wave trains are observed in all directions, 
as expected when both $\kappa_1$ and $\kappa_2$ are nonzero. 

\begin{figure}[tbh]
\begin{center}
\includegraphics[width=2.75in]{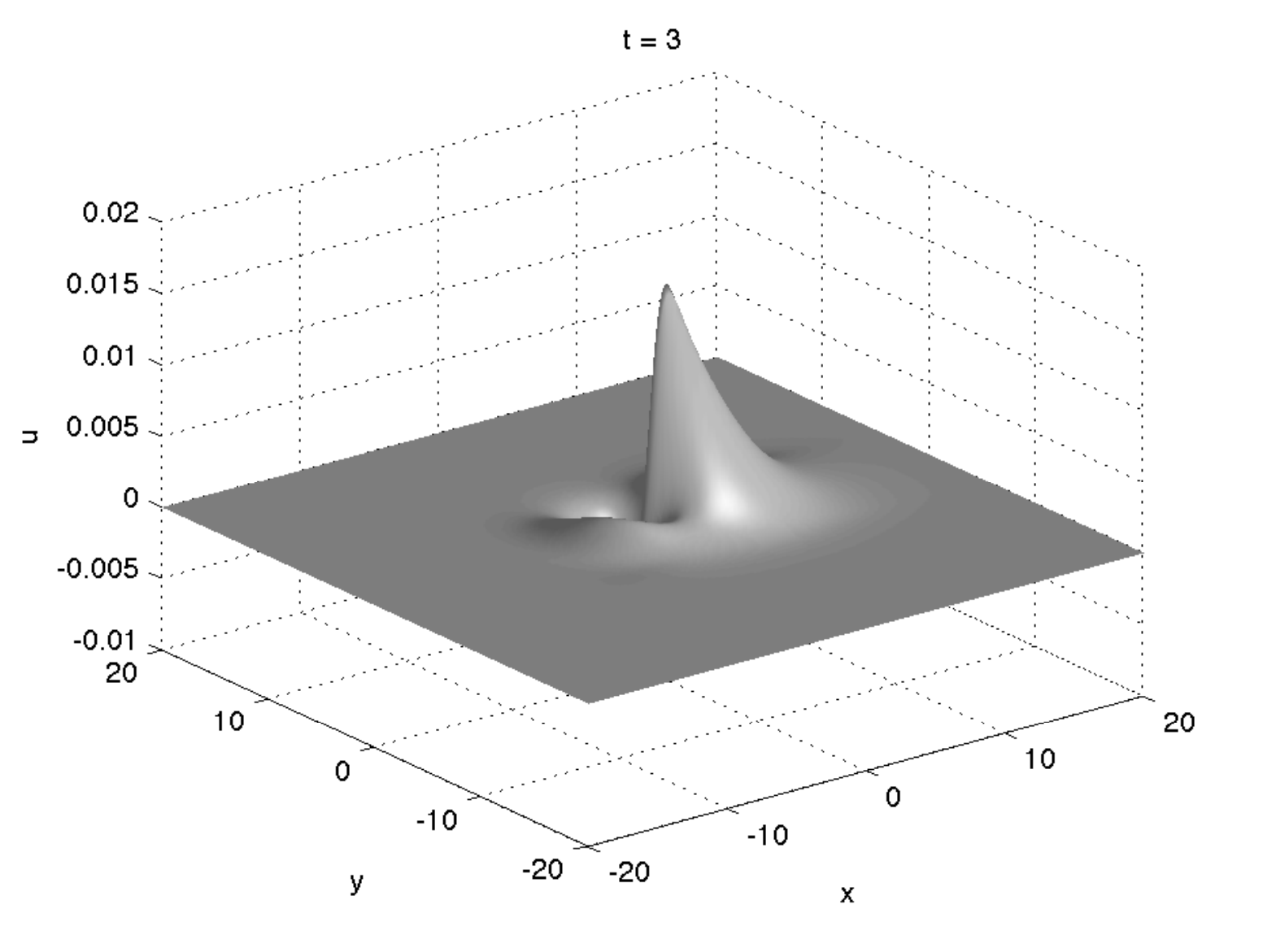}
\includegraphics[width=2.75in]{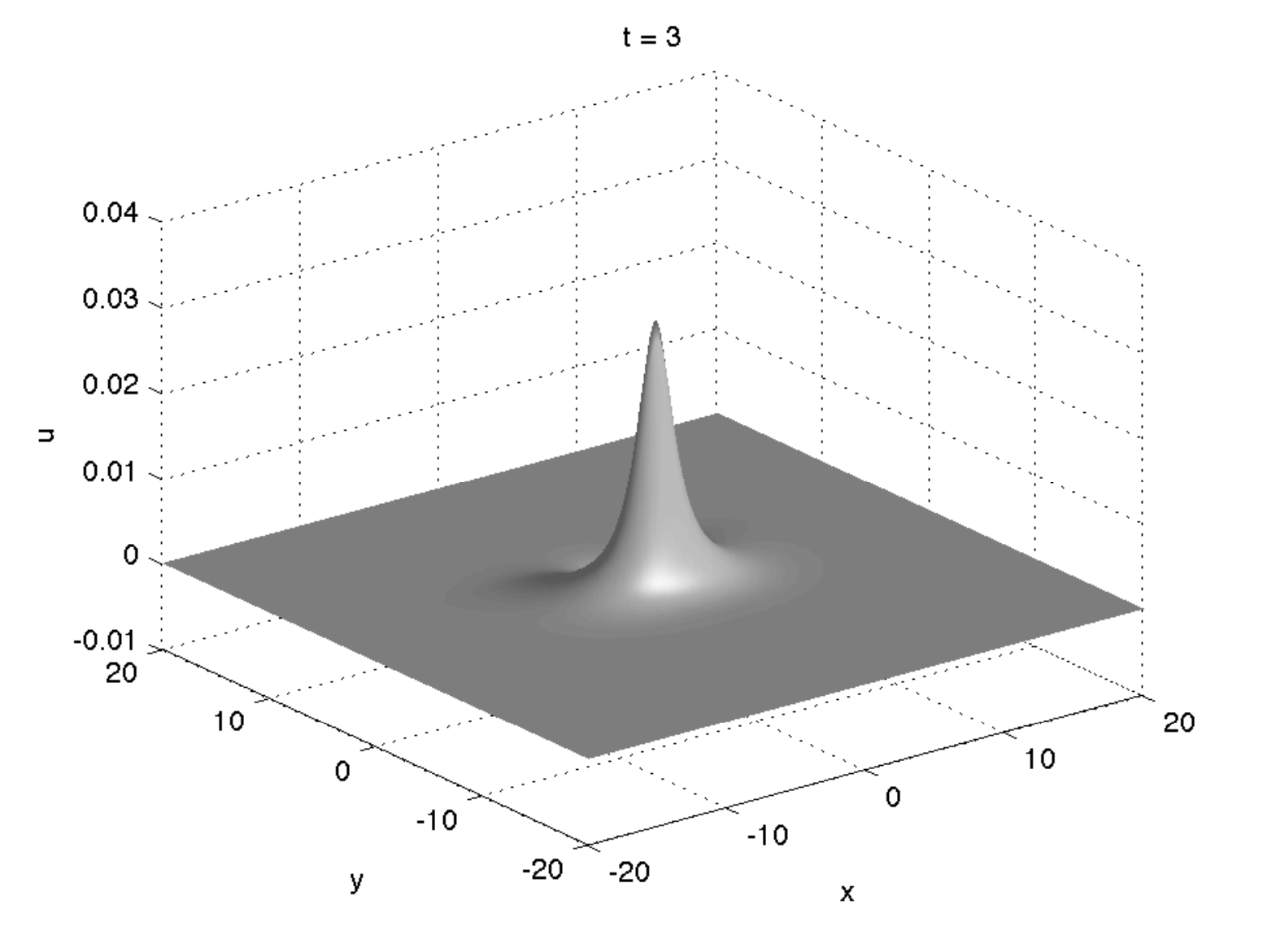}\\
\includegraphics[width=2.75in]{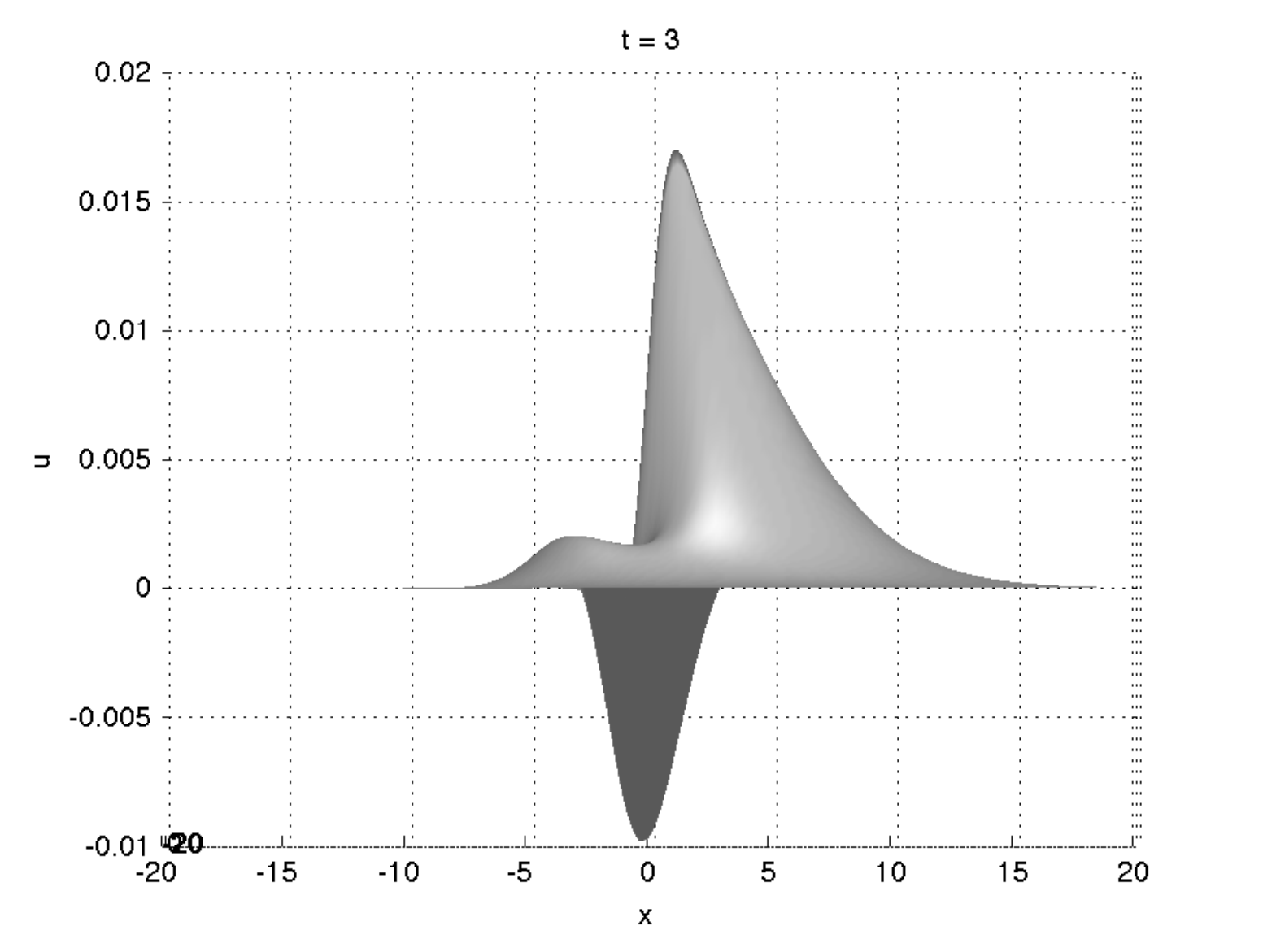}
\includegraphics[width=2.75in]{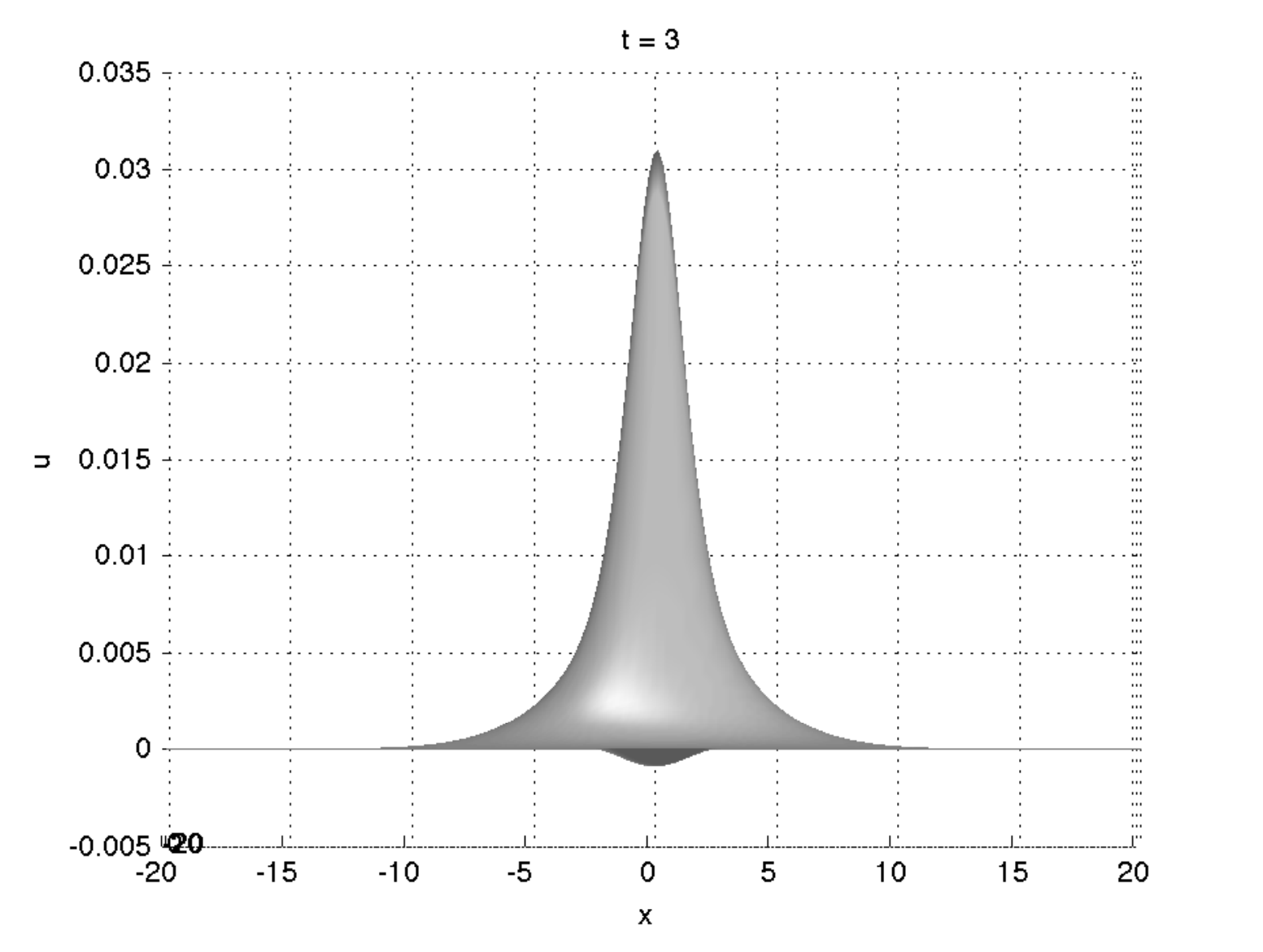}\\
(a)\includegraphics[width=2.75in]{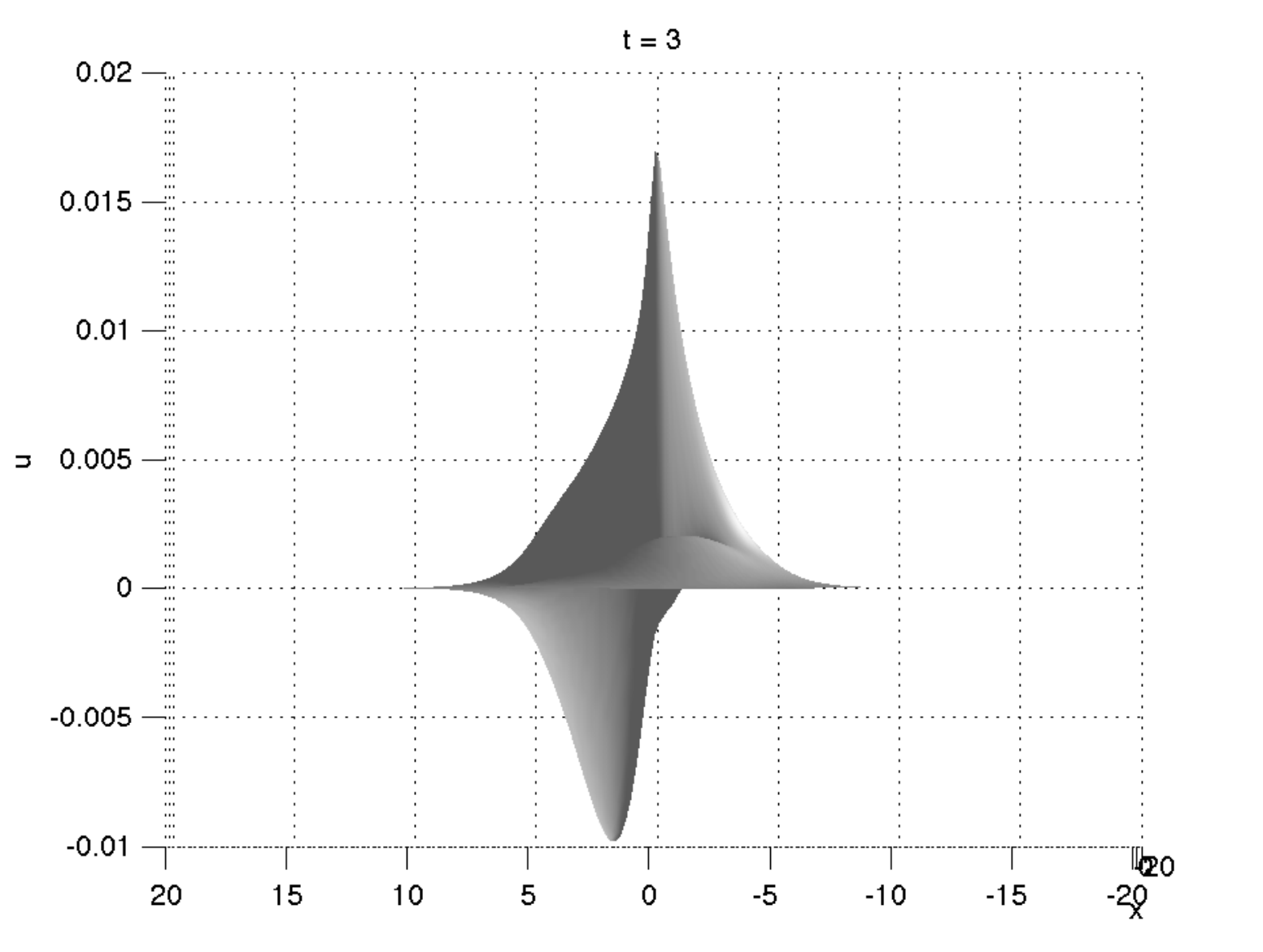}
(b)\includegraphics[width=2.75in]{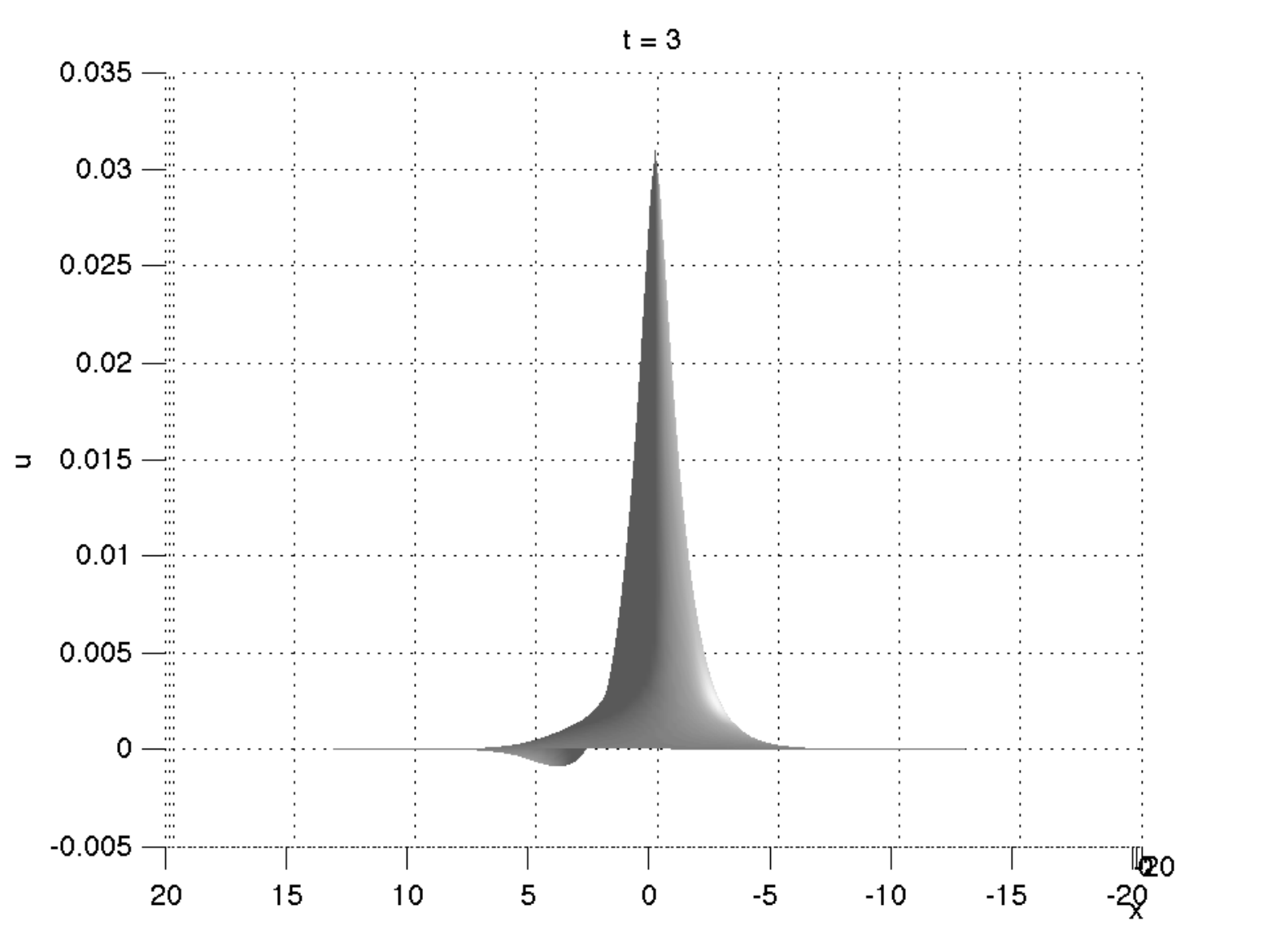}
%\includegraphics[width=2.75in]{dispersion_kx1_ky1_t3}
%\includegraphics[width=2.75in]{dispersion_kx0_ky1_t3}\\
%\includegraphics[width=2.75in]{dispersion_kx1_ky1_t3_x_axis}
%\includegraphics[width=2.75in]{dispersion_kx0_ky1_t3_x_axis}\\
%(a)\includegraphics[width=2.75in]{dispersion_kx1_ky1_t3_y_axis}
%(b)\includegraphics[width=2.75in]{dispersion_kx0_ky1_t3_y_axis}
\end{center}
\caption{Same as Figure \ref{fig:2D_qp_u_k_ne_0} and \ref{fig:2D_qp_u_k_ne_0_frontal} with different choices of $\bs{\kappa}$. Left panel, (a), ${\bs \kappa}= (1,\,1)$, and  right panel, (b), ${\bs \kappa}= (0,\,1)$. The final time is $t=3$ for both (a) and (b). Top to bottom are first-component of $\bs{u}$ in two dimensions, with frontal view from $x$-axis, and frontal view from $y$-axis.}
\label{fig:2D_qp_u_k_ne_0_more}
\end{figure} 

%\begin{figure}[tbh]
%\begin{center}
%(a)\includegraphics[width=2.75in]{dispersion_kx1_ky1_t3_x_axis-eps-converted-to.pdf}
%(b)\includegraphics[width=2.75in]{dispersion_kx0_ky1_t6_x_axis-eps-converted-to.pdf}
%\end{center}
%\caption{}
%\label{fig:2D_qp_u_k_ne_0_more_x}
%\end{figure} 
%
%
%\begin{figure}[tbh]
%\begin{center}
%(a)\includegraphics[width=2.75in]{dispersion_kx1_ky1_t3_y_axis-eps-converted-to.pdf}
%(b)\includegraphics[width=2.75in]{dispersion_kx0_ky1_t6_y_axis-eps-converted-to.pdf}
%\end{center}
%\caption{}
%\label{fig:2D_qp_u_k_ne_0_more_y}
%\end{figure} 

\section{Smooth initial data} 

Just as the particle algorithms developed for the SW equation \cite{bib:DCDS03, bib:Camassa05, bib:chl06, bib:Camassa07, bib:cl08}, the $N$-particle system in this paper can be seen as a Lagrangian numerical algorithm for solving the model PDEs (\ref{eq:EPDIFF}). The peakon of the particle method for the SW equation behaves like a member of a functional basis. This basis is advantageous not only for approximating rough initial data, but smooth data \cite{bib:chl06} can also be handled relatively well. This feature extends to the two-dimensional setting. 
%
%
%We have shown previously that similar to the peakon for the the SW equation, the Green functions for the $N$- particle algorithm are also basis functions that well represent rough data. 
In this section, we present an example with smooth initial data and follow numerically the ensuing solutions. We show that the $N$-particle algorithm can be used as a numerical method for solving the model PDEs (\ref{eq:EPDIFF}) in alternative to the traditional Eulerian methods for smooth solutions, if certain technical issues, such as the computational cost, can be overcome.  

In Section \ref{sec:7.1}, we introduce an operator-splitting pseudospectral algorithm. The method is called operator-splitting, because two sets of equations, one elliptic and one hyperbolic, are alternatively solved, other than solving a non-local integral-differential equation of $\bs{m}$. This operator-splitting method is introduced to assess the particle algorithm for handling smooth solution. In Section \ref{sec:6}, we have shown that the particle method is suitable for solutions with jump-derivatives at the peaks or with sharpening peaks. In Section \ref{sec:7.1},  we introduce the pseudospectral method to compare with the particle algorithm, in particular for problems with smooth initial data and smooth solutions at later times. We are interested in knowing how well the smooth solutions can be represented by a finite number of particles when particles cluster at some places and coarsen at the others, because in \cite{bib:Camassa05}, we showed that for the one-dimensional case, particle clustering might induce instability for the algorithm and cause blow-up, while particle coarsening would cause saw-tooth-like roughness for smooth solutions.

We remark that since the operator-splitting approach solves two sets of equations in alternating steps, the convergence property of the method to the true solution is a rather delicate problem, due to the splitting error. Even in the one-dimensional case \cite{bib:JCP_RB,bib:particle_RB,bib:DRP,bib:DRP2}, the method is not guaranteed to converge (although numerical convergences are established for both one and two dimensional algorithms). Nevertheless, the primary advantage for introducing the operator-splitting methods is to avoid solving a non-local integral-differential equation. The convergence and the error bound for the operator splitting are interesting open questions on their own right, even for the one-dimensional case, and thus belongs to a dedicated study and paper. 
\subsection{An operator-splitting pseudospectral method for the model system~(\ref{eq:EPDIFF})}\label{sec:7.1}
The two-dimensional version of equations~(\ref{eq:EPDIFF}), for which $\bs{x}=(x, y)^T$, $\bs{u}=(u_1, u_2)^T$, and $\bs{m}= (m_1, m_2)^T$, in component form is
\begin{equation}
\begin{split}
\frac{\partial m_1}{\partial t} + (u_1)_xm_1+(u_2)_xm_2+u_1(m_1)_x+u_2(m_1)_y+m_1((u_1)_x+(u_2)_y) &= 0,\\
\frac{\partial m_2}{\partial t} + (u_1)_ym_1+(u_2)_ym_2+u_1(m_2)_x+u_2(m_2)_y+m_2((u_1)_x+(u_2)_y) &= 0.
\end{split}
\end{equation}
After collecting terms of the above equations, together with equation (\ref{eq:elliptic}), we obtain
\begin{equation}\label{eq:2d_moment}
\begin{split}
\frac{\partial m_1}{\partial t} + (u_1m_1)_x+(u_2m_1)_y +m_1(u_1)_x+m_2(u_2)_x &= 0,\\
\frac{\partial m_2}{\partial t} + (u_1m_2)_x+(u_2m_2)_y +m_1(u_1)_y+m_2(u_2)_y &= 0,
\end{split}
\end{equation}
where
\begin{equation}\label{eq:2d_elliptic}
\begin{split}
m_1&=\left(1- a^2(\partial_{xx}+\partial_{yy})\right)^{ b}u_1,\\
m_2&=\left(1- a^2(\partial_{xx}+\partial_{yy})\right)^{ b}u_2.
\end{split}
\end{equation}

\begin{figure}[tbh]
\begin{center}
(a)\includegraphics[width=2.75in]{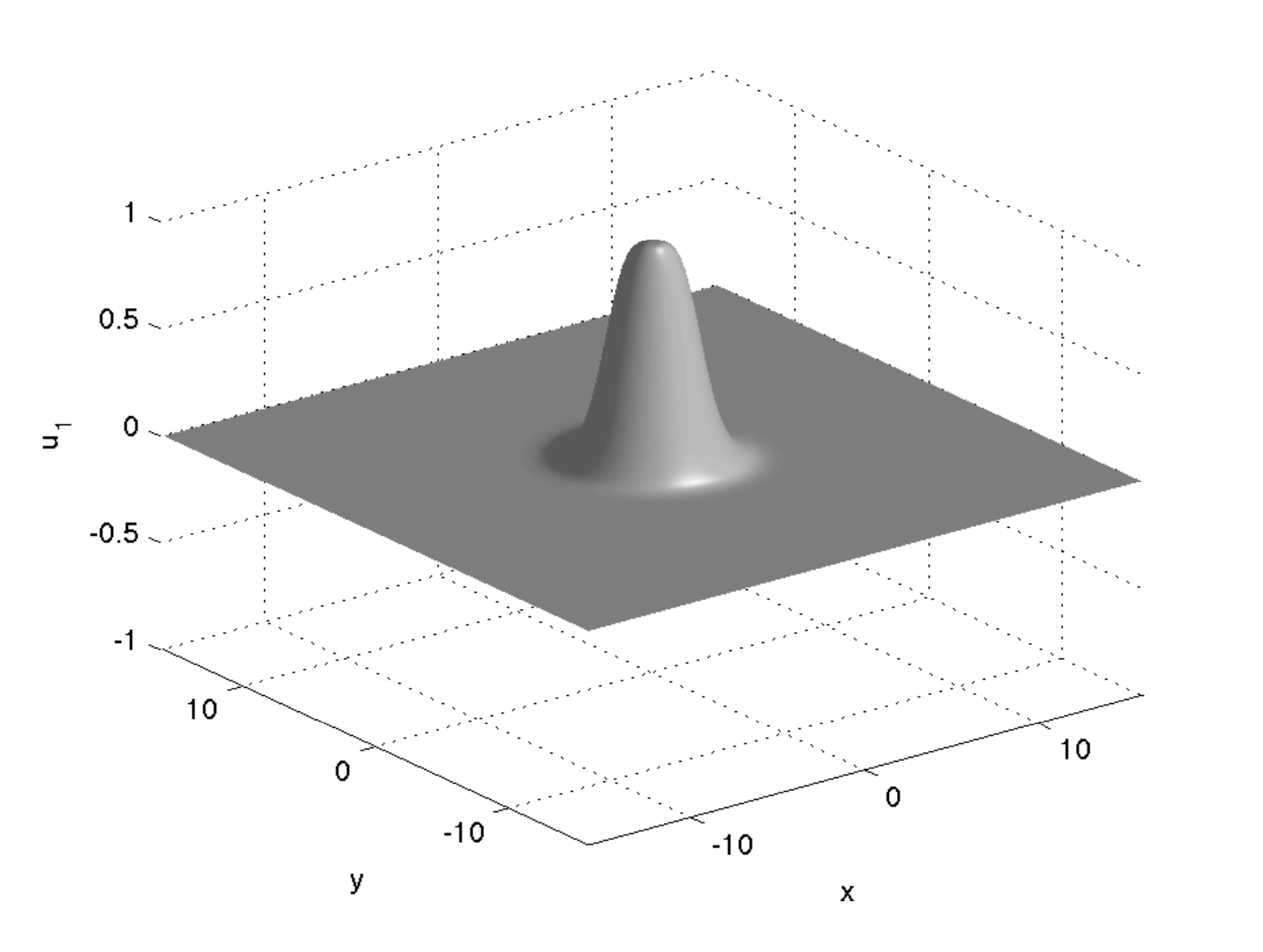}
%(b)\includegraphics[width=6in]{sech35-eps-converted-to.pdf}
(b)\includegraphics[width=2.75in]{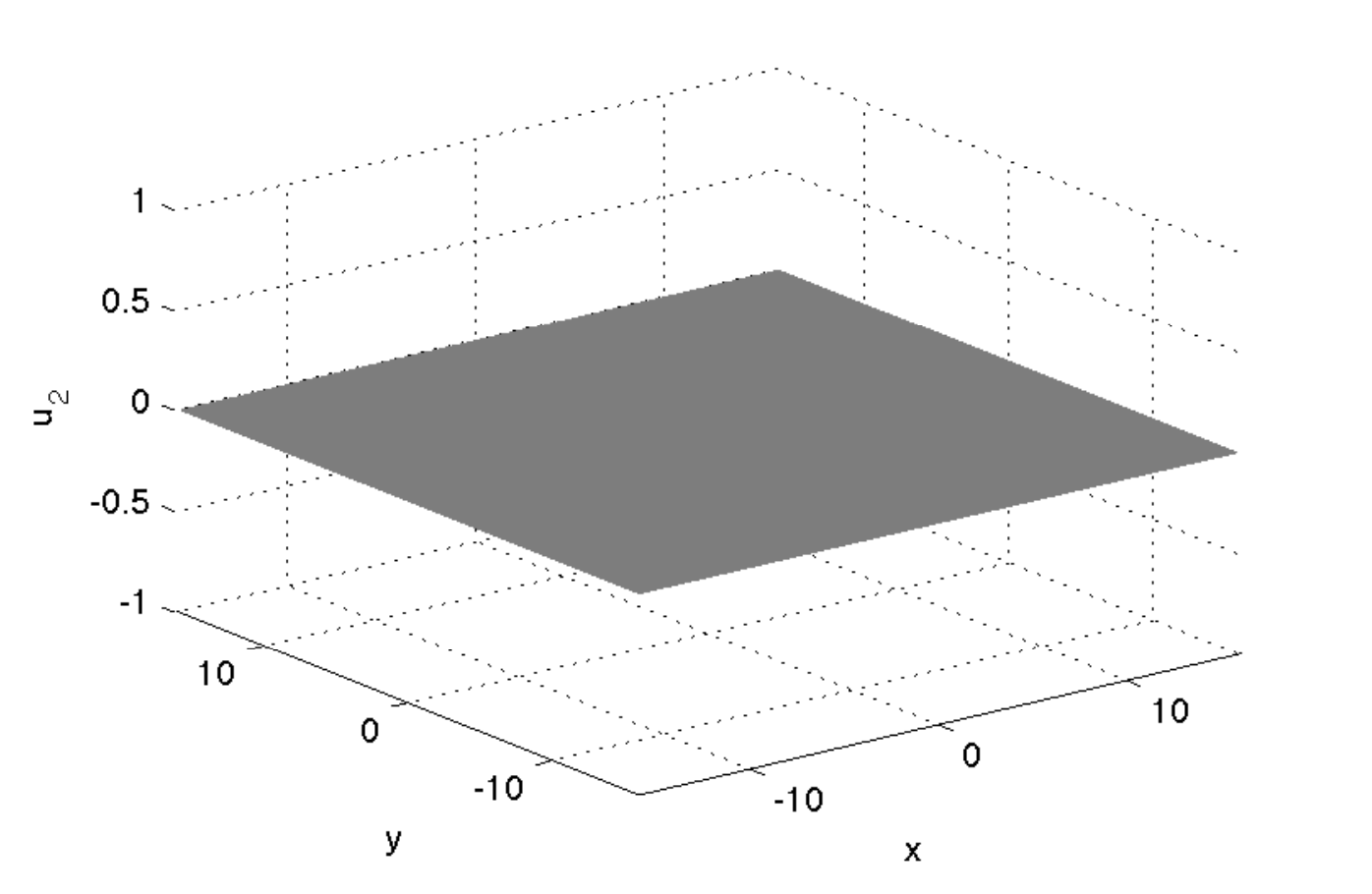}
\end{center}
\caption{The initial conditions for the comparison between particle and operator-splitting algorithms. 
%PDEs (\ref{eq:2d_moment}) and (\ref{eq:2d_elliptic}). 
(a) First component $u_1(x, y)=\sech\left(\ds{(x^2+y^2)}/{4}\right)$. (b) Second component $u_2(x, y)=0$. }
\label{fig:inital-data-smooth}
\end{figure}

\begin{figure}[tbh]
\begin{center}
(a)\includegraphics[width=2.75in]{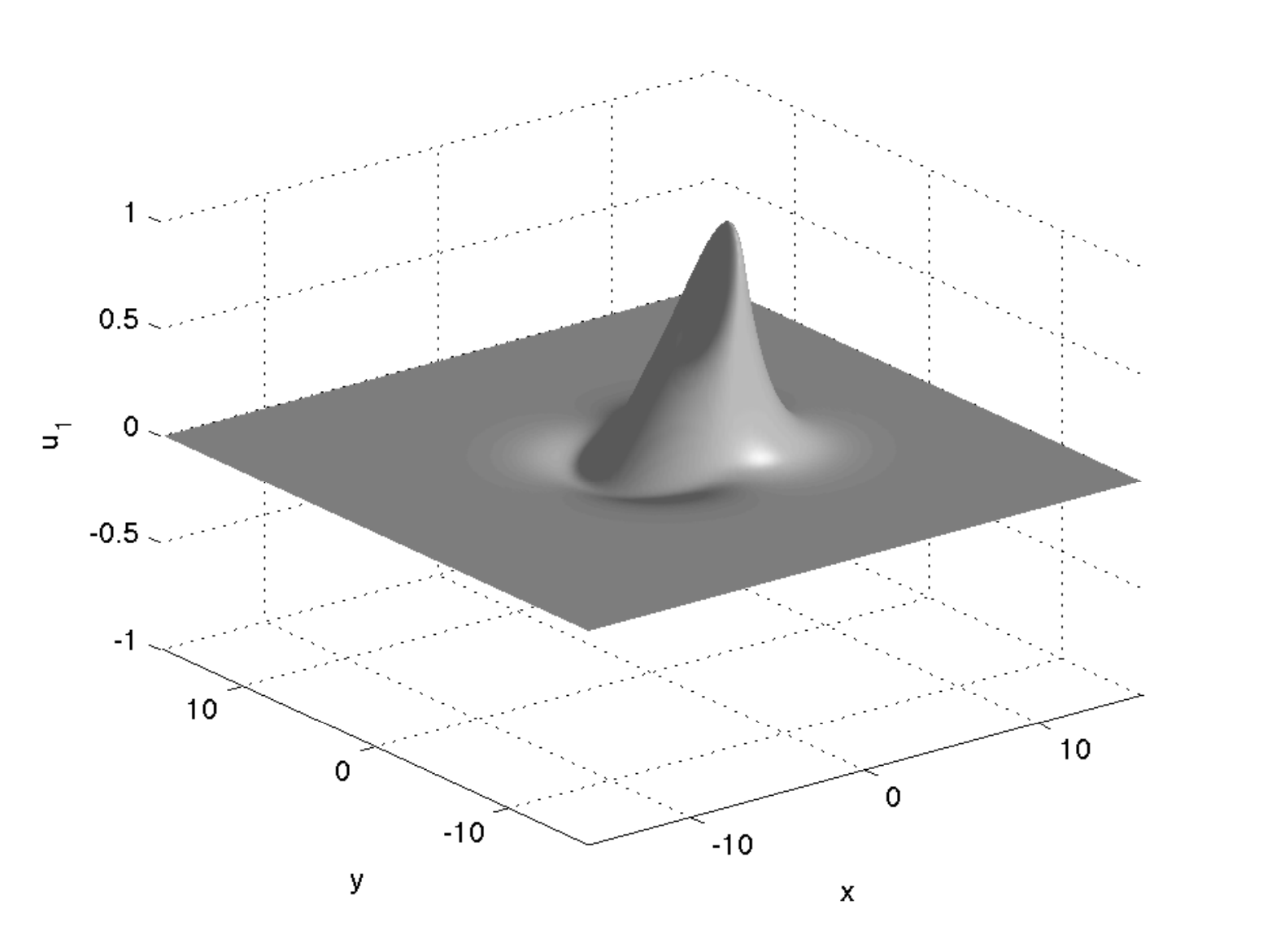}
%(b)\includegraphics[width=6in]{sech35-eps-converted-to.pdf}
(b)\includegraphics[width=2.75in]{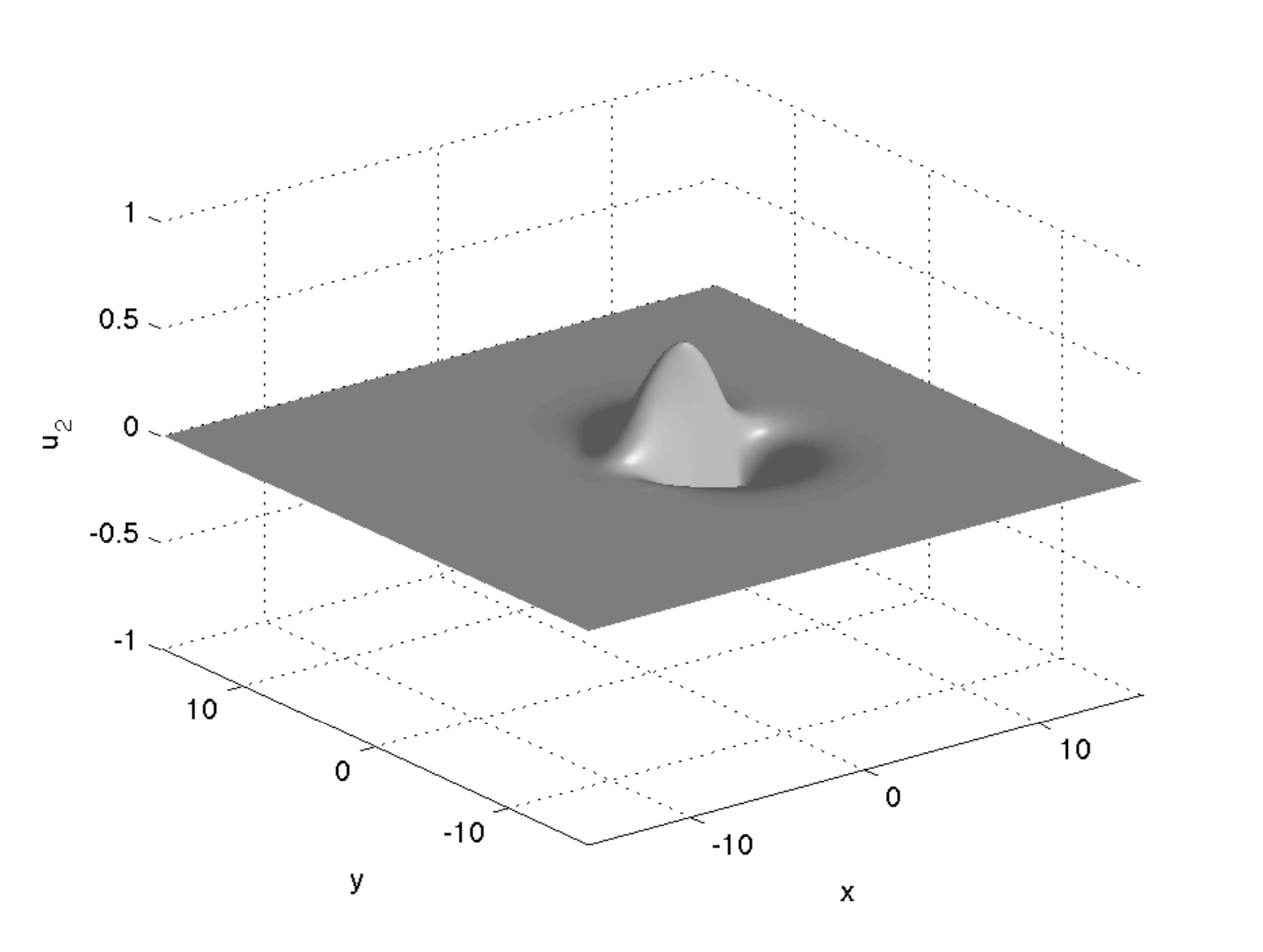}
\end{center}
\caption{Simulations for the initial data shown in Figure \ref{fig:inital-data-smooth} at $t=3$, by using the operator-splitting pseudospectral method, conical case $ b=3/2$, with $d x= d y=0.125$ and $\Delta t=0.0125$. The computational domain is $[-16, 16]\times [-16, 16]$. (a) The first component  $u_1$. (b) The second component  $u_2$.}
\label{fig:t3_pseudo}
\end{figure}

\begin{figure}[tbh]
\begin{center}
(a)\includegraphics[width=2.75in]{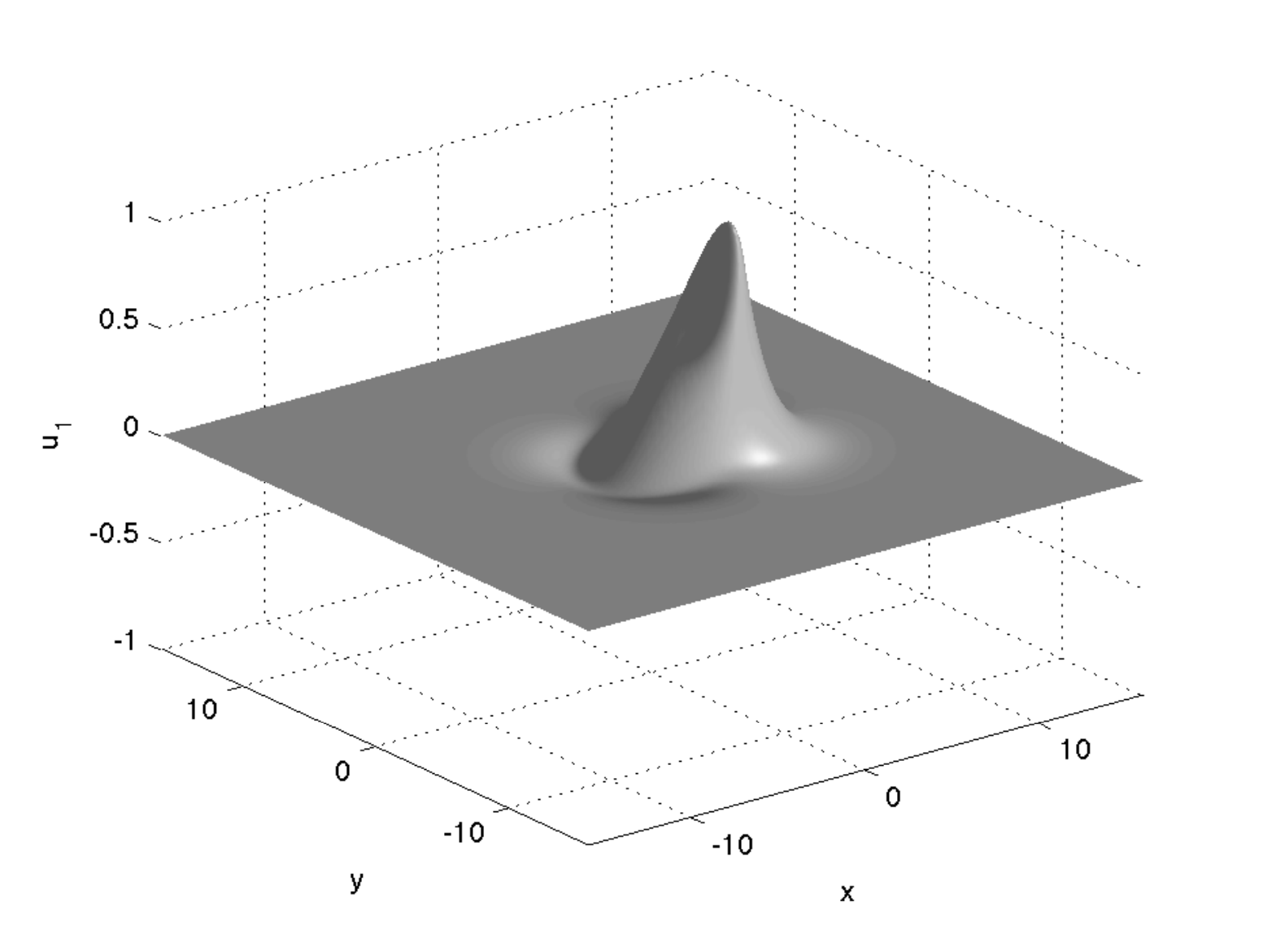}
(b)\includegraphics[width=2.75in]{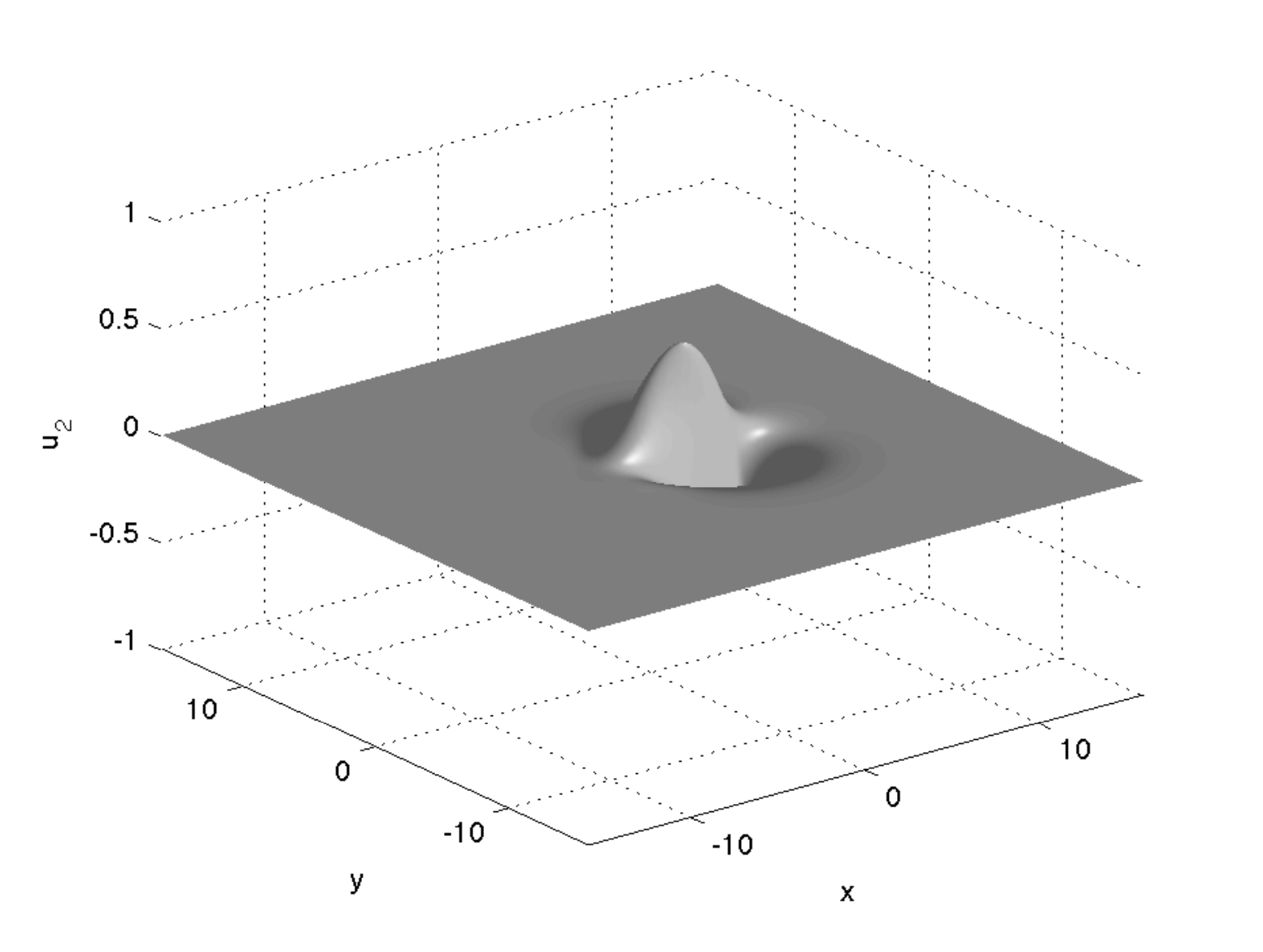}
\end{center}
\caption{Same as Figure \ref{fig:t3_pseudo}, but by using the $N$-particle algorithm with $81\times 81$ particles placed in a $[-8, 8]\times [-8, 8]$ domain initially. The time step is $\Delta t = 0.1$. The velocities are reconstructed on a $[-16, 16]\times [-16, 16]$ domain with $dx=dy=0.1$. (a) The first component $u_1$. (b) The second component $u_2$. The difference between Figure \ref{fig:t3_pseudo} and Figure \ref{fig:t3_particle}, for both components in the maximum norm, is at the order of $O(10^{-3})$. The two figures are virtually identical.}
\label{fig:t3_particle}
\end{figure}

\begin{figure}[tbh]
\begin{center}
(a)\includegraphics[width=2.75in]{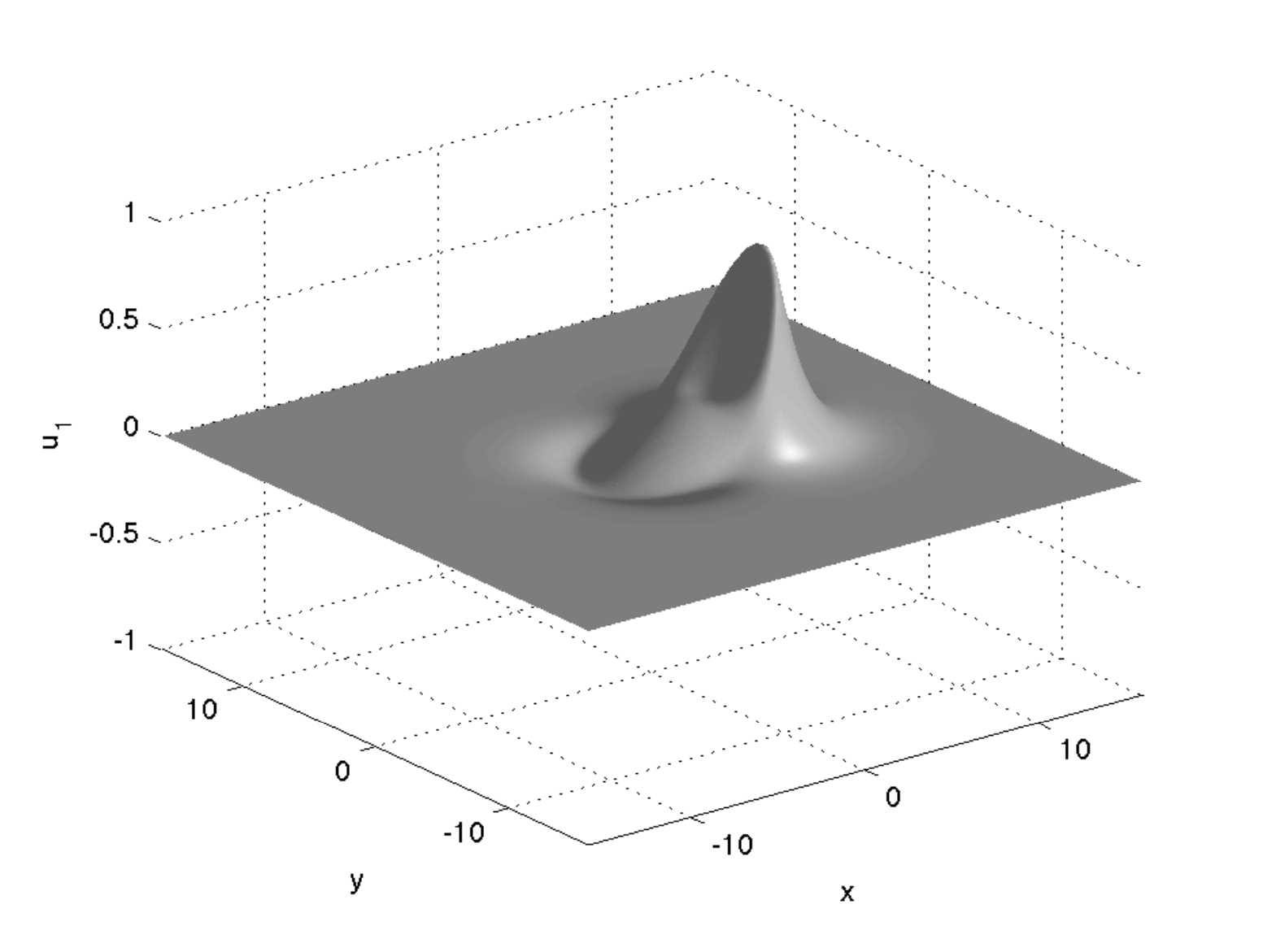}
%(b)\includegraphics[width=6in]{sech35-eps-converted-to.pdf}
(b)\includegraphics[width=2.75in]{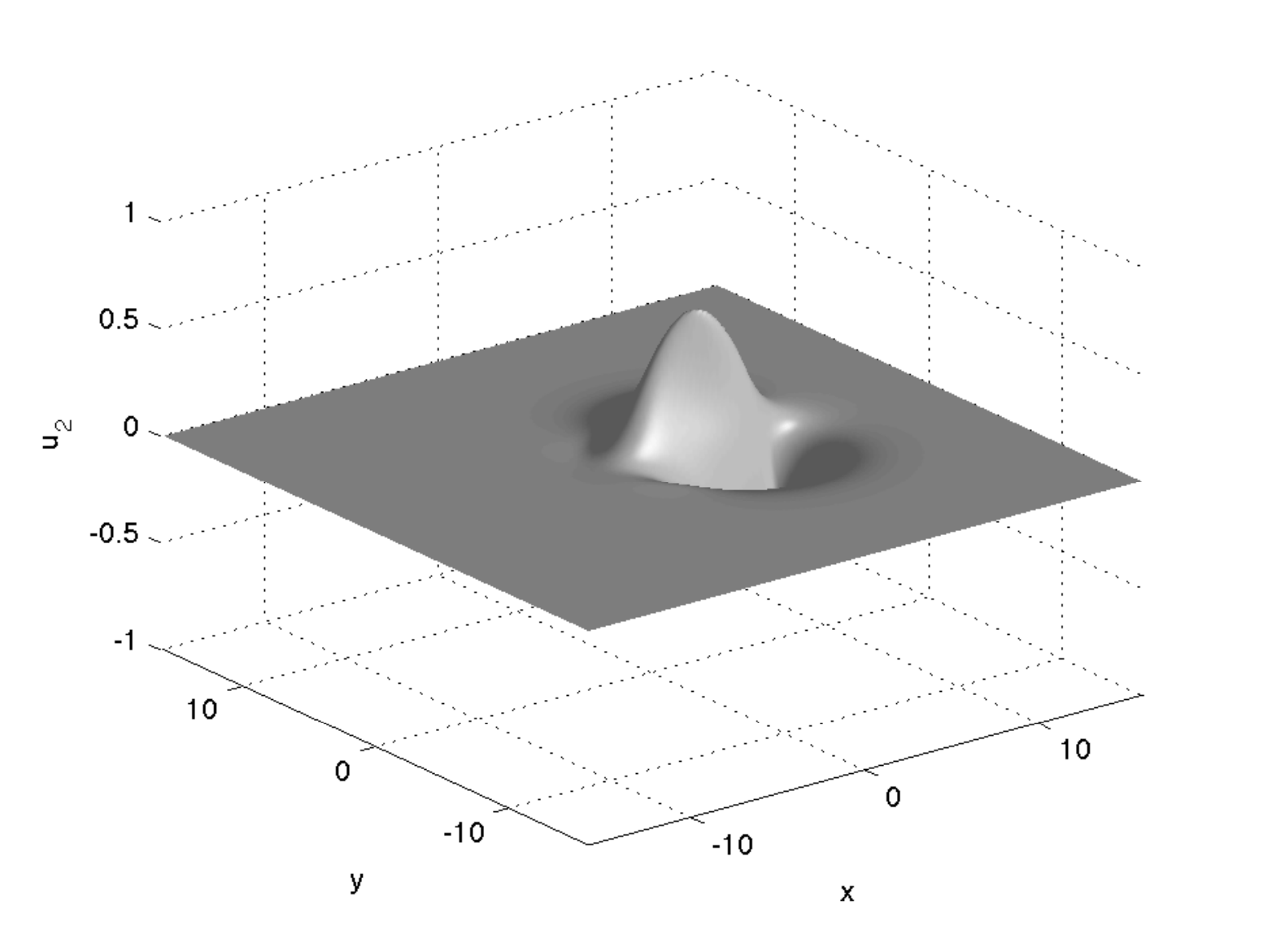}
\end{center}
\caption{Same as Figure \ref{fig:inital-data-smooth} but $t=5$, by using the operator-splitting pseudospectral method. $ b=3/2$ for the simulations. $d x= d y=0.25$ and $\Delta t=0.025$. The computational domain is $[-32, 32]\times [-32, 32]$. (a) The first component  $u_1$. (b) The second component $u_2$.}
\label{fig:t5_pseudo}
\end{figure}

\begin{figure}[tbh]
\begin{center}
(a)\includegraphics[width=2.75in]{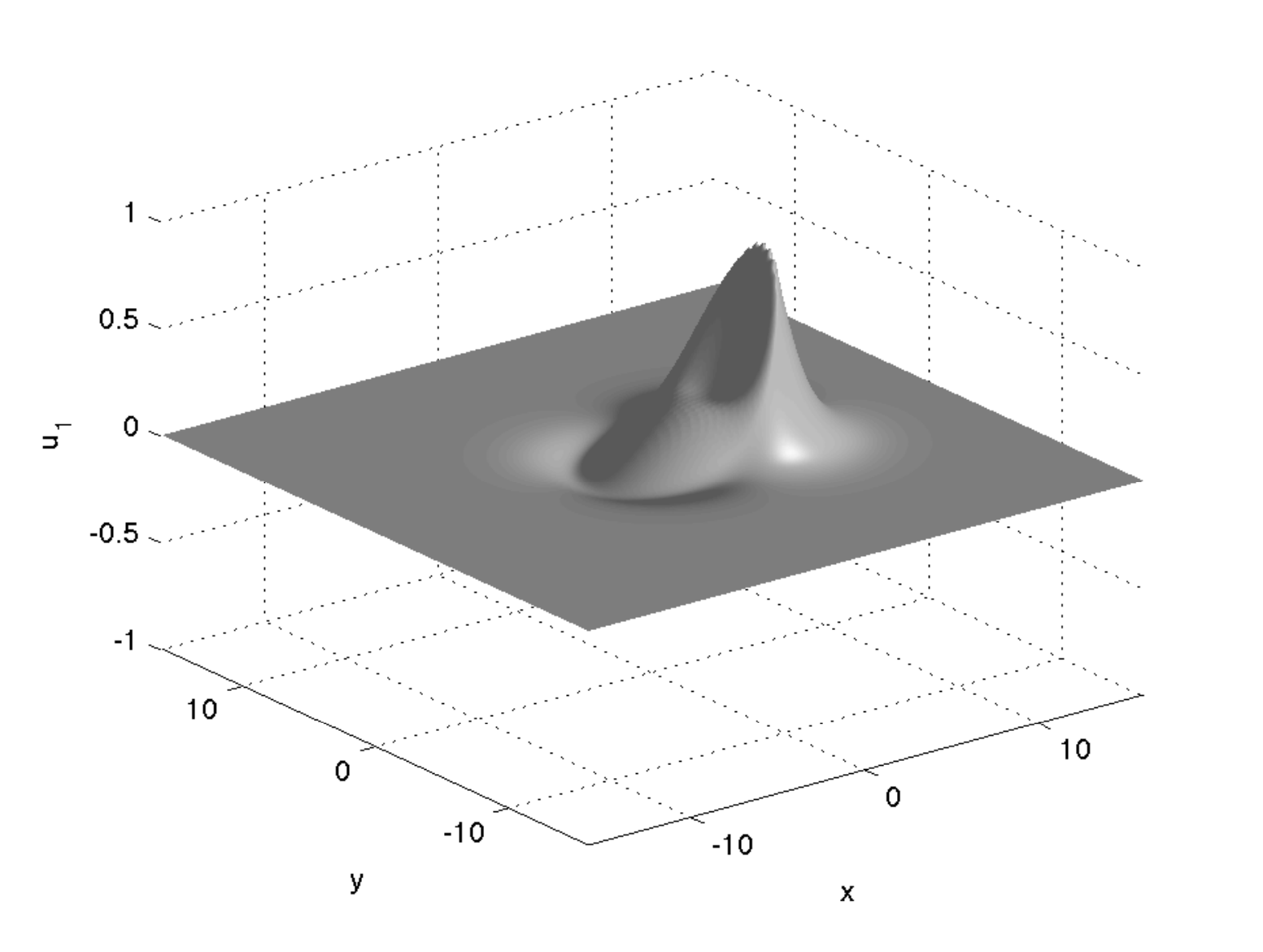}
%(b)\includegraphics[width=6in]{sech35-eps-converted-to.pdf}
(b)\includegraphics[width=2.75in]{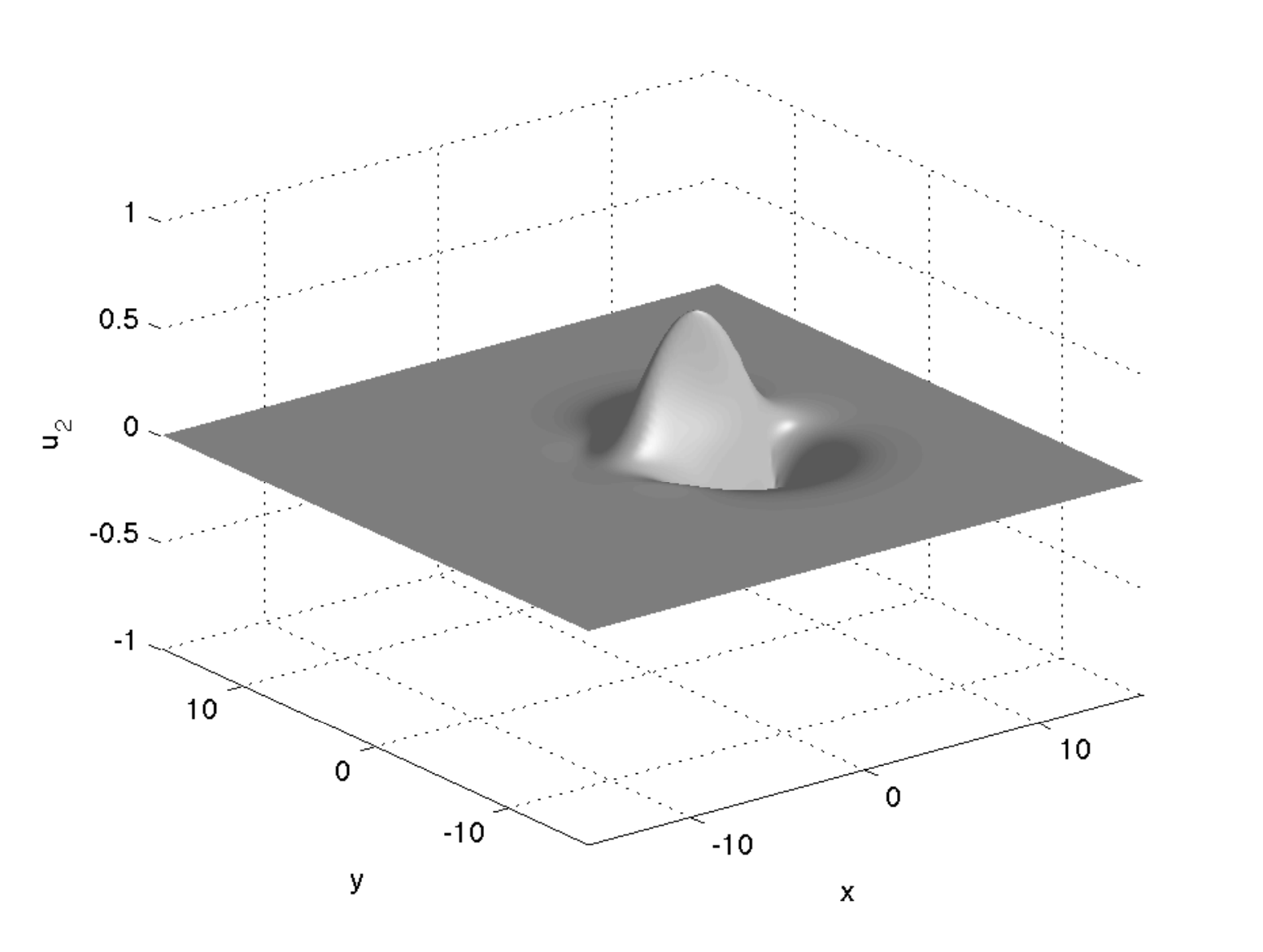}
\end{center}
\caption{Same  as Figure \ref{fig:t5_pseudo},  by the $N$-particle algorithm with $81\times 81$ particles placed in a $[-8, 8]\times [-8, 8]$ domain initially. The time step is $\Delta t = 0.1$. The velocities are reconstructed on a $[-16, 16]\times [-16, 16]$ domain with $dx=dy=0.1$. (a) First component  $u_1$. (b) Second component $u_2$. There is visible saw-tooth like roughness near the tip of the $u_1$ wave. }
\label{fig:t5_particle}
\end{figure}

\begin{figure}[tbh]
\begin{center}
(a)\includegraphics[width=2.75in]{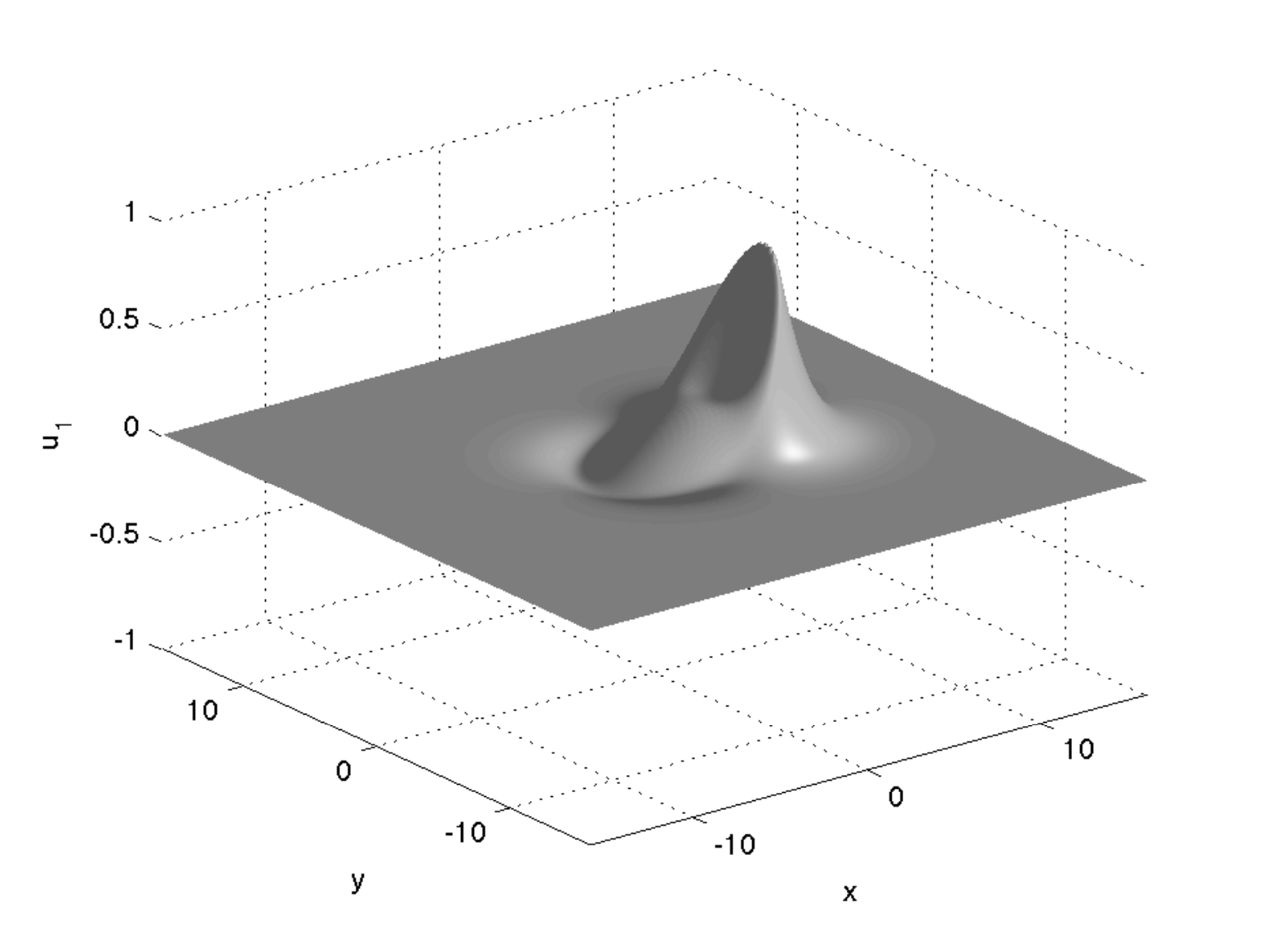}
%(b)\includegraphics[width=6in]{sech35-eps-converted-to.pdf}
(b)\includegraphics[width=2.75in]{t5_u2_particle_81x81-eps-converted-to.pdf}
\end{center}
\caption{Same simulations as Figure \ref{fig:t5_particle}, obtained by using the $N$-particle algorithm with $101\times 101$ particles placed in a $[-8, 8]\times [-8, 8]$ domain initially. The time step is $\Delta t = 0.1$. The velocities are reconstructed on a $[-16, 16]\times [-16, 16]$ domain with $dx=dy=0.08$. (a) The first component $u_1$. (b) The second component $u_2$. The tip of the $u_1$ wave is smoother than that in Figure \ref{fig:t5_particle}.}
\label{fig:t5_particle_finer}
\end{figure}

\begin{figure}[tbh]
\begin{center}
(a)\includegraphics[width=2.75in]{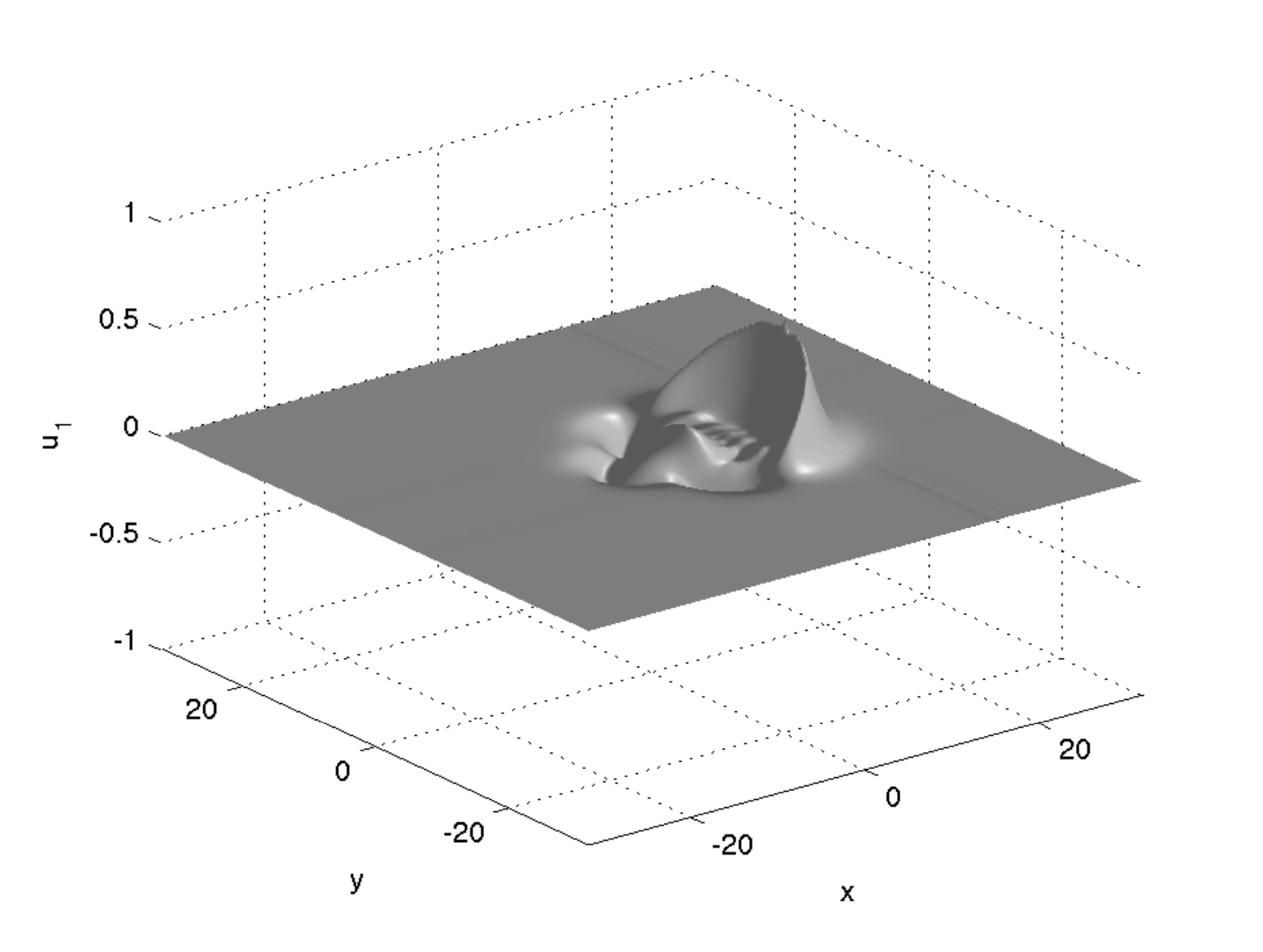}
%(b)\includegraphics[width=6in]{sech35-eps-converted-to.pdf}
(b)\includegraphics[width=2.75in]{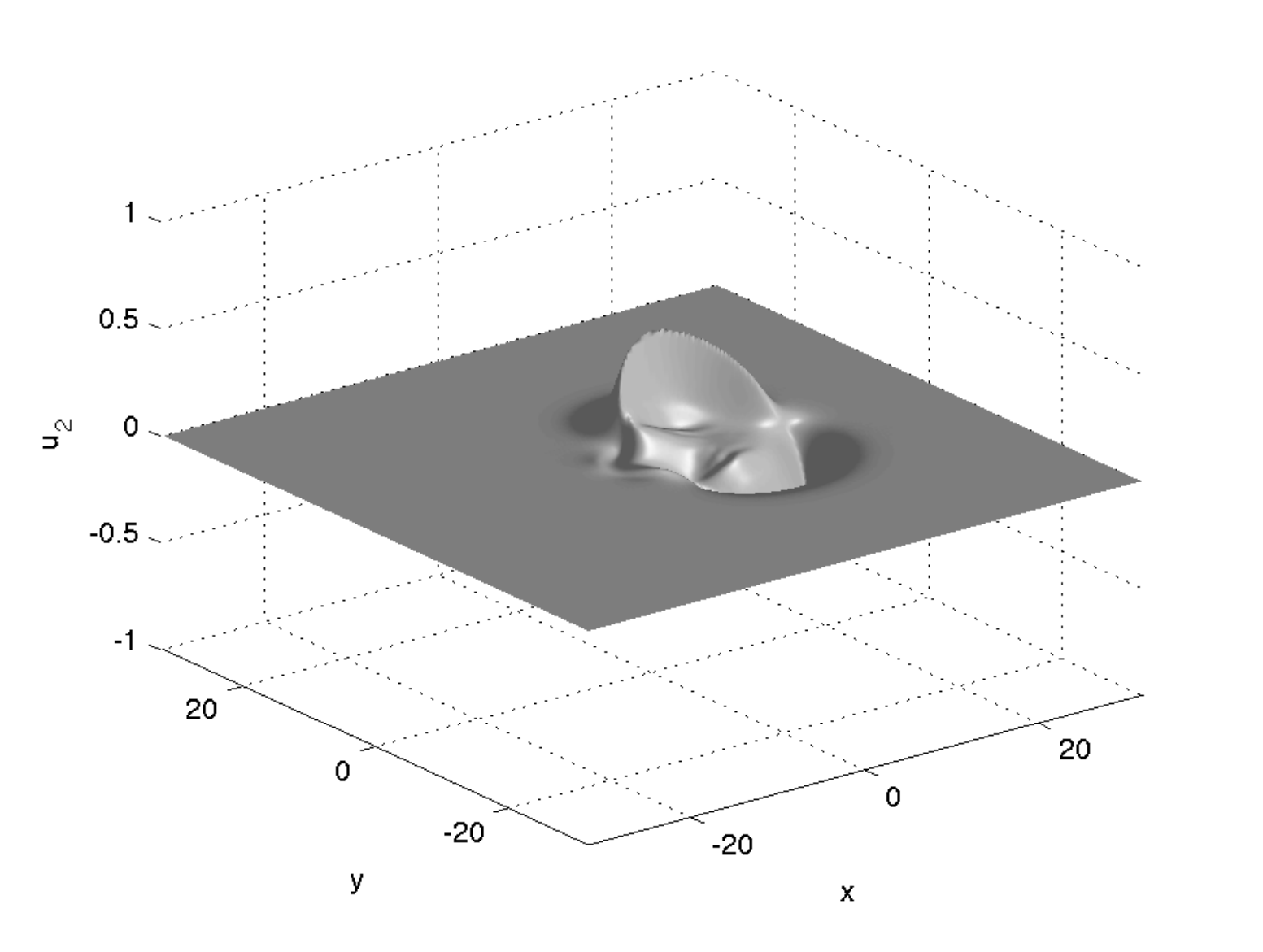}
\end{center}
\caption{Same as Figure \ref{fig:inital-data-smooth} at $t=20$, by using the operator-splitting pseudospectral method.  $ b=3/2$ for the simulations. $d x= d y=0.25$ and $\Delta t=0.025$. The computational domain is $[-32, 32]\times [-32, 32]$. (a) The first component  $u_1$. (b) The second component $u_2$.}
\label{fig:t20_pseudo}
\end{figure}

We propose an operator-splitting pseudospsectral method for solving the equations~(\ref{eq:EPDIFF}) by alternating between solving equations~(\ref{eq:2d_moment}) and (\ref{eq:2d_elliptic}).
In detail, the resulting algorithm consists of the following steps:
%
%\vskip 0.2in
%
%\noindent
%{\bf  Algorithm:}  
%
\vskip 12pt

\noindent
{\bf Step 1.}  Given smooth initial data $u_1^{0}$ and  $u_2^{0}$, we compute $m_1^{0}$ and $m_2^{0}$  by using equation (\ref{eq:m-Fourier}).

\vskip 12pt

\noindent
{\bf Step 2.} Integrate equation (\ref{eq:2d_moment}) by using the two-stage, second-order Runge-Kutta method (\ref{eq:2ndRK}). All derivatives are computed by the pseudospectral method. For example, the $j^{th}$ row of the partial derivative of $u_1m_1$ with respect to $x$ is computed by 
\beq
\left[(u_1m_1)_x\right]_{j} = \mathcal{F}^{-1}\left\{ik_{x}\mathcal{F}\left\{[u_1m_1]_{j}\right\} \right\},
\eeq
where $\mathcal{F}$ is the one-dimensional Fast Fourier Transform, $\mathcal{F}^{-1}$ is the inverse Fast Fourier Transform, $k_x$ is the corresponding wavenumbers, and $i=\sqrt{-1}$. 

\vskip 12pt 

\noindent
{\bf Step 3.} After the integration over $\Delta t$ in {\bf Step 2}, obtain $m_1^{1}$ and $m_2^{1}$. Compute $u_1^{1}$ and $u_2^{1}$ by using equation (\ref{eq:m-Fourier}) again.

\vskip 12pt 

\noindent
{\bf Step 4.}
Return to {\bf Step 2} and  {\bf Step 3} for computing  $m_1^{n}$ , $m_2^{n}$, $u_1^{n}$, and $u_2^{n}$, where $n=2, 3, \dots$.

\vskip 12pt 

\noindent
We remark that the proposed algorithm is a two-dimensional extension of the operator-splitting algorithms developed for the SW equation \cite{bib:DRP, bib:DRP2}. Similar to those one-dimensional solvers, an implicit iteration between equation (\ref{eq:2d_moment}) and (\ref{eq:2d_elliptic}) can be implemented for the current algorithm to guarantee the convergence of numerical solutions. 

We now are now in position to integrate the equations~(\ref{eq:EPDIFF}) for smooth initial condition. Consider the smooth initial data $\bs{u}=(u_1, u_2)^{T}$, where  

\begin{equation}\label{eq:smooth_IC}
\begin{split}
u_1(x,y) & = \sech\left(\ds\frac{x^2+y^2}{4}\right), \\
u_2(x,y) & = 0. 
\end{split}
\end{equation}
Figure \ref{fig:t3_pseudo}(a) and (b) are simulations for $u_1$ and $u_2$ at $t=3$, respectively. The operator-splitting  pseudospectral method is used with $d x=d y=0.125$ and $\Delta t=0.0125$.  The computational domain is $[-32, 32]\times [-32, 32]$. $ b=3/2$ for the simulations. Figure \ref{fig:t3_particle}(a) and (b) are the same simulations as Figure \ref{fig:t3_pseudo}(a) and (b), but are obtained by using the $N$-particle algorithm with $51\times 51$ particles placed in a $16\times 16$ domain initially. The time step is $\Delta t = 0.1$ for the particle algorithm. The velocities are reconstructed on a $32\times 32$ domain with $dx=dy=0.16$. Figure \ref{fig:t3_pseudo} and Figure \ref{fig:t3_particle} are virtually indistinguishable.  The difference between Figure \ref{fig:t3_pseudo} and Figure \ref{fig:t3_particle}, for both components in the maximum norm, is at the order of $O(10^{-3})$. We remark that for all simulations in this section we use the $N$-particle algorithm defined in equation (\ref{eq:N-particle}), for which the Green function is not normalized for convenience of comparison. 

Figures \ref{fig:t5_pseudo} and \ref{fig:t5_particle} depict the result of simulations at $t=5$ for the same initial data (\ref{eq:smooth_IC}), obtained by the pseudospectral and $N$-particle algorithms, respectively. The number of particles is $81\times81$, the same as that for the simulations at $t=3$.  At first glance, Figures \ref{fig:t5_pseudo} and \ref{fig:t5_particle} seem identical. However, if we blow up the region around the peak of $u_1$ in  Figure \ref{fig:t5_particle}(a), we can detect a saw-tooth-like roughness. This is because many particles have moved away from this region at this time, and the smooth wave cannot be represented by too few ``conon"-particles. When we increase the particle number from $81\times 81$ to $101\times 101$, this saw-tooth roughness becomes less visible, as shown in Figure \ref{fig:t5_particle_finer}. However, with the increment of particles from 81 to 101 in one direction, the computational cost increases about six-fold. This is because the cost of double summation is $O(N^2)$ (with $N$ the number of particles), and we reduce the temporal step size to a half to maintain the stability of the two-stage RK scheme. Together with overhead, the elapsed CPU-time for the $101\times 101$ mesh grid is 6 times more than the $81\times 81$ one. Therefore, long-time simulations for this example by using the $N$-particle algorithm are not feasible without introducing a fast algorithm, such as the fast multipole method, or by taking 
advantage of massive parallelization.  This is beyond the current scopes of the paper. Nevertheless, developing fast algorithms for the $N$-particle method can be 
implemented and we expect to report on this in the near future.  
Before ending this section, we demonstrate the ability of handling smooth data for the operator-splitting pseudospectral method. We evolve the initial data (\ref{eq:smooth_IC}) until $t=20$. Figure \ref{fig:t20_pseudo}(a) and (b) show the result of these algorithm simulations for $u_1$ and $u_2$, respectively.  

\section{Discussion and concluding remarks}

In this paper, we have studied a class of multidimensional PDEs for a parametric family of 
elliptic operators, and extended the class to a dispersive deformation which, to the best of our knowledge, has not been investigated in the literature. We have used the Lagrangian formulation as the most natural avenue for deriving finite-dimensional particle systems  discretizing the PDE system. These particle systems for non-smooth kernels (Green functions) of the invertible elliptic operator govern nontrivial dynamics worth examining in further detail. Within the class we have focussed on, the regularity of the Green functions is determined by the power of the elliptic operator, which we denoted by $b$. If $b=3/2$, the Green function has a finite jump in its radial derivative reminiscent of the peakon solution for the SW equation. In fact, by using this ``conon" case, in which the motion of the two-dimensional particles is restricted in a one-dimensional channel, we show that this choice of non-smooth kernel reduces the two-dimensional particle system to the {\it completely integrable} one-dimensional case, even though the ``single channel" solution retains its dependence on two spatial dimensions. With this reduction, complete intergrability persists for the dispersive deformation, giving rise to traveling wave solutions which are smooth along their direction of travel.

We have also studied particle collisions restricted to a line for various parameters $b$ of the Green function kernel~(\ref{eq:elliptic}). We have found that for sufficiently smooth kernels, when $b\ge3$, two particles head-on collisions in finite time are avoided. This is in contrast to their less-smooth counterparts with $b < 3$.   

A pseudospectral scheme for solving the PDEs under study was introduced to provide an independent means of numerically computing smooth solutions of the PDE's we studied. By comparing solutions obtained with this scheme with those from the $N$-particle algorithm, we show that the $N$-particle system can potentially be used as a Lagrangian method for solving the model PDEs, in particular when weak solutions are considered. Nevertheless, for long time simulations and smooth initial data, it is clearly necessary to develop fast summation algorithms for the $N$-particle  to achieve realistic computational costs.    

We do not investigate particle dynamics for $b\leq 1$ in this paper. Since the Green function of the elliptic operator corresponding to this power is unbounded at the support point, it would be necessary to regularize these kernels to implement an $N$-particle algorithm. The regularization results in a smooth modified kernel. In principle, the behavior of the regularized Green function should in principle be similar to that of this function with a parameter $b$ in the range $b\ge 3$. We leave this to future work. Also left to forthcoming investigations are the implications of the sensitivity to the singularity of the Green functions for image matching applications, and the convergence of the point particle approximation to solutions of the model PDEs under various singular kernels.

Finally, we remark that the flexibility in the choice of elliptic operators connecting the ``primary" field $\bs{u}$ and the ``auxiliary" field $\bs{m}$ could be exploited to move beyond the realm of interesting mathematical PDE's and towards more physically grounded models such as the Euler equations for ideal fluids \cite{MM13, CHJM14}. Doing so could provide valid alternatives to vortex methods for numerical simulations of 2D and 3D  Euler equations, as well as analytical tools which may prove useful in theoretical investigations of these equations.  
 
\section{Acknowledgements}

RC acknowledges the support of NSF DMS-0509423, CMG-0620687, DMS-0908423, DMS-1009750, and RTG DMS-0943851.

\appendix

\section{Lagrangian and Eulerian formulations}\label{sec:derivation}

Suppose that $J(\bs{\xi,t})$ is the Jacobian determinant of the diffeomorphism $\bs{x}=\bs{q}(\bs{\xi,t})$ parametrized by time $t$, 
$$
J(\bs{\xi,t}) \equiv \det \left( {\partial x_i \over \partial \xi_j }\right). 
$$
The conjugate field $\bs{p}(\bs{\xi}, t)$ is defined by 
\beq\label{eq:mp}
\bs{m}\big(\bs{q}(\bs{\xi,t}),t\big) \equiv \,{\bs{p}(\bs{\xi},t)\over  J(\bs{\xi,t}) } .
\eeq
The Yukawa operator $\mathcal{L}$ gives
\[
\bs{m}=\mathcal{L}\bs{u} \,\,\, ,\text{or}\,\,\, \bs{u} =\bs{G}*\bs{m}.
\]
Let $\bs y= \bs q(\bs \eta, t)$. We have 

\begin{equation}\label{eq:qdot}
\begin{split}
\bs{u} = {d \bs{q} \over d t} =& G_{ b-n/2}*\bs{m} \\
=& \int_{\mathbb{R}^{n}} 
G_{ b-n/2}\big(|\bs{x} -\bs{y})|\big) 
\, \bs{m} (\bs{y},t) \, d \bs{y} \,\\
=& \int_{\mathbb{R}^{n}} 
G_{ b-n/2}\big(|\bs{q}(\bs{\xi},t) -\bs{q}(\bs{\eta},t)|\big) 
\, \bs{m} (\bs{q}(\bs{\eta},t),t) \, d \bs{q}(\bs{\eta},t) \,\\
=& \int_{\mathbb{R}^{n}} 
G_{ b-n/2}\big(|\bs{q}(\bs{\xi},t) -\bs{q}(\bs{\eta},t)|\big) 
\, \bs{p} (\bs{\eta},t) \, d V_\eta.
\end{split}
\end{equation} 
We now show that if we define 
\[
{d \bs{p} \over d t} = - \int_{\mathbb{R}^{n}} 
G'_{ b-n/2}\big(|\bs{q}(\bs{\xi},t) -\bs{q}(\bs{\eta},t)|\big) 
{\bs{q}(\bs{\xi},t) -\bs{q}(\bs{\eta},t) \over |\bs{q}(\bs{\xi},t) -\bs{q}(\bs{\eta},t)|}
\, \, \bs{p}(\bs{\xi},t) \cdot \bs{p} (\bs{\eta},t) \, d V_\eta \, ,\\
\]
we recover the model PDEs (\ref{eq:EPDIFF}). To this end, we follow the diffeomorphism variable transformation to compute
\begin{equation}\label{dpdtRight}  
\begin{split}
{d \bs{p} \over d t} =& - \int_{\mathbb{R}^{n}} 
G'_{ b-n/2}\big(|\bs{q}(\bs{\xi},t) -\bs{q}(\bs{\eta},t)|\big) 
{\bs{q}(\bs{\xi},t) -\bs{q}(\bs{\eta},t) \over |\bs{q}(\bs{\xi},t) -\bs{q}(\bs{\eta},t)|}
\, \, \bs{p}(\bs{\xi},t) \cdot \bs{p} (\bs{\eta},t) \, d V_\eta \\
=& - \int_{\mathbb{R}^{n}} 
G'_{ b-n/2}\big(|\bs{x} -\bs{y}|\big) 
{\bs{x} -\bs{y} \over |\bs{x} -\bs{y}|}
  J (\bs{\xi}, t) \, \, \bs{m}(\bs{x},t) \cdot \bs{m} (\bs{y},t) \, d \bs{y} \\
=& - \bs{m}_j(\bs {x},t)\nabla_{\bs{x}}\int_{\mathbb{R}^{n}} 
G_{ b-n/2}\big(|\bs{x} -\bs{y}|\big)
   \, \, \bs{m}_j (\bs{y},t) \, d \bs{y} \, \, J (\bs{\xi}, t) \\
=& - \bs{m}_j(\bs {x},t)\nabla_{\bs{x}} \bs{u}_j J (\bs{\xi}, t) \\
=&-\bs{m}\cdot (\nabla \bs{u})^T  J (\bs{\xi}, t).
\end{split}
\end{equation} 
On the other hand, by the definition of the conjugate field (\ref{eq:mp}), we have  
\begin{equation}\label{dpdtLeft}  
\begin{split}
\frac{d\bs{p}}{d t} =&\frac{d}{dt} \left(\bs{m}(\bs{q}(\bs{\xi}), t)J(\bs{\xi}, t) \right)\\
=& \frac{d}{dt}\bs{m}(\bs{q}(\bs{\xi}), t)\, J(\bs{\xi}, t) + \bs{m}(\bs{q}(\bs{\xi}), t)\frac{d J(\bs{\xi}, t)}{dt}\\
=& \left( \bs{m}_t + \frac{d\bs{q(\bs{\xi}, t)}}{dt}\cdot\nabla\bs{m}(\bs{x}, t)\right) J(\bs{\xi}, t)  + \bs{m}(\nabla \cdot\bs{u})\big)J(\bs{\xi}, t) \\
=&\big(\bs{m}_t+(\bs{u}\cdot \nabla) \bs{m} +\bs{m}(\nabla \cdot\bs{u})\big)J(\bs{\xi}, t).                 
\end{split}
 \end{equation} 
Here we use the well known property of determinant differentiation (\ref{eq:det_t}). From equations (\ref{dpdtRight}) and (\ref{dpdtLeft}) the Eulerian 
form of the model equations (\ref{eq:EPDIFF}) follows.

Next, we derive the evolution equation of the determinant $J(\bs{\xi},t)$ in equation  (\ref{eq:jchars}). From the determinant differentiation (\ref{eq:det_t}) and equation (\ref{eq:qdot}), we have 
\beq
\begin{split}
\frac{d J}{dt} &= J\,\nabla\cdot{\bs u} = J\,\nabla\cdot  \int_{\mathbb{R}^{2}}G_{ b-1}\big(|\bs{x} -\bs{q}(\bs{\eta}, t)|\big) 
\, \bs{p} (\bs{\eta},t) \, d V_{\eta}\\
& = J\, \int_{\mathbb{R}^{2}}G'_{ b-1}\big(|\bs{x} -\bs{q}(\bs{\eta}, t)|\big) \frac{\big(\bs{x} -\bs{q}(\bs{\eta}, t)\big)\cdot\bs{p}(\bs{\eta}, t)}{|\bs{x} -\bs{q}(\bs{\eta}, t)|}\, d V_{\eta}.
\end{split}
\eeq
If $J$ is evaluated at $\bs{x} = \bs{q}(\bs{\xi}, t)$, then 
\beq
\frac{d J}{dt}  = J\, \int_{\mathbb{R}^{2}}G'_{ b-1}\big(|\bs{q}(\bs{\xi}, t) -\bs{q}(\bs{\eta}, t)|\big) \frac{\big(\bs{q}(\bs{\xi}, t) -\bs{q}(\bs{\eta}, t)\big)\cdot\bs{p}(\bs{\eta}, t)}{|\bs{q}(\bs{\xi}, t) -\bs{q}(\bs{\eta}, t)|}\, d V_{\eta}.
\eeq

\section{Notations for numerical implementation}\label{sec:notations}

For numerical implementation, we represent the $N$-particle systems in the following matrix-vector forms. Let 
\beq
\bs{Q}^\alpha = \left(\begin{array}{cccc}
q^\alpha_1\\ q^\alpha_2\\ \vdots\\ q^\alpha_N
\end{array} \right),
\quad\bs{P}^\alpha = \left(\begin{array}{cccc}
p^\alpha_1\\ p^\alpha_2\\ \vdots\\ p^\alpha_N
\end{array} \right),
\eeq
The system of equations for $\bs{Q}$ can be written as
\beq\label{eq:Q}
\frac{d\bs{Q}^\alpha}{dt} = \bs{A}\bs{P}^\alpha,\quad \alpha=1,2,
\eeq
where 
\beq
 \bs{A} = \left( \begin{array}{cccc}
G_{ b-1}(|\bs{q}_1-\bs{q}_1|) & G_{ b-1}(|\bs{q}_1-\bs{q}_2|) &\cdots & G_{ b-1}(|\bs{q}_1-\bs{q}_N|)\\
G_{ b-1}(|\bs{q}_2-\bs{q}_1|) & G_{ b-1}(|\bs{q}_2-\bs{q}_2|) &\cdots & G_{ b-1}(|\bs{q}_2-\bs{q}_N|) \\
\vdots & \vdots & \ddots &\vdots\\
G_{ b-1}(|\bs{q}_N-\bs{q}_1|) & G_{ b-1}(|\bs{q}_N-\bs{q}_2|) &\cdots & G_{ b-1}(|\bs{q}_N-\bs{q}_N|)
\end{array} \right),
\eeq
while the system for $\bs{P}$ is 
\beq\label{eq:P}
\frac{d\bs{P}^{\alpha}}{dt} = -\left[\bs{I}^{P^1}\bs{B}^{\alpha}\bs{P}^{1} +\bs{I}^{P^2}\bs{B}^{\alpha}\bs{P}^{2}\right],\quad \alpha=1, 2,
\eeq
with
\beq
\bs{B}^{\alpha} = \left( \begin{array}{cccc}
0& G_{ b-1}'(|\bs{q}_1-\bs{q}_2|)\frac{q^{\alpha}_1-q^{\alpha}_2}{|\bs{q_1-q_2}|} &\cdots & G_{ b-1}'(|\bs{q}_1-\bs{q}_N|)\frac{q^{\alpha}_1-q^{\alpha}_N}{|\bs{q_1-q_N}|}\\
G_{ b-1}'(|\bs{q}_2-\bs{q}_1|)\frac{q^{\alpha}_2-q^{\alpha}_1}{|\bs{q_2-q_1}|} & 0 &\cdots & G_{ b-1}'(|\bs{q}_2-\bs{q}_N|)\frac{q^{\alpha}_2-q^{\alpha}_N}{|\bs{q_2-q_N}|} \\
\vdots & \vdots & \ddots &\vdots\\
G_{ b-1}'(|\bs{q}_N-\bs{q}_1|)\frac{q^{\alpha}_N-q^{\alpha}_1}{|\bs{q_N-q_1}|}& G_{ b-1}'(|\bs{q}_N-\bs{q}_2|)\frac{q^{\alpha}_N-q^{\alpha}_2}{|\bs{q_N-q_2}|} &\cdots & 0
\end{array} \right),
\eeq
and
\beq
\bs{I}^{P^1} = \left( \begin{array}{cccc}
p_1^1& 0 &\cdots &0\\
0& p_2^1 &\cdots & 0 \\
\vdots & \vdots & \ddots &\vdots\\
0 & 0&\cdots & p_N^1
\end{array} \right), \quad
\bs{I}^{P^2} = \left( \begin{array}{cccc}
p_1^2& 0 &\cdots &0\\
0& p_2^2 &\cdots & 0 \\
\vdots & \vdots & \ddots &\vdots\\
0 & 0&\cdots & p_N^2
\end{array} \right).
\eeq

\section{Lyapunov function and stable manifold}\label{app:L}

For $ b\ge 3$, the Lyapunov (energy) function for the system (\ref{eq:q-z-1}) \& (\ref{eq:q-z-2}) that satisfies 
\beq
\frac{dV}{dt}(q,z) =0,
\eeq
is
\beq
V(q, z) = 2Hz^2+\frac{1}{2}G_{ b-1}(q).
\eeq
The solution through the point $(\eta_1, \eta_2)$ is given by the curve $V(q, z)=V(\eta_1, \eta_2)$. The real curves $V(q, z)\equiv h$, where $h$ is some constant, are given by equations
\beq
z=\pm \sqrt{\frac{2h-G_{ b-1}(q)}{4H}},
\eeq
for all $q$ for which $2h-G_{ b-1}(q) \ge 0$. This implies that when $h=V_0=V(0, 0) = \ds\frac{1}{2}G_{ b-1}(0)$, the stable manifold is 
\beq
z= \sqrt{\frac{G_{ b-1}(0)-G_{ b-1}(q)}{4H}}.
\eeq
Here we simply recover equation (\ref{eq:stable-1}).  Moreover, for $P=0$, equation (\ref{eq:H_collision}) becomes
\beq\label{eq:H-P=0}
H(q, z)=\frac{1}{4z^2}\left(G_{ b-1}(0) - G_{ b-1}(q)\right),
\eeq
for any $q$ and $z$. Substituting the above $H$ into the Lyapunov function, we obtain
\beq\label{eq:general-V}
V(q,z) = \frac{1}{2} G_{ b-1}(0) =V_0,
\eeq
for any $q$ and $z$.

\section{Exact solution of the head-on collision}\label{sec:exact} 

The Hamiltonian that generates the system (\ref{eq:ODE-nu-3/2}) is 
\beq\label{eq:Hamiltonian-head-on}
H_A=\frac{1}{2}\left(p_1^2+p_2^2\right)+p_1p_2e^{-|q_1-q_2|} =\frac{1}{2}\left(c^2+c^2\right)
\eeq
From equation (\ref{eq:one-particle-scaled}), we have $P(0)=c_1+c_2$ and $p(0)=c_1-c_2$. Also, because initially the locations of the two particles are $[x_1, 0]$ and $[x_2, 0]$, we have $Q(0) = x_1+x_2$ and $q(0) =  x_1-x_2=\gamma$ is the initial distance between the two particles. Therefore, with these initial data, the initial value problems (\ref{eq:ODE-nu-3/2}) can be solved. In particular, the second pair of equation (\ref{eq:ODE-nu-3/2}) can be solved by eliminating $p$ in the $\dot{q}$ equation by noting that  the Hamiltonian that generates  equation(\ref{eq:ODE-nu-3/2}) is
\beq\label{eq:Hamiltonian-head-on2}
H= \frac{1}{2}P^2\left(1+e^{-|q|}\right)+\frac{1}{2}p^2\left(1-e^{-|q|}\right) = c_1^2+c_2^2,
\eeq
and $P(t)= c_1+c_2$. The solution of the second pair of equation (\ref{eq:ODE-nu-3/2}) is shown in \cite{chh94} and is equal to
\beq\label{eq:sol1}
\begin{split}
q &= -\log\left[\frac{4\gamma(c_1-c_2)^2 e^{(c_1-c_2)t}}{\left(\gamma e^{(c_1-c_2)t}+4c_1^2\right)\left(\gamma e^{(c_1-c_2)t}+4c_2^2\right)} \right],\\
p &=\pm\frac{\gamma(c_1-c_2)(e^{-(c_1-c_2)t} -4c_1c_2)}{\gamma e^{-(c_1-c_2)t}+4c_1c_2}.
\end{split}
\eeq
The solutions of (\ref{eq:sol1}) for head-on collision (particle-antiparticle collision) has $c_1=-c_2=c$ and thus can be simplified to
\beq\label{eq:sol2}
\begin{split}
q &= -2\log\left[\frac{4c\sqrt{\gamma}e^{ct}}{\gamma e^{2ct}+4c^2} \right],\\
p &=\pm2c\frac{\gamma e^{-2ct} +4c^2}{\gamma e^{-2ct}-4c^2}.
\end{split}
\eeq 
If we choose $\gamma=4c^2$, the particle-antiparticle collision occurs at time $t=0$ at $x=0$, and 
\beq\label{eq:sol3}
\begin{split}
q &= -2\log \sech(ct),\\
p &=\pm\frac{2c}{\tanh(ct)}.
\end{split}
\eeq 

The constructed exact solution can be compared with numerical solution of equations (\ref{eq:q-z-1}) and (\ref{eq:q-z-2}).
For $ b=3/2$, the Green function is $G_{1/2}(q) = e^{-|q|}$. Note that the radial derivative has a finite jump and thus we write equations (\ref{eq:q-z-1}) and (\ref{eq:q-z-2}) as
\beq\label{eq:q-z-ex}
\begin{split}
\dot{q} & = 4Hz,\\
\dot{z} & =\frac{1}{2}\sgn(q) e^{-|q|},
\end{split}
\eeq
To compare the exact solution (\ref{eq:sol3}) with that obtained by solving equation (\ref{eq:q-z-ex}), we shift the collision time to $t_c > 0$  so that 
\beq\label{eq:sol4}
\begin{split}
q(t-t_c) &= -2\log \sech(c(t-t_c)),\\
p(t-t_c) &=\pm\frac{2c}{\tanh(c(t-t_c))}.
\end{split}
\eeq 
If we choose $c=2$, the initial separation of the particles at $t=0$ is $4c^2=16$, and the collision time $t_c$ satisfies
\beq\label{eq:tc}
\sech(2t_c) = e^{-8},\quad \text{or} \quad t_c=\frac{1}{2}\sech^{-1}(e^{-8}) \approx 4.346573576213.
\eeq 
For this choice of $c$ and $t_c$, the initial relative momentum is 
\beq\label{eq:p0}
p(-t_c) = \frac{4}{\tanh(-2t_c)}\approx -4.
\eeq
Hence the initial data for equation (\ref{eq:q-z-ex}) are
\beq\label{eq:init-comp}
q_0 = 16,\quad z_0 = -\frac{1}{4}.
\eeq

\begin{figure}[tbh]
\begin{center}
(a)\includegraphics[width=2.75in]{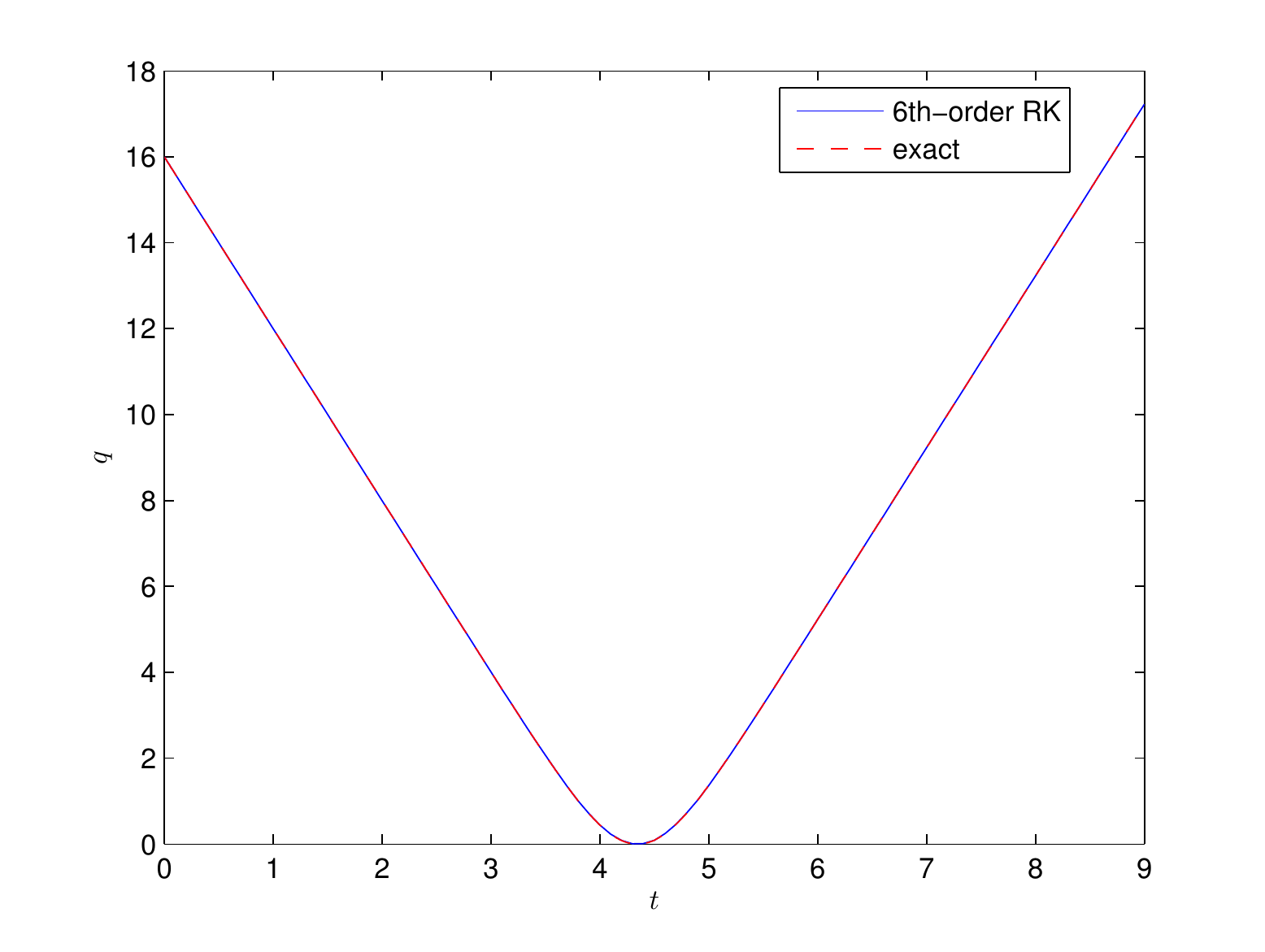}
(b)\includegraphics[width=2.75in]{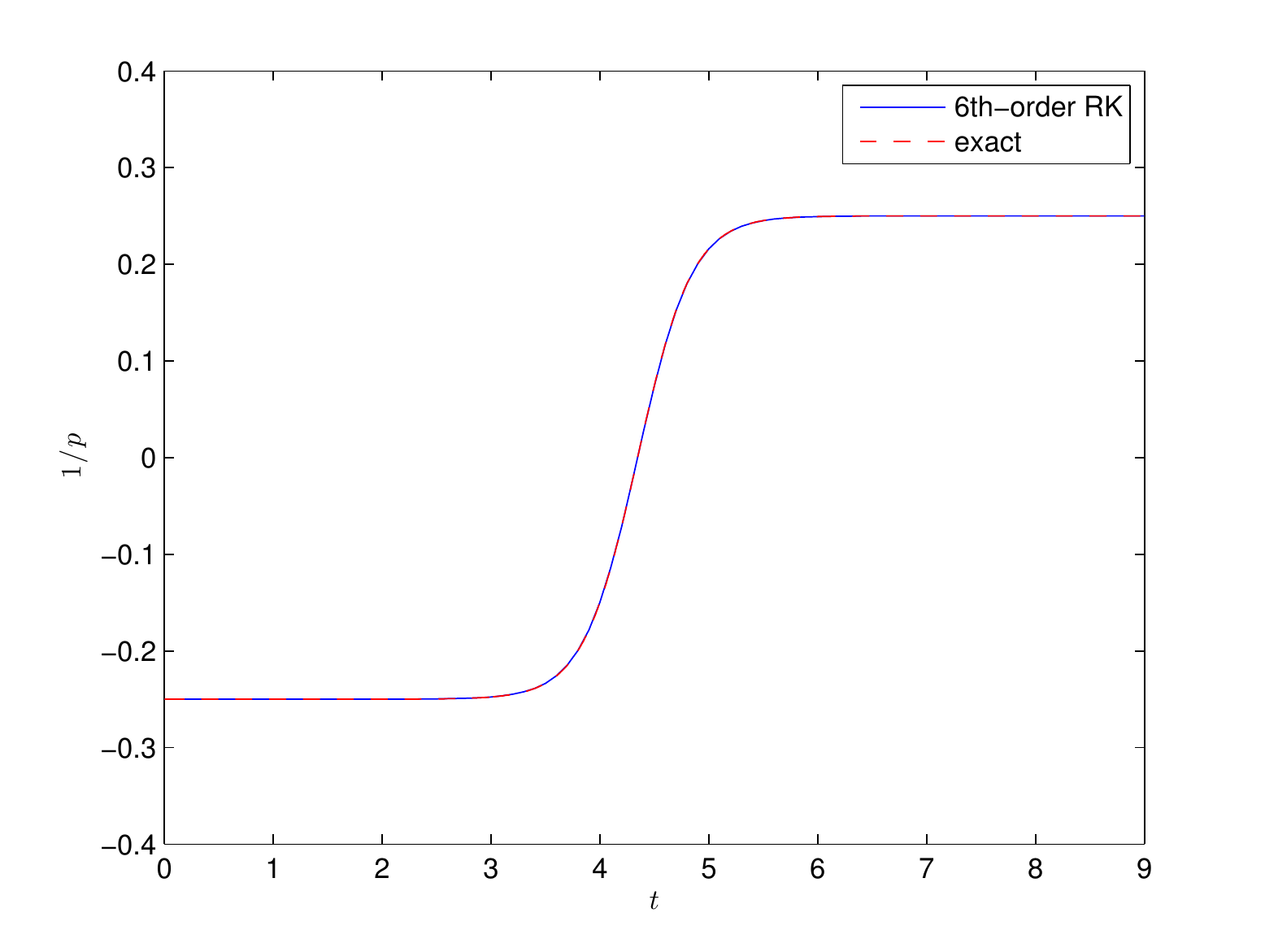}
\end{center}
\caption{Comparison between the exact solutions computed by equation (\ref{eq:sol4}) and the numerical solutions of equation (\ref{eq:q-z-ex}) obtained by using the sixth-order Runge-Kutta method with $q_0 = 16$ \& $z_0=1/4$. The time step is $\Delta t= 3.125$e-5. (a) The distance between the particle and the antiparticle versus time. The difference in 2-norm between the exact solutions and the numerical solutions is $2.9506$e-05. (b) The (inverse) relative momentum versus time. The difference in 2-norm between the exact solutions and the numerical solutions is $8.7397$e-07.}
\label{fig:comp_exact_RK6}
\end{figure}
Figure \ref{fig:comp_exact_RK6} compares the exact solutions (\ref{eq:sol4}) and the numerical solutions of equation (\ref{eq:q-z-ex}) obtained by using the sixth-order Runge-Kutta method. Equation (\ref{eq:init-comp}) is used as the initial data and the time step is $\Delta t=3.125$e-05. The differences between the two calculations in 2-norm, defined in (\ref{eq:2norm}), for $q$ and $z$ are $2.9506$e-05 and  $8.7397$e-07, respectively.


\begin{thebibliography}{10}

\bibitem{bib:as65}
\textrm{Abramowitz,~M. and Stegun,~I.A.}.
\textit{Handbook of Mathematical Functions}.
\textrm{Dover, New York, 1965.}

\bibitem{bib:bmty05}
\textrm{Beg~M. F., Miller~M. I., Trouv\'e,~A., and Younes,~L.}.
\textrm{Computing large deformation metric mappings via geodesics flows of diffeomorphisms}.
\textit{Int. J. Comp. Vis.}, {\bf 61}(2), 139-157, 2005.

\bibitem{bib:butcher64}
\textrm{Butcher,~J. C.}.
\textrm{On Runge-Kutta processes of high order}.
\textit{J. Austral. Math. Soc.},{\bf 4}(6), 179-194, 1964.

\bibitem{bib:DCDS03}
\textrm{Camassa~R.},
\textrm{Characteristics and initial value problem of a completely integrable shallow water equation}.
\textit{DCDS-B}
\textrm{\bf 3}, 115-139, 2003.

\bibitem{bib:ch93}
\textrm{Camassa~R. and Holm,~D. D.}.
\textrm{An integrable shallow water equation with peaked solitons}.
\textit{Phys. Rev. Lett.}, {\bf 71}, 1961-1964, 1993.

\bibitem{chh94}
\textrm{Camassa~R., Holm~D. D., and Hyman,~J. M.}.
\textrm{A new integrable shallow water equation}.
\textit{Advan.  Appl. Mech.}, {\bf 31},1-33, 1994.

\bibitem{bib:Camassa05}
\textrm{Camassa~R., Huang~J., and Lee~L.}.
\textrm{On a completely integral numerical scheme for a nonlinear shallow-water wave equation}.
\textit{J. Nonlin. Math. Phys.}
\textrm{{\bf 12}, 146-162, 2005}.

\bibitem{bib:chl06}
\textrm{Camassa,~R., Huang,~J., and Lee,~L.}.
\textrm{Integral and integrable algorithms for a nonlinear shallow-water wave equation}.
\textit{J. Comput. Phys.},
\textrm{{\bf 216}, 547-572, 2006}.

\bibitem{bib:Camassa07}
\textrm{Camassa~R. and Lee~L.}.
\textrm{A completely integrable particle method for a nonlinear shallow-water wave equation in periodic domains}.
\textit{DCDIS-A},
\textrm{{\bf 14(S2)}, 1-5, 2007}.

\bibitem{bib:cl08}
\textrm{R. Camassa and L. Lee}.
\textrm{Complete integrable particle methods and the recurrence of initial states for a nonlinear
shallow-water wave equation}.
\textit{J. Comp. Phys.},
\textrm{{\bf 227}, 7206-7221, 2008}.

\bibitem{bib:cdm12}
{Chertock,~A., Du Toit,~P., and Marsden,~J. E.}.
\newblock{Integration of the EPDIFF equation by particle methods},
\newblock{\em ESAIM:M2AN}, {\bf 46}, 515-534, 2012.

\bibitem{bib:DRP}
\textrm{Chiu,~P.H., Lee,~L., and Sheu,~T.W.H.}.
\textrm{A dispersion-relation-preserving algorithm for a nonlinear shallow-water wave equation}.
\textrm{ \em J. Comput. Phys.}, {\bf 228}, 8034-8052, 2009.

\bibitem{bib:DRP2}
\textrm{Chiu,~P.H., Lee,~L., and Sheu,~T.W.H.}.
\textrm{A sixth-order dual preserving scheme for the SW equation.}.
\textrm{\em J. Comput. Appl. Math.}, {\bf 223}, 2767-2778, 2010.

\bibitem{bib:JCP_RB}
\textrm{Camassa,~R.,  Chiu,~P.H. Lee,~L., and Sheu,~T.W.H.},
\textrm{Viscous and inviscid regularizations in a class of evolutionary partial differential equations},
\textit{ J. Comput. Phys.}, 
\textrm{{\bf 229}, 6676-6687, 2010.}

\bibitem{bib:particle_RB}
\textrm{Camassa,~R.,  Chiu,~P.H. Lee,~L., and Sheu,~T.W.H.},
\textrm{A particle method and numerical investigation of a quasi-linear partial differential equation},
\textit{Comm. Pure and Appl. Math.}, 
\textrm{{\bf 10}, 1503-1512, 2011}

\bibitem{bib:dgm98}
{Dupuis,~P., Grenander,~U., Miller,~M. I.}.
\textrm{Variational problems on flows of diffeomorphisms for image matching}.
\textit{Q. Appl. Math.}, {\bf 56}, 587-600, 1998.

\bibitem{bib:hmr98}
{Holm,~D. D., Marsden,~J. E., and Ratiu,~T. S.}.
\newblock{Euler-Poincar\'e models of ideal fluids with nonlinear dispersion}.
\newblock{Phys. Rev. Lett.}, {\bf 80}(19), 4173-4176, 1998.

%\bibitem{bib:hs03}
%\textrm{Holm,~D. D. and Staley,~M. F.}.
%\newblock{Wave structure and nonlinear balances in a family of evolutionary PDEs}.
%\newblock{\em SIAM J. Applied Dynamical System}, {\bf 2}(3), 323-380, 2003.
%
%\bibitem{bib:hs13}
%\textrm{Holm,~D. D. and Staley,~M. F.}.
%\newblock{Interaction dynamics of singular wave fronts}.
%\newblock{\em arXiv:1301.1460v1}, 2013. %S170-S178, 2004.

\bibitem{bib:hrty04}
\textrm{Holm,~D. D., Ratnanather,~J. T., Trouv\'e,~A., and Younes,~L.}.
\newblock{Soliton dynamics in computational anatomy}.
\newblock{\em Neuroimage}, {\bf 23}, S170-S178, 2004.

\bibitem{CHJM14}
{Cotter,~C. J., Holm,~D. D., Jacobs,~H. O., and Meier,~D. M.}.
\newblock{A jetlet hierarchy for ideal fluid dynamics}.
\newblock{\em J. Phys. A: Math}, {\bf 47}, 352001, 2014.

%\bibitem{bib:hjz09}
%\textrm{Huang,~J., Jia,~J., Zhang,~B.}.
%\newblock{FMM-Yukawa: An adaptive fast multipole method for screened Coulomb interactions}.
%\newblock{\em Compt. Phys. Comm.}, {\bf 180}, 2331-2338, 2009.

\bibitem{bib:jm00}
{Joshi,~S. and Miller,~M. I.}.
\textrm{Landmark matching via large deformation diffeomorphisms}.
\textit{IEEE Trans. Image Processing}, {\bf 9}, 1357-1370, 2000.

\bibitem{bib:kl14}
\textrm{Kuang,~D. and Lee,~L.}
\newblock{A feedback control geodesic landmark shooting algorithm for template matching and pattern recognition}.
\newblock{preprint}.

\bibitem{Landau}
{Landau,~L.D. and~Lifshitz,~E.M.}
\textit{Mechanics}.
\newblock{Pergamon, 2nd edition, Oxford, 1969.}

\bibitem{image1}
{Marsland~S. and Twining,~C.}.
\newblock{Constructing diffeomorphic representations for the groupwise 
analysis of non-rigid registrations of medical images}.
\newblock{\em  IEEE Transactions on Medical Imaging}, {\bf 23}(8), 1006-1020, 2004.

\bibitem{NP}
{McLachlan,~R.I. and Marsland~S.}.
\newblock{N-particle dynamics of the Euler equations for planar diffeomorphisms}.
\newblock{\em Dyn. Sys.}, {\bf 22}(3), 269-290, 2007.

\bibitem{bib:my01}
{Miller,~M.I. and Younes,~L.}.
\textrm{Group actions, homeomorphisms, and matching: A general framework}.
\textit{International Journal of Computer Vision}, {\bf 41}(1/2), 61-84, 2001.

\bibitem{anatomy}
{Miller,~M.I., Trouv\'e,~A., and Younes,~L.}.
\newblock{On metrics and the Euler-Lagrange equations of 
computational anatomy}.
\newblock{\em Annual Reviews in Biomedical Engineering}, {\bf 4}, 375-405, 2002.

\bibitem{bib:mty06}
{Miller,~M. I., Trouv\'e,~A., and Younes,~L.}.
\newblock{Geodesic shooting for computational anatomy}.
\newblock{\em J Math Imaging}, {\bf 24}, 209-228, 2006.

\bibitem{MZM06}
\textrm{Mohseni,~K., Zhao,~H., and Marsden,~J. }.
\textrm{Shock regularization for the Burgers equation}.
\textit{AIAA Paper 2006-1516, 44th AIAA Aerospace Science Meeting and Exhibit} \{Reno, Nevada, Jan, 9-12, 2006.\}

\bibitem{bib:mumford}
\textrm{Mumford,~D. and  Desolneux,~A.}.
\textit{Pattern Theory: The Stochastic Analysis of Real World Signals}.
\textrm{A K Peters, Lid, Natick, MA, 2010}.

\bibitem{MM13}
{Mumford,~D. and Michor,~P.}.
\newblock{On Euler's equation and ``EPDIFF''}.
\newblock{\em J. Geom. Mech.}, {\bf 6}(3), 319--344, 2013.

\bibitem{bib:perko}
{Perko,~L.}
\textit{Differential Equations and Dynamical Systems}.
\newblock{Springer-Verlag, 2nd edition, New York, 1996.}

\bibitem{bib:t95}
{Trouv\'e,~A.}.
\textrm{An infinite dimensional group approach for physics based model}.
\textrm{Technical report, 1995}.

\bibitem{image2}
{Vaillant~M., Miller,~M. I., Younes`L., and Trouv\'e~A.}.
\newblock{Statistics on diffeomorphisms via tangent space representations}.
\newblock{\em  NeuroImage}, {\bf 23},  S161-S169, 2004



\end{thebibliography}
\end{document}